\documentclass[11pt,reqno]{amsart}
\usepackage[a4paper,margin=26mm]{geometry}
\usepackage[T1]{fontenc}
\usepackage{lmodern,microtype}
\usepackage{amsmath,amssymb,amsthm,mathtools,mathrsfs,enumitem}
\usepackage{comment}
\usepackage{tikz}
\usepackage{float}
\usepackage{hyperref}
\usetikzlibrary{arrows.meta,calc}
\hypersetup{colorlinks=true,linkcolor=blue,citecolor=blue,urlcolor=blue,
 pdftitle={Macroscopic return-time spectra and Hausdorff dimension transitions},
 pdfauthor={Jinjun Li},pdfsubject={Return-time spectra and critical phase transitions}}
\allowdisplaybreaks[2]
\numberwithin{equation}{section}
\newtheorem{theorem}{Theorem}[section]
\newtheorem{lemma}[theorem]{Lemma}
\newtheorem{maintheorem}{Theorem}

\newtheorem{proposition}[theorem]{Proposition}
\newtheorem{corollary}[theorem]{Corollary}
\theoremstyle{definition}

\theoremstyle{remark}
\newtheorem{remark}[theorem]{Remark}
\newcommand{\N}{\mathbb N}\newcommand{\R}{\mathbb R}
\newcommand{\cH}{\mathcal H}
\newcommand{\dimH}{\dim_{\rm H}}
\newcommand{\lowD}{\mathsf D^-}\newcommand{\upD}{\mathsf D^+}
\newcommand{\DD}{\mathsf D}\newcommand{\NN}{\mathbb N}
\newcommand{\one}{\mathbf1}
\DeclareMathOperator{\card}{\#}

\title[Macroscopic return-time spectra]{Macroscopic return-time spectra and Hausdorff dimension transitions}

\author[J. J. Li]{Jinjun Li}
\address{School of Mathematics and Information Science, Guangzhou University, Guangzhou 510006, P. R. China}
\email{lijinjun@gzhu.edu.cn}

\date{October 2026}

\begin{document}
\begin{abstract}
Can the complete macroscopic spectrum of return times determine the
Hausdorff dimension of the corresponding set of initial points?
We classify all common lower and upper macroscopic
intermediate-dimension spectra of zero-density integer sets and
realize every admissible spectrum on a full-dimensional family of
zero-frequency orbits under logarithmic-power bounds on inverse-branch
contraction.  For first-return branches with stretched-exponential
contraction $e^{-n^\delta}$, $0<\delta<1$, complete spectrum levels
exhibit a sharp zero/full dimension transition at the counting
endpoint $\delta$. At the critical exponent, exact laws
$N_y(t)\sim c_0t^\delta(\log t)^\beta$ yield zero dimension for
$\beta\le0$ and full dimension for $\beta>0$, even with the
complete spectrum fixed.  An explicit two-branch intermittent
interval map realizes this transition through genuine first returns.
Above the critical exponent, imposing an exact counting law and a
complete spectrum does not reduce the dimension of an intrinsic
mean-energy level in an affine return tower.  For nonlinear branches,
we obtain quantitative bounds and an exact conditional-pressure formula
under thermodynamic regularity assumptions, verified for an explicit
non-affine family.
\end{abstract}
\maketitle
\setcounter{tocdepth}{2}
\begingroup
\small
\setlength{\baselineskip}{11pt}
\tableofcontents
\endgroup
\enlargethispage{2pt}

\section{Introduction}\label{sec:intro}
\subsection{Return frequency, temporal geometry, and phase-space dimension}
Return-time statistics describe the temporal organization of visits
to a target, while Hausdorff dimension measures the size of the
corresponding sets of initial points.  We ask to what extent a
complete macroscopic return-time spectrum determines this geometric
size.  In the zero-frequency regime, exact counting and mean-energy
constraints can distinguish levels with the same complete spectrum.
Let \(T:X\to X\) be a dynamical system and \(Y\subset X\) a fixed target set. For each \(y\in X\), define the set of positive visiting times and its associated counting function by
\[
 R(y,Y)=\{n\ge1:T^ny\in Y\},\qquad
 N_y(t)=\#(R(y,Y)\cap[1,t]).
\]
The \emph{return frequency} (or visit frequency) of $y$ to $Y$ is
$
 p(y)=\lim_{t\to\infty}\frac{N_y(t)}{t},
$
whenever the limit exists.  It measures the proportion of physical
time spent visiting $Y$, not the distribution of individual visits.

Classical recurrence theory provides fundamental information on
visiting times.
Poincar\'e's recurrence theorem \cite{Poincare1890} guarantees
infinitely many returns to a set of positive measure for almost
every point in that set under a probability-preserving transformation.
In an ergodic probability-preserving system, for a measurable set
$Y$ with $\mu(Y)>0$, Birkhoff's theorem \cite{Birkhoff} gives
$p(y)=\mu(Y)$ for $\mu$-almost every $y$, while Kac's formula
\cite{Kac} identifies the mean first-return time, with respect
to the normalized restriction of $\mu$ to $Y$, as $1/\mu(Y)$. Khintchine's recurrence theorem \cite{Khintchine} and
Furstenberg's multiple recurrence theorem
\cite{Furstenberg1977,Furstenberg} further reveal arithmetic
structure in the return times of sets of positive measure.
These classical results describe fundamental aspects of recurrence,
but do not by themselves determine the large-scale distribution
of visiting times along individual exceptional orbits or the
Hausdorff dimension of the corresponding level sets.

Zero visit frequency is typical for finite-measure targets in
conservative infinite-measure ergodic systems \cite{Aaronson}, and
can also occur along exceptional orbits of probability-preserving
systems.  Return-time statistics \cite{HLV,BSTV} and multifractal
results for returns to shrinking neighborhoods \cite{OlsenReturn}
address related but different questions: their observables are
return-time distributions or changing targets, whereas here the
target $Y$ is fixed.  An exact counting law $N_y(t)\sim c_0t^b$ does not determine
\emph{where} visits occur within distant time scales.  Even when
this law and the complete macroscopic spectrum are both fixed,
intrinsic energy constraints can distinguish refined levels with
different Hausdorff dimensions.  We distinguish these refined levels
from the complete spectrum level itself.

To study the large-scale geometry of return-time sets, we use
macroscopic intermediate dimensions.  The underlying theory goes
back to the macroscopic Hausdorff dimension of Barlow and Taylor
\cite{BT89,BT92}, based on coverings of distant shells; see also
Daw and Seuret \cite{DawSeuret}.  We use the lower and upper
macroscopic intermediate dimensions of Li and Xi \cite{LX}, denoted
by $\mathsf D^-_\theta$ and $\mathsf D^+_\theta$, for
$0\le\theta\le1$.  For integer sets, $\theta$ records clustering
across scales, and $\mathsf D^+_1$ is the upper counting exponent.

Li and Xi \cite{LX} develop the macroscopic theory, including
projection and capacity results, whereas Li, Xi and Zhang \cite{LXZ}
study discrete intermediate dimensions with a different aggregation
of shells.  In the bounded-set setting, Banaji and Rutar
\cite{BanajiRutar} classified attainable intermediate-dimension
functions using a different dimensional framework.  Here the
covering functional sums over the entire tail of distant shells:
we characterize its common lower and upper spectra for zero-density
integer sets, then realize these spectra by \emph{genuine first
returns} and determine the Hausdorff dimensions of the corresponding
orbit sets.  In particular, the dynamical dimension transitions do
not follow from the static attainable-spectrum classifications.
The precise definitions appear in Section~\ref{sec:prelim}.

Three questions guide the paper.  Which functions on $[0,1]$ occur
as common lower and upper spectra of zero-density integer sets?
Under what geometric assumptions can every such function be
realized on a full-dimensional family of orbits?  Once a spectrum
is fixed, how do exact counting and energy constraints affect the
dimensions of the resulting refined levels?

For a fixed system, the complete spectrum level is a fixed set,
and Theorem~\ref{main:transition}(i) determines its dimension.
The finer phenomenon concerns \emph{different complete refined
levels within the same spectrum level}.  For the interval map in
Proposition~\ref{prop:intervalmodel}, for example, every $c_0>0$
and $\beta>0$ give
\[
\begin{aligned}
 \text{$\dim_H\{y\in Y:N_y(t)\sim c_0t^\delta\}=0$ \quad and \quad $
 \dim_H\{y\in Y:N_y(t)\sim c_0t^\delta(\log t)^\beta\}=1.$}
\end{aligned}
\]
Both dimensions remain unchanged if one also requires the same
complete spectrum $\varphi_\delta$, defined in \eqref{yingshe}.
Thus the logarithmic transition is already a statement about exact
counting levels, and persists under the additional spectrum constraint.

The mechanism depends on the distinction between physical return
time and geometric contraction.  If $F=T^\tau$ is a first-return
map to an interval $Y$, a branch with roof $n$ takes $n$ iterates
before the next return, whereas its logarithmic geometric contraction is
$-\log|\psi_n'|$, where $\psi_n$ is the inverse return branch.
For intermittent maps with polynomial branch derivatives this
quantity grows logarithmically in $n$ \cite{LSV,Young99}; long waits can then
be inserted without dimension loss.  When
$-\log|\psi_n'|\asymp n^\delta$, the same insertion scheme can
have dimension zero.  We replace it by balanced fractal clusters
and bound the prescribed contraction at every prefix.  A
guarded-block lemma converts this bound into a dimension estimate
for arbitrary Euclidean balls.

The energy-level formulas are related to the multifractal theory
of Birkhoff averages and constrained pressure.  Olsen developed
frameworks for deformed Birkhoff averages and multifractal pressure
\cite{Olsen2003,Olsen2017}; related theories for countable Markov
systems were developed by Iommi, Jordan and others
\cite{IommiShift,IJ,IJQuotients,JT}.  The entropy maximization and
pressure expressions used here belong to this established
framework.  Our contribution is to preserve the dimension of an
entire energy level while simultaneously prescribing an exact
physical-time counting law and a complete macroscopic spectrum.
The nonlinear extension retains the full geometric potential,
and its thermodynamic assumptions are verified for an explicit
non-affine branch family.

\subsection{Main results}
Call a function \(h:[0,1]\to[0,1]\) scale-admissible if \(h\) is nondecreasing on \([0,1]\) and the function
$\theta\longmapsto\theta(1-h(\theta))$
is nondecreasing on \([0,1]\).
Such a function is continuous at every positive parameter, but may
jump at zero.  The condition is precisely the compatibility relation
forced by the scale comparison inequalities for macroscopic
intermediate dimensions; it is therefore the intrinsic restriction on
possible spectra rather than an auxiliary assumption.  The first theorem is purely a classification of integer
sets, with no assumption on an underlying dynamical system: the known
scale comparison is the only restriction, even when zero natural
density is imposed.
\begin{maintheorem}[Theorem~\ref{thm:attain}]
\label{main:admissibility}
There exists a zero-density integer set $E$ with
\[
 \mathsf D^-_\theta(E)=\mathsf D^+_\theta(E)=h(\theta)
 \quad(0\le\theta\le1)
\]
if and only if $h$ is scale-admissible.
\end{maintheorem}
Necessity is the comparison inequality of \cite[Proposition~2.6]{LX}.
Sufficiency is proved by supporting scale curves and a diagonal
construction with uniform estimates for the covering sums over
all distant shells.  It
includes both the possible jump at zero and the endpoint $h(1)=1$.
Theorem A gives the complete static classification of possible
macroscopic spectra.  The next theorem addresses the genuinely
 dynamical realization problem: whether such spectra can be produced
by first-return structures on large dimensional sets.

For the dynamical statements,  let \(K\) be the limit set generated by countably many inverse first-return branches, and let \(K_\infty\subseteq K\) denote the set of points whose orbits return to \(Y\) infinitely often. The common geometric
assumptions (G1) and (G2) in Section~\ref{sec:guards} require uniform
contraction and bounded distortion, an interior guard branch, and
finite-subsystem approximation of $\dim_H K$.  Separately, (R)
requires every positive integer to occur as a genuine first-return
length.  Theorems B and C impose different, mutually exclusive
contraction regimes: the logarithmic-power upper bound (L$_q$),
respectively the two-sided stretched-exponential bounds
\eqref{eq:twosidedstretched}.  Thus the hypotheses on return lengths
and those on branch geometry are not conflated.

Unlike the preceding realization problem for abstract integer sets, the return times here must arise from actual dynamical orbits, with successive returns governed by the first-return structure of the system.
\begin{maintheorem}[Theorem~\ref{thm:allprofiles}]
\label{main:realization}
Assume the structural conditions \textup{(G1), (G2)},
the genuine-return condition \textup{(R)}, and the logarithmic-power
contraction condition \textup{(L$_q$)} for some $q\ge1$.
Then every scale-admissible $h$ is realized on a full-dimensional
zero-frequency set:
\[
 \dim_H\left\{y\in K_\infty:p(y)=0,\quad
 \mathsf D^-_\theta(R(y,Y))=\mathsf D^+_\theta(R(y,Y))=h(\theta)
 \text{ for all }\theta\right\}=\dim_H K.
\]
\end{maintheorem}

Under stretched-exponential contraction, the endpoint exponent
instead determines a zero/full Hausdorff dimension transition.
At the critical counting exponent, an exact counting law detects a
second-order logarithmic transition that is invisible to the complete
macroscopic spectrum.

 For a scale-admissible function
$h:[0,1]\to[0,1]$, define
\[
\begin{aligned}
	\mathscr E(h)=\left\{y\in K_\infty:\;p(y)=0,\quad \mathsf D^-_\theta(R(y,Y))
	=\mathsf D^+_\theta(R(y,Y))
	=h(\theta)
	\quad\text{for all }\theta\in[0,1]\right\}.
\end{aligned}
\]
Fix \(0<\delta<1\), and let \(\varphi_\delta\) denote the member \(a=\delta\) of the one-parameter family \(\varphi_a\) defined in \eqref{yingshe} below.
For $c_0>0$ and $\beta\in\mathbb R$, define
\begin{equation}\label{eq:critical-log-level}
	\begin{aligned}
		\mathscr C_\beta(c_0)
		=\Bigl\{y\in K_\infty:\;
		N_y(t)&\sim c_0t^\delta(\log t)^\beta
		\quad(t\to\infty),\\
		&\mathsf D^-_\theta(R(y,Y))
		=\mathsf D^+_\theta(R(y,Y))
		=\varphi_\delta(\theta)
		\quad\text{for all }\theta\in[0,1]
		\Bigr\}.
	\end{aligned}
\end{equation}

The first transition depends only on the counting endpoint $h(1)$.
The second resolves exact logarithmic counting levels within the
critical spectrum level, and also determines their dimensions
without the spectrum constraint.

\begin{maintheorem}[
	Theorems~\ref{thm:generalphase}
	and~\ref{thm:logcritical}]
	\label{main:transition}
	Assume that $\mathrm{(G1)}$, \(\mathrm{(G2)}\), and \(\mathrm{(R)}\) hold, and that the two-sided stretched-exponential bounds \eqref{eq:twosidedstretched} hold for some \(0<\delta<1\).
	\begin{enumerate}[label=\textup{(\roman*)},
		leftmargin=*]
		\item For every scale-admissible $h$,
		\[
		\dim_H\mathscr E(h)=
		\begin{cases}
			0, & h(1)<\delta,\\
			\dim_H K, & h(1)\ge\delta.
		\end{cases}
		\]
		
		\item For every $c_0>0$ and $\beta\in\mathbb R$,
		the set $\mathscr C_\beta(c_0)$ is nonempty, and
		\[
		\dim_H\mathscr C_\beta(c_0)=
		\begin{cases}
			0, & \beta\le0,\\
			\dim_H K, & \beta>0.
		\end{cases}
		\]
Moreover,
		\[
		\begin{aligned}
			\dim_H\bigl\{y\in K_\infty:
			N_y(t)\sim c_0t^\delta(\log t)^\beta\bigr\} =
			\begin{cases}
				0, & \beta\le 0,\\[2pt]
				\dim_H K, & \beta>0.
			\end{cases}
		\end{aligned}
		\]
	\end{enumerate}
\end{maintheorem}
Part (i) is Theorem~\ref{thm:generalphase}; part (ii) is
Theorem~\ref{thm:logcritical}.  Exact power counts
$N_y(t)\sim c_0t^b$ have dimension zero for $b\le\delta$ and full
dimension for $b>\delta$ (Theorem~\ref{thm:exactpower}).
Thus a critical counting exponent, an exact leading coefficient, and an
arbitrarily small positive logarithmic correction impose different
geometric constraints.  Proposition~\ref{prop:intervalmodel} supplies
an explicit two-branch interval map with a neutral fixed point and
genuine return derivatives
$|\psi_n'|\asymp n^{\delta-1}e^{-n^\delta}$.
For this model all full-dimensional conclusions have dimension one.
The original map is $C^1$ at the neutral fixed point; the induced
branches have uniform contraction and distortion.

Above the critical exponent, fixing the mean return energy yields
intermediate dimensions even while the exact count and complete
spectrum stay fixed.  For the affine full-branch tower with a fixed parameter \(\kappa>0\) and contraction ratios
$r_n=Z(\kappa)^{-1}e^{-\kappa n^\delta}$, write
$Z(\lambda)=\sum_{n\ge1}e^{-\lambda n^\delta}$. For $y\in K_\infty$, let
$0=t_0(y)<t_1(y)<t_2(y)<\cdots$
be its successive return times to $Y$, and set
$\tau_i(y)=t_i(y)-t_{i-1}(y)$ for $i\ge1$.
These are the genuine first-return lengths associated
with the coding described in Section~\ref{sec:guards}.
The joint levels in the next theorem prescribe three quantities at
once: an exact physical-time count $N_y(t)\sim c_0t^b$, the mean
return energy $\lim_{k\to\infty}k^{-1}\sum_{i=1}^k\tau_i(y)^\delta=E$,
and the common macroscopic spectrum $h_{q,b}$ defined below.
\begin{maintheorem}[Theorem~\ref{thm:intrinsicpressure}]
\label{main:secondorder}
Fix $\delta<b<1$, $c_0>0$, and $q>(1-b)/(b-\delta)$, and put
\[
 d_{q,b}=\frac1{1+(1-b)(q+1)},\qquad h_{q,b}(0)=d_{q,b},
\]
\[
 h_{q,b}(\theta)=\max\left\{d_{q,b},
             \left(1-\frac{1-b}{\theta}\right)_+\right\}\quad(\theta>0).
\]
For each prescribed mean \(\delta\)-energy \(E>1\), the joint level set
$\mathscr L(c_0,b,E;q)$ defined in \eqref{eq:secondorder-levelset}
has Hausdorff dimension
\[
 \frac{\lambda(E)E+\log Z(\lambda(E))}
      {\kappa E+\log Z(\kappa)},
 \]
 where \(\lambda(E)>0\) is the unique solution of
$ -\frac{Z'(\lambda(E))}{Z(\lambda(E))}=E.$
This equals the dimension of the energy level without the two
physical-time constraints.
\end{maintheorem}
As the prescribed mean return energy \(E\) varies over \((1,\infty)\), the Hausdorff dimension of the corresponding joint level set attains every value in \((0,1]\).
For nonlinear branch geometry, Theorem~\ref{thm:nonlinear-stability}
gives quantitative bounds, while Theorem~\ref{thm:nonlinear-conditional}
replaces the affine quotient by a conditional-pressure equation that
retains the full geometric Lyapunov exponent.  The two-branch model
in Section~\ref{sec:stretched} acts on a single interval; the
constrained-pressure examples in Section~\ref{sec:nonlinearpressure}
are realized by genuine first-return towers, and are not claimed to
arise from one smooth interval transformation.
The needed thermodynamic hypotheses are verified for an explicit
non-affine, non-conjugate branch family in
Proposition~\ref{prop:nonlinear-example}; its dimension function is
analytic and has a unique maximum.

\subsection{Organization}
Section~\ref{sec:prelim} proves the classification of integer-set
spectra.  Section~\ref{sec:inducing} establishes the guarded insertion
lemma and its full-dimensional dynamical consequences.
Section~\ref{sec:stretched} proves the energy obstruction, the general
spectrum transition, the exact logarithmic transition, and the
interval-map application.  Section~\ref{sec:secondorder} proves the affine formula for complete
energy levels and then treats prescribed-schedule Cantor subsets;
Section~\ref{sec:nonlinearpressure} treats nonlinear geometry and
conditional pressure.

\medskip
\noindent\textbf{Asymptotic notation.}
For nonnegative quantities $A$ and $B$, we write $A\lesssim B$
if $A\le CB$ for some constant $C>0$, and $A\gtrsim B$ if
$B\lesssim A$. We write $A\asymp B$ if both $A\lesssim B$ and
$B\lesssim A$, whereas in an asymptotic statement $A\sim B$ means
$A/B\to1$ in the indicated limit. Equivalently, $A=O(B)$ means
$A\lesssim B$, while $A=o(B)$ means $A/B\to0$ in the indicated limit.
Unless stated otherwise, implicit constants may depend on fixed
parameters but not on the scale variables; any explicitly asserted
uniformity also applies to these constants.

\section{Macroscopic intermediate dimensions and scale constraints}\label{sec:prelim}
For a subset $E\subset\N$, let $S_n=[2^{n-1},2^n)\cap\N$. For $0\le\theta\le1$, $s\ge0$ and $N\ge1$ put $\cH_N^{s,\theta}(E)=\sum_{n\ge N}\cH_{n,N}^{s,\theta}(E)$, where
\begin{align*}
\cH_{n,N}^{s,\theta}(E)
=\inf\left\{\sum_i(r_i/2^n)^s:
 E\cap S_n\subset\bigcup_i B(x_i,r_i),\quad
 x_i\in S_n,\quad 1\le r_i\le2^{n-N\theta}\right\}.
\end{align*}
With the extended logarithm conventions $\log0=-\infty$, $\log\infty=\infty$, define
\begin{align*}
\lowD_\theta(E)&=\inf\biggl\{s\ge0:\liminf_{N\to\infty}\frac{\log\cH_N^{s,\theta}(E)}{N\log2}\le0\biggr\},\\
\upD_\theta(E)&=\inf\biggl\{s\ge0:\limsup_{N\to\infty}\frac{\log\cH_N^{s,\theta}(E)}{N\log2}\le0\biggr\}.
\end{align*}
These are the macroscopic intermediate dimensions introduced by Li and Xi \cite{LX}; see also the earlier discrete intermediate dimensions of Li, Xi and Zhang \cite{LXZ}, which use a different aggregation of shells. They agree with macroscopic Hausdorff dimension at $\theta=0$, are bounded by $1$ for subsets of $\N$, and are invariant under bounded perturbations and finite thickening of integer sets. The lower dimension requires the sum over all shells $n\ge N$, not an estimate for just one shell.

The following elementary observation explains why positive visit frequency forces every macroscopic intermediate dimension to take its maximal value.
\begin{lemma}\label{lem:density}
If $\limsup_{R\to\infty}R^{-1}\card(E\cap[1,R])>0$, then
$\lowD_\theta(E)=\upD_\theta(E)=1$ for every $\theta\in[0,1]$.
\end{lemma}
\begin{proof}
Infinitely many dyadic shells contain at least $c2^n$ points. For an admissible cover of one such shell, $\sum_i r_i\ge c'2^n$ because a ball of radius $r_i\ge1$ covers at most $3r_i$ integers. As $s<1$ and $r_i\le2^{n-N\theta}$,
\[
\sum_i(r_i/2^n)^s\ge2^{N\theta(1-s)}\sum_i(r_i/2^n)\ge c'2^{N\theta(1-s)}.
\]
There are infinitely many such shells for each $N$, so $\cH_N^{s,\theta}(E)=\infty$ whenever $s<1$. Conversely, unit-ball covers give $\cH_N^{s,\theta}(E)\lesssim\sum_{n\ge N}2^{n(1-s)}\to0$ for $s>1$.
\end{proof}

The comparison below is standard for macroscopic intermediate dimensions: the dimension inequality is Proposition~2.6 of Li and Xi \cite{LX}, and its covering-content form follows by the same subdivision argument. We recall it to fix the precise scale constraint used later.
\begin{lemma}\label{lem:compare}
For $0\le\theta_1<\theta_2\le1$ and $0<s\le s'\le1$,
\[
\cH_N^{s',\theta_2}(E)\le C\,2^{N[(\theta_2-\theta_1)(1-s)-\theta_2(s'-s)]}\cH_N^{s,\theta_1}(E).
\]
Consequently, for either choice of sign $\pm$ and $0<\theta_1<\theta_2\le1$,
\begin{equation*}
0\le \DD^\pm_{\theta_1}(E)\le \DD^\pm_{\theta_2}(E)\le
1-\frac{\theta_1}{\theta_2}\bigl(1-\DD^\pm_{\theta_1}(E)\bigr).
\end{equation*}
\end{lemma}
\begin{proof}
Fix a shell $n\ge N$ and put $U_i=2^{n-N\theta_i}$ for $i=1,2$.
A ball of radius $r\le U_2$ remains admissible, and its normalized
$s'$-contribution is at most $2^{-N\theta_2(s'-s)}$ times its
$s$-contribution.  If $U_2<r\le U_1$, subdivide its intersection
with the shell into $O(r/U_2)$ balls of radius $U_2$, with centres
in the shell.  The ratio of the new contribution to the old one
is at most
\[
 C2^{-n(s'-s)}r^{1-s}U_2^{s'-1}
 \le C2^{N[\theta_2(1-s')-\theta_1(1-s)]}.
\]
This last factor also bounds the ratio in the first case, since
$(\theta_2-\theta_1)(1-s)\ge0$.  Taking infima over covers and
summing over $n\ge N$ proves the content inequality.
The dimensions are nondecreasing because increasing $\theta$
restricts the permitted covers.  For $\theta_1>0$, choose
$s'=1-(\theta_1/\theta_2)(1-s)$, so that the exponential factor
equals one, and pass to the critical exponents in the definitions.
This gives the asserted comparison for both signs.
\end{proof}

In particular, Lemma~\ref{lem:compare} shows that every function
occurring as both the lower and upper macroscopic intermediate-dimension
spectrum must be scale-admissible, as defined in the Introduction.

\subsection{Two-scale packet geometry}

Fix $0<b<a<1$ and a positive integer $M$.  For sufficiently large $j$ let
\[
\ell_j=\lceil j^4\rceil,\qquad
Q_j=\left\lfloor\frac{2^{bj}}{\ell_j}\right\rfloor,\qquad
G_j=\left\lceil 2^{(a-b)j}\ell_j\right\rceil.
\]
A packet at level $j$ consists of $\ell_j$ successive integer points, with adjacent gaps in $[1,M]$.  There are $Q_j$ packets, separated by jumps of size $G_j$ (up to an additive $O(M\ell_j)$ variation). Place this configuration, denoted by $E_j$, at $2^j$ and assume $E_j\subset [2^j,2^j+C2^{aj}]$.  Set $E=\bigcup_{j\ge j_0}E_j$. All constants below may depend on $a,b,M$.
The elementary packet construction separates the shell count from the spatial spread of return times.

For $0<b<a<1$ set
\begin{equation}\label{zcqx}
P_{a,b}(0)=0,\qquad
P_{a,b}(\theta)=\min\left\{b,\left(1-\frac{1-a}{\theta}\right)_+\right\},\quad0<\theta\le1.
\end{equation}
Here and throughout, \(u_+=\max\{u,0\}\).
These elementary spectra will serve as support curves in the classification
of admissible functions.  Their dynamical realizations are treated in
Section~\ref{sec:inducing}.
\begin{lemma}\label{lem:packets}
The set $E$ has zero natural density and
\[
\lowD_\theta(E)=\upD_\theta(E)=P_{a,b}(\theta),\qquad0\le\theta\le1.
\]
Furthermore, for every $0<\theta\le1$ and every $0<s<P_{a,b}(\theta)$ there is $\varepsilon>0$ such that
$
\cH_N^{s,\theta}(E)\ge 2^{\varepsilon N}
$ for all sufficiently large $N.$
For every \(s>P_{a,b}(\theta)\), the quantity \(\mathcal H_N^{s,\theta}(E)\) decays exponentially in \(N\). 
\end{lemma}
\begin{proof}
Write $R_j=2^j$. The number of points in $E_j$ is
$
K_j\asymp Q_j\ell_j\asymp R_j^b$ and $ \operatorname{diam}E_j\asymp Q_jG_j\asymp R_j^a.$
An interval $I$ of length $r\ge1$ meets at most $C(1+r/G_j)$ packets, and each packet has $O(\ell_j)$ points. Consequently
\begin{equation}\label{eq:localcount}
\card(E_j\cap I)\le C\left(\ell_j+rR_j^{-(a-b)}\right).
\end{equation}
The total number of points up to $R$ is $O(R^b)$; thus the natural density vanishes. For $j$ large, $E_j$ is contained in the single shell $S_{j+1}=[2^j,2^{j+1})$.

For the upper dimensional bound, it suffices to consider \(P_{a,b}(\theta)<s<1\), since \(P_{a,b}(\theta)<1\) and \(\mathcal H_N^{s,\theta}(E)\) is nonincreasing in \(s\). Let \(n=j+1\ge N\). Unit-radius balls give
\begin{equation}\label{eq:unit}
\cH_{n,N}^{s,\theta}(E)\le C2^{(b-s)j}.
\end{equation}
Put $U=2^{n-N\theta}$. If $U\le C2^{aj}$, covering the containing interval with $O(2^{aj}/U)$ balls of radius $U$ gives
\begin{equation}\label{eq:wide}
\cH_{n,N}^{s,\theta}(E)\le C2^{(a-1)j+N\theta(1-s)}.
\end{equation}
If $U>C2^{aj}$, one ball of radius $C2^{aj}$ (or an inessential bounded multiple of such balls near a shell boundary) gives
\begin{equation}\label{eq:one}
\cH_{n,N}^{s,\theta}(E)\le C2^{-(1-a)sj}.
\end{equation}
For $s>b$, summing \eqref{eq:unit} over $j\ge N-1$ decays exponentially. Alternatively, when $\theta>1-a$ and $s>1-(1-a)/\theta$, sum \eqref{eq:wide} until $U\asymp2^{aj}$ and \eqref{eq:one} thereafter. The resulting bound is
$C(2^{N[a-1+\theta(1-s)]}+2^{-c_sN})$, which also decays exponentially. When $0\le\theta\le1-a$, one-ball covers suffice and imply dimension zero.

For the lower bound fix $0<\theta\le1$, and $0<s<\min\{b,1-(1-a)/\theta\}$. Take the shell $n=N+1$ (thus $j=N$), and place normalized counting mass $\nu_N$ on $E_N$. By \eqref{eq:localcount}, for every allowed radius $1\le r\le U=2^{n-N\theta}$,
\begin{align*}
\frac{\nu_N(B(x,r))}{(r/2^n)^s}
&\le C\left(N^4 2^{(s-b)N}r^{-s}
 +2^{(s-a)N}r^{1-s}\right)\\
&\le C\left(N^4 2^{-(b-s)N}+2^{[1-a-\theta(1-s)]N+O(1)}\right)
 \le C2^{-\varepsilon N},
\end{align*}
where $\varepsilon>0$ is smaller than both positive quantities $b-s$ and $\theta(1-s)-(1-a)$. Here $0<s<1$ is used to maximize $r^{1-s}$ at $r=U$.
Every admissible cover of $E\cap S_{N+1}$ therefore obeys
$1\le C2^{-\varepsilon N}\sum_i(r_i/2^{N+1})^s$.
Hence $\cH_N^{s,\theta}(E)\ge c2^{\varepsilon N}$, proving the exponential lower bound and the required lower dimension. Finally, for $\theta=0$ all $E_j$ are covered by a bounded number of balls of radius $O(2^{aj})$; their total normalized $s$-content is $O(2^{-(1-a)sN})$, so both dimensions vanish.
\end{proof}

\begin{remark}\label{rem:why}
If all $K_j\asymp2^{bj}$ visits were separated by individual excursions of length $\asymp 2^{(a-b)j}$, the sum of logarithmic contractions of the forced branches would be $\asymp j2^{bj}$, comparable to (or larger than) the free induced entropy scale. Grouping $\ell_j=j^4$ visits into each packet reduces the number of forced long branches to $Q_j\asymp2^{bj}/j^4$, and their logarithmic geometric contribution to $O(2^{bj}/j^3)$.
\end{remark}

\subsection{Bounded-gap clusters}
A finite integer set $B=\{b_1<\cdots<b_m\}\subset\N$ is called a
\emph{bounded-gap cluster with gap bound $M$} if $M\ge1$ and
$b_{i+1}-b_i\le M$ for $1\le i<m$.
For a family of such clusters, a \emph{uniform gap bound} means that
the same constant $M$ works for every cluster. We consider sets of the form
$E=\bigcup_{j\ge j_0}B_j$, with clusters located at geometrically
increasing scales.
A concrete example is obtained by fixing $0<a<1$ and $\rho>1$,
and choosing $j_0$ sufficiently large. For every $j\ge j_0$, let
\[
 E=E_{a,\rho}=\bigcup_{j=j_0}^{\infty}B_j  \quad \text{where $B_j=\bigl\{\lfloor\rho^j\rfloor,\ldots,
 	\lfloor\rho^j\rfloor+\lfloor\rho^{aj}\rfloor-1\bigr\}.$}
\]
Each $B_j$ is a block of $\lfloor\rho^{aj}\rfloor$ consecutive
integers beginning at $R_j=\lfloor\rho^j\rfloor$, so the uniform
gap bound is $M=1$. Its diameter is $O(R_j^a)$, and successive
blocks are separated by distances comparable to $R_j$. Moreover,
$\#(E_{a,\rho}\cap[1,T])=O(T^a)$, so $E_{a,\rho}$ has zero natural
density. The next lemma gives its macroscopic intermediate-dimension
spectrum and extends the conclusion to general bounded-gap clusters.

For $0<a<1$, define
\begin{equation}\label{yingshe}
\varphi_a(0)=0,\qquad
\varphi_a(\theta)=\left(1-\frac{1-a}{\theta}\right)_+\quad(0<\theta\le1).
\end{equation}
The single-parameter spectrum \(\varphi_a\) is the limiting case of the two-parameter family \(P_{a,b}\) as \(b\uparrow a\). Its realization by bounded-gap clusters is treated separately below.
\begin{lemma}\label{lem:clusters}
Let $E=\bigcup_{j\ge j_0}B_j\subset\N$, where the $B_j$ are nonempty
bounded-gap clusters with a uniform gap bound. Put $R_j=\min B_j$ and
suppose $R_j\asymp\rho^j$ for some $\rho>1$. Assume that
$\#B_j\asymp R_j^a$, that each $B_j$ lies in an interval of length
$O(R_j^a)$, and that different clusters lie in disjoint,
macroscopically separated intervals. Then
\[
\lowD_\theta(E)=\upD_\theta(E)=\varphi_a(\theta),\qquad 0\le\theta\le1.
\]
\end{lemma}
\begin{proof}
At $\theta=0$, one ball of radius $O(R_j^a)$ covers each cluster, and $\sum_j R_j^{-(1-a)s}<\infty$ for every $s>0$. For $\theta>0$, in shell $n$ the largest allowed radius is $U_{n,N}=2^{n-N\theta}$. Covering a cluster by balls of radius $\min\{U_{n,N},CR_j^a\}$ yields, for $0<s<1$, an upper bound (up to constants)
\[
\min\{2^{(a-s)n},\quad
 2^{(a-1)n+N\theta(1-s)}+2^{-(1-a)sn}\}.
\]
Summing over occupied shells shows $\cH_N^{s,\theta}(E)\to0$ whenever $s>\varphi_a(\theta)$. For the lower bound when $\theta>1-a$, in some shell $n=N+O(1)$ at least a fixed proportion of a cluster is present, and an admissible covering satisfies $\sum r_i\gtrsim 2^{an}$. Since $s<1$ and $r_i\le U_{n,N}$,
\[
\cH_N^{s,\theta}(E)\gtrsim
2^{-sn}U_{n,N}^{s-1}2^{an}
\asymp2^{N[a-1+\theta(1-s)]}.
\]
For $s<\varphi_a(\theta)$ this grows exponentially. If $\theta\le1-a$, the lower bound zero is automatic. The bounded gaps between logarithmic scales of consecutive clusters ensure the lower estimate applies at every outer scale $N$.
\end{proof}

\subsection{Realization of admissible macroscopic dimension spectra}\label{sec:profile}
The two-parameter family supplies the elementary support curves needed to approximate more general spectra. The corresponding realization question is analogous in spirit, though not in its covering definition or admissibility condition, to the intermediate-dimension theory for bounded sets introduced by Falconer, Fraser and Kempton \cite{FFK} and the classification of attainable functions by Banaji and Rutar \cite{BanajiRutar}. Here the construction takes place first in integer time; the dynamical statement additionally invokes the guarded insertion hypotheses of Section~\ref{sec:inducing}.

For $0<b<a<1$, we use the supporting curves $P_{a,b}$ defined by \eqref{zcqx}.
Lemma~\ref{lem:packets} realizes \(P_{a,b}\) by a subset of \(\NN\) with zero natural density, with matching upper and lower macroscopic intermediate dimensions. Its proof gives exponentially separated lower and upper covering estimates at every positive parameter, uniformly over a fixed compact range of \(a,b\). The special case $a=1$ is defined separately by
$P_{1,b}(\theta)=b$ for \emph{every} $0\le\theta\le1$; it is obtained
by spreading the packets across a fixed fraction of each dyadic shell.
This is not the pointwise continuation at $\theta=0$ of the formula
for $a<1$: the latter always has $P_{a,b}(0)=0$. The all-zero spectrum is realized by a subpower number of visits per shell.

Two-parameter supporting curves provide a representation of any admissible dimension spectrum.
\begin{lemma}\label{lem:support}
Let \(h\) be scale-admissible with \(h(1)<1\). For each
\(t\in(0,1)\) with \(h(t)>0\), put
$
b_t=h(t)$ and $ a_t=1-t(1-h(t)).$
Then \(b_t<a_t<1\) whenever \(h(t)>0\) and \(t<1\), and
\[
P_{a_t,b_t}(\theta)\le h(\theta)\quad(0\le\theta\le1),\qquad
P_{a_t,b_t}(t)=h(t).
\]
If \(h(1)<1\), then
\begin{equation}\label{eq:sup}
h(\theta)=\max\left\{h(0),\ \sup_{t\in\mathbb Q\cap(0,1)}P_{a_t,b_t}(\theta)\right\},\quad\theta>0,
\end{equation}
where zero-valued spectra are omitted. The supremum is interpreted as zero if empty. For \(\theta=0\), the right side is \(h(0)\).
\end{lemma}
\begin{proof}
For \(\theta\ge t\), monotonicity gives \(P_{a_t,b_t}(\theta)\le b_t=h(t)\le h(\theta)\). For \(\theta<t\), the monotonicity of \(\theta(1-h(\theta))\) gives \(1-t(1-h(t))/\theta\le h(\theta)\), so the same conclusion follows. Equality at \(\theta=t\) is immediate. Finally \(h\) is continuous at every positive parameter, so a dense sequence of rational supporting points recovers its values. The constant spectrum \(h(0)\) handles the endpoint at zero. If \(h(1)=1\), scale-admissibility forces \(h(\theta)=1\) for every \(\theta>0\), a special case dealt with separately.
\end{proof}

The following estimate isolates the uniformity needed when the packet
parameters vary from one dyadic shell to the next.  In particular, it
keeps track of the outer index $N$ in the sum over all shells
$n\ge N$, rather than using a single-shell estimate.
\begin{lemma}\label{lem:uniformdiagonal}
Let the packet of Lemma~\ref{lem:packets} in shell $S_{j+1}$ have
parameters $0<b_j<a_j\le1$, with packet size $\ell_j=\lceil j^4\rceil$
and with a sufficiently small fixed packing constant when $a_j=1$.
Assume
\[
 j\min\{b_j,a_j-b_j,1-b_j\}\ge\sqrt j,
 \qquad j(1-a_j)\ge\sqrt j\quad\text{if }a_j<1.
\]
The constants in the following inequalities are independent of
$j,b_j,a_j$ (once $j$ is large).  For $n=j+1\ge N$ and $0<s<1$,
\begin{equation}\label{eq:uniformupper}
 \cH_{n,N}^{s,\theta}(E_j)
 \le Cj^{C_0}\min\left\{2^{(b_j-s)j},
 \begin{cases}
 2^{(a_j-1)j+N\theta(1-s)},&2^{j-N\theta}\le2^{a_jj},\\
 2^{-s(1-a_j)j},&2^{j-N\theta}>2^{a_jj}
 \end{cases}\right\}.
\end{equation}
If $\nu_j$ is normalized counting measure on $E_j$, then for every
$1\le r\le 2^{j+1-N\theta}$,
\begin{equation}\label{eq:uniformlower}
 \frac{\nu_j(B(x,r))}{(r/2^j)^s}
 \le Cj^{C_0}\left(
  2^{(s-b_j)j}+2^{(1-a_j)j-N\theta(1-s)}
 \right).
\end{equation}
Here $C_0$ is absolute; bounded shifts from $j$ to $j+1$ are absorbed
in $C$.  A baseline block of $2^{o(j)}$ consecutive integers
has zero spectrum and contributes $2^{-cj+o(j)}$ to every
$s$-content with $s>0$ if it is contained in a time interval
of length $2^{o(j)}$.
\end{lemma}
\begin{proof}
For varying parameters the packet construction still has
$K_j=2^{b_jj+O(\log j)}$ and span
$2^{a_jj+O(\log j)}$, uniformly under the displayed separation
assumptions.  The normalized local count bound reads
$\nu_j(B(x,r))\le Cj^{C_0}(2^{-b_jj}+r2^{-a_jj})$.
Multiplying by $(2^j/r)^s$, maximizing $r^{-s}$ at $r=1$
and $r^{1-s}$ at the permitted upper radius gives
\eqref{eq:uniformlower}.  Covering by unit balls, by
$O(2^{a_jj}/U)$ balls of maximal radius $U$, or by one ball
when $U$ exceeds the packet span, gives \eqref{eq:uniformupper}.
The baseline assertion follows by covering the baseline with one
ball and summing over shells.
\end{proof}

The cyclic schedule in the classification is required to work at
every large outer scale, not merely along a selected subsequence.
The following observation makes this uniformity explicit.
\begin{lemma}\label{lem:epoch-recurrence}
Suppose $J_k\ge k^3$ and $J_{k+1}-J_k\ge k^2$ for every $k$, and,
throughout the epoch
$J_k\le j<J_{k+1}$, the types $1,\ldots,k$ occur successively
and periodically, one type per integer shell.  Then the active epoch
index $k(j)$ satisfies $k(j)=o(j)$.  Moreover, for each fixed
type $i$ and all sufficiently large integers $N$, there is a shell
index $j_N\ge N$ of type $i$ satisfying
$
 0\le j_N-N\le 2k(N)+2=o(N).
$
Consequently any packet lower bound with exponent
$2^{c j+o(j)}$ at the occupied shell of type $i$ transfers to
$2^{c'N}$, for every $0<c'<c$, at all sufficiently large outer
scales $N$, provided the permitted covering radii are compared using
that same $N$.
\end{lemma}
\begin{proof}
The inequality $J_k\ge k^3$ implies $k(j)\le j^{1/3}$.
If $N$ lies in an epoch of index $k\ge i$, the desired type
occurs within at most $k$ positions unless that epoch ends
first.  In the latter case fewer than $k$ positions remain in the
current epoch.  The next epoch has length at least $(k+1)^2$, so
it contains type $i$ within its first $k+1$ positions.  Adding the
two waiting times gives the bound $2k+2$.  Since $k\le N^{1/3}$,
we have $j_N=N+o(N)$, and
any fixed strict exponential margin in the single-shell estimate
survives the shift, with the outer-scale parameter $N$ retained.
\end{proof}

\begin{theorem}\label{thm:attain}
A function \(h:[0,1]\to[0,1]\) is simultaneously the lower and upper macroscopic intermediate-dimension spectrum of some \(E\subset\NN\) if and only if \(h\) is scale-admissible. Moreover, \(E\) can be chosen to have zero natural density.
\end{theorem}
\begin{proof}
Necessity is Lemma~\ref{lem:compare}. We explain the diagonal realization, including the endpoints. The function $h\equiv0$ is realized by the baseline described below. First suppose $0<h(1)<1$. Enumerate the nonzero supporting packets in \eqref{eq:sup} and add, if \(0<h(0)<1\), the constant packet type \((a,b)=(1,h(0))\). Choose integer epochs $J_k$ so rapidly increasing that
$J_k\ge k^3$, $J_{k+1}\ge J_k+k^2$, each $J_k$ is divisible by $k$, and the
separation inequalities below hold for all types $i\le k$ and
$j\ge J_k$.  On the epoch $J_k\le j<J_{k+1}$ cycle through types
$1,\ldots,k$ in order (if the constant type is present it is included
in the cycle; a finite list can be repeated to make an infinite enumeration).  By Lemma~\ref{lem:epoch-recurrence}, every fixed type occurs at a
shell $j=N+o(N)$ after every sufficiently large outer index $N$;
the largest active type index is $o(j)$.
Choose the epochs so that for all active nonconstant types at level
$j$,
$
\min\{b_i,a_i-b_i,1-a_i,1-b_i\}\,j\ge\sqrt j,
$
whenever the listed quantities are positive. Every shell also contains a baseline block of \(\lceil2^{\sqrt j}\rceil\) consecutive points, which has the zero spectrum and guarantees many freely selectable return symbols in the dynamical realization.

For a packet of type \(i\) in shell \(j\), use \(\ell_j\asymp j^4\) consecutive points per packet, \(Q_j\asymp2^{b_i j}/\ell_j\) packets, and interpacket spacing \(G_j\asymp 2^{(a_i-b_i)j}\ell_j\). A fixed small multiplicative constant keeps the construction within the shell, also for \(a_i=1\). The local estimate
\begin{equation}\label{eq:packetlocal}
\#(E_j\cap I)\le C\bigl(\ell_j+|I|2^{-(a_i-b_i)j}\bigr)
\end{equation}
contains an absolute constant independent of active type, after enlarging the harmless polynomial factor. Lemma~\ref{lem:uniformdiagonal} makes the estimates uniform, up to factors $2^{o(j)}$, for all active types. The baseline has an entirely negligible count at every fixed positive exponent and is contained in an interval of length \(2^{\sqrt j}\).

For fixed \(\theta>0\) and \(s<h(\theta)\), choose a supporting type \(i\) for which \(P_{a_i,b_i}(\theta)>s\). Its occurrence in a shell \(j=N+o(N)\), guaranteed by
Lemma~\ref{lem:epoch-recurrence}, together with
\eqref{eq:packetlocal}, gives \(\cH_N^{s,\theta}(E)\ge2^{c_sN}\) for all large \(N\). Conversely, if \(s>h(\theta)\), for every active type either \(b_i\le h(\theta)\), when unit-ball covers yield \(2^{-(s-h(\theta))j+o(j)}\), or \(b_i>h(\theta)\), in which case \(P_{a_i,b_i}(\theta)\le h(\theta)\) forces \(1-a_i\ge\theta(1-h(\theta))\); the maximal-radius cover gives a summable bound \(2^{-c_sj+o(j)}\), uniformly for \(j\ge N\). Thus both dimensions equal \(h(\theta)\). At \(\theta=0\), nonconstant packets are each covered by a single ball with contribution at most \(C2^{-s\sqrt j}\), whose tail tends to zero; the constant type has dimension exactly \(h(0)\), by the per-shell mass estimate at \(a=1\). This gives the endpoint value.

If \(h(1)=1\), then \(h(\theta)=1\) for every \(\theta>0\). To realize \(h(0)<1\), use in the shells not reserved for the constant type the varying packet parameters \(a_j=1-j^{-1/2}\), \(b_j=1-2j^{-1/2}\); the preceding estimates give spectrum zero at zero and one at every positive parameter. If \(h(0)=1\), distribute \(\asymp2^j/j^4\) points approximately evenly over a fixed fraction of each shell, yielding the constant-one spectrum while the natural density still tends to zero. Finally, the all-zero function is realized by the baseline alone. All cases have zero natural density.
\end{proof}

\section{Dimension-preserving insertions and dynamical realization}\label{sec:inducing}

\subsection{Tempered insertions and Euclidean dimension}\label{sec:guards}

We first formulate the geometric input independently of any particular
intermittent map.  Let $Y$ be a compact interval and let
$F:\bigcup_{n\ge1}Y_n\longrightarrow Y$ be a piecewise smooth full-branch
expanding map. Its inverse branches are denoted by $\psi_n:Y\to Y_n$.
The partition intervals have disjoint interiors, and we assume:
\begin{enumerate}[label=\textup{(G\arabic*)},leftmargin=*]
\item\label{hyp:G1} $\sup_n\|\psi_n'\|_\infty\le\lambda<1$, and there
is $D<\infty$ such that
\[
 \sup_w\sup_{x,z\in Y}\bigl|\log|\psi_w'(x)|
       -\log|\psi_w'(z)|\bigr|\le D,
\]
where $w$ ranges over finite words.
\item\label{hyp:G2} Some $\psi_c(Y)$ is compactly contained in
$\operatorname{int}Y$. For $K_M$ the limit set generated by
$\{\psi_1,\ldots,\psi_M\}$, assume
$\sup_M\dim_H K_M=\dim_H K$, where $K$ is the countable-system
limit set.
\end{enumerate}
Conditions \textup{(G1)} and \textup{(G2)} are the common
\emph{geometric} hypotheses: the former gives uniform contraction
and distortion, and the latter provides Euclidean separation by a fixed
guard branch and approximation of the full limit-set dimension by
finite subsystems. Neither condition imposes a polynomial bound on
branch derivatives. The latter is recorded separately as the optional
condition \textup{(P)} below.

The following condition concerns the \emph{original dynamics}, rather
than the size of its inverse branches:
\par\smallskip
\noindent\textbf{(R) Genuine return lengths.}
The map \(F\) is the first-return map of an underlying dynamical system \(T\) to \(Y\). Moreover, for every \(n\ge1\), there exists a full branch \(Y_n\) on which the first-return time is exactly \(n\).
\par\smallskip
Thus, the branch label \(n\) represents the actual number of iterates of \(T\) between successive returns to \(Y\), rather than an arbitrarily assigned return-time parameter.

We keep branch growth separate from the common geometric assumptions
\textup{(G1), (G2)} and the genuine-return assumption \textup{(R)}.
The optional polynomial lower bound \textup{(P)} and the more general
logarithmic-power condition \textup{(L$_q$)} are
\begingroup
\makeatletter
\renewcommand{\tagform@}[1]{\maketag@@@{(\ignorespaces#1\unskip\@@italiccorr)\hspace{0.9em}}}
\makeatother
\begin{align}
\tag{P}\label{eq:polynomialbranch}
 &\inf_{y\in Y}|\psi_n'(y)|\ge A(n+1)^{-p}
   &&(n\ge1;\ A,p>0),\\
\tag{L$_q$}\label{eq:logpowerbranch}
 &-\log\inf_{y\in Y}|\psi_n'(y)|
      \le C(1+\log(n+1))^q
   &&(n\ge1;\ C>0,\ q\ge1).
\end{align}
\endgroup
Condition \textup{(P)} implies \textup{(L$_1$)} (and conversely),
after adjusting constants, whereas \textup{(L$_q$)} is generally weaker
when $q>1$. The general realization theorem in this section uses
\eqref{eq:logpowerbranch}; the polynomial condition is only an easily checked
special case. In contrast, the stretched-exponential regime in
Section~\ref{sec:stretched} uses the independent two-sided bounds
\eqref{eq:twosidedstretched}. Those bounds are incompatible with
\textup{(P)} when $\delta>0$, since $e^{-cn^\delta}=o(n^{-p})$
for every fixed $p>0$.

The Pomeau--Manneville inducing schemes treated in
Subsection~\ref{sec:lsv} satisfy (G1), (G2), (R), and (P).
The piecewise affine full-branch tower below illustrates the guard
and finite-subsystem assumptions, while the stretched-exponential
tower and the two-branch interval model satisfy (R), (G1), and (G2)
without (P). A smooth conjugate preserves corresponding geometric
estimates if the conjugacy and its inverse have derivatives bounded
above and away from zero. We do not assume any of these properties
for an arbitrary nonuniformly expanding map without verification.

We use the Euclidean metric on $Y$.  For an induced itinerary
$\omega=(\omega_k)$, the return times in the original tower are
$t_k=\sum_{i=1}^k\omega_i$; this identification is assumed to be a genuine
first-return coding, not merely a freely assigned roof function on an
unrelated symbolic space.

The insertion rule must respect complete guarded blocks.  This restriction
ensures separation at the first free choice where two itineraries differ.
The following estimate treats all Euclidean radii, including those crossed
by a single very long prescribed return.

For countably many inverse branches, bounded distortion does not
imply a uniform bound on the number of level-$n$ cylinders meeting a
Euclidean ball: one ball may meet infinitely many cylinders of very
different diameters.  The correct geometric input is instead the
separation of two trajectories at their \emph{first differing free
block}.  We record this estimate explicitly before constructing a
measure on the surviving itineraries.

\begin{lemma}\label{lem:guarded-separation}
Let $\mathcal B$ be a finite family of free block words whose compact
images $\psi_B(Y)$ are pairwise disjoint, and suppose
$\inf_Y|\psi_B'|\ge e^{-b}$ for all $B\in\mathcal B$.
Prescribed words are inserted between free blocks by a causal rule,
so that the prescribed word before a given free block depends only
on earlier free blocks.  Assume that, after any prefix of $k$ free
blocks, the total negative logarithm of the lower derivative bounds
of all prescribed symbols up to the next free block is at most a
nondecreasing function $A(k)$.  Then two coded points $x,y$ whose
first differing free block is the $j$th one satisfy
\begin{equation}\label{eq:guard-separation}
 |x-y|\ge g\exp\{-(j-1)b-A(j-1)\},
\end{equation}
where $g>0$ depends only on the fixed free block family and an
optional common initial prescribed word.  In particular, the bound
is uniform over all free choices and stages.
\end{lemma}
\begin{proof}
Let $g_0$ be the minimum distance between images of distinct free
blocks; it is positive by compactness, finiteness and separation.
Before the first differing free block the two itineraries have the
same free and prescribed words, by causality.  The inverse map
associated with this common prefix has derivative bounded below by
$\exp\{-(j-1)b-A(j-1)\}$, using the chain rule and the assumed
lower derivative estimates.  Applying the mean value theorem to the
two distinct block images gives the conclusion with $g=g_0$; a fixed
common initial prefix merely decreases $g$ by a fixed positive factor.
\end{proof}

\begin{lemma}\label{lem:rigorous-forcing}
Assume \ref{hyp:G1} and \ref{hyp:G2}.  Set
$\gamma(n)=\max\{1,-\log\inf_Y|\psi_n'|\}.$
Let $\mathcal B$ be a finite family of equal-length words whose compact
images $\psi_B(Y)$ are pairwise separated.  Write $m=\#\mathcal B$ and
suppose
\begin{equation}\label{eq:blockbounds}
 0<a\le b,\qquad e^{-b}\le\inf_Y|\psi_B'|
 \le\sup_Y|\psi_B'|\le e^{-a}<1 \quad(B\in\mathcal B).
\end{equation}
Prescribed words are inserted only between complete blocks of
$\mathcal B$, and each insertion depends only on earlier free blocks.
Suppose a nondecreasing deterministic function $A(k)$ bounds, uniformly
in the free choices, the total $\gamma$-weight of prescribed symbols
before the $(k+1)$st free block, including partial insertions, and
\begin{equation}\label{eq:prefixpenalty}
                         A(k)=o(k).
\end{equation}
Then the resulting set has Hausdorff dimension at least $\log m/b$.
More generally, for any positive probability vector $(p_B)_{B\in\mathcal B}$,
the image of every full-measure subset of the Bernoulli block space has
dimension at least
\begin{equation}\label{eq:weighted-insertion}
 \frac{-\sum_Bp_B\log p_B}
      {\sum_Bp_B[-\log\inf_Y|\psi_B'|]}.
\end{equation}
In particular, for every $M\ge c$ and $\epsilon>0$, a sufficiently long
finite guarded block family over $\{1,\ldots,M\}$ can be chosen so that
\emph{every} insertion rule satisfying these conditions produces a set
of dimension at least $s_M-\epsilon$, where $s_M=\dim_H K_M$.
\end{lemma}
\begin{proof}
Let $g>0$ and the first-difference estimate
\eqref{eq:guard-separation} be as in
Lemma~\ref{lem:guarded-separation}.  All separation constants are
independent of the free itinerary and of the insertion stage; the
finite block family remains fixed.

Put $D_k=kb+A(k)$ and $\rho_k=(g/3)e^{-D_k}$.  Since $A$ is
nondecreasing, \eqref{eq:guard-separation} implies that, for each
coded point $x$, the ball $B(x,\rho_k)$ contains no point whose
first $k$ free blocks differ from those of $x$.  Give the free blocks
the uniform Bernoulli law and push it forward to the constructed
Euclidean set.  Then
\[
 \nu(B(x,\rho_k))\le m^{-k},\qquad
 -\log\rho_k=kb+A(k)+O(1)=bk+o(k).
\]
It is important that $\rho_{k+1}/\rho_k$ need not be uniformly
bounded below: the increment $A(k+1)-A(k)$ can be unbounded, even
when $A(k)=o(k)$.  What we do have is
$
 \frac{\log\rho_{k+1}}{\log\rho_k}\longrightarrow1.
$
For $\rho_{k+1}\le r<\rho_k$, monotonicity gives
$\nu(B(x,r))\le\nu(B(x,\rho_k))$.
Interpolating between consecutive logarithmic radii proves
\[
 \underline d_\nu(x)
 :=\liminf_{r\downarrow0}\frac{\log\nu(B(x,r))}{\log r}
 \ge\frac{\log m}{b}.
\]
The local-dimension form of the mass-distribution principle proves
the unweighted lower bound in the ambient Euclidean metric.

For the weighted assertion put
$b_B=-\log\inf_Y|\psi_B'|$.  Let $\mathbb P$ be the Bernoulli law
with positive block probabilities $p_B$, and write
$\nu=\pi_*\mathbb P$ for its push-forward under the coding map.
For a coded point with free blocks $B_1,B_2,\ldots$,
replace $\rho_k$ by
$
 r_k(x)=\frac g3\exp\left\{-\sum_{j=1}^k b_{B_j}-A(k)\right\}.
$
The first-difference argument shows that
$\nu(B(x,r_k(x)))\le\prod_{j=1}^kp_{B_j}$.
The strong law gives, for $\mathbb P$-almost every coding sequence,
\[
 -\frac1k\log\mathbb P([B_1\cdots B_k])\longrightarrow
     -\sum_Bp_B\log p_B,\qquad
 -\frac1k\log r_k(x)\longrightarrow\sum_Bp_B b_B.
\]
Because the alphabet is finite, the $b_B$ are bounded; together
with $A(k)=o(k)$ this implies
$\log r_{k+1}(x)/\log r_k(x)\to1$.
The same monotone interpolation yields
\[
 \underline d_\nu(x)\ge
 \frac{-\sum_Bp_B\log p_B}{\sum_Bp_Bb_B}
 \quad\text{for }\nu\text{-almost every }x.
\]
Restricting the Bernoulli law to any measurable full-measure subset
preserves these estimates.  The mass-distribution principle proves
\eqref{eq:weighted-insertion}, without any bounded-overlap assumption
on all cylinders at one symbolic level.

To obtain the last assertion, choose a finite-system ergodic measure
with entropy $h$ and Lyapunov exponent $\chi$ such that $h/\chi$
is arbitrarily close to $s_M$; this is the finite conformal pressure
variational principle.  If $s_M=0$ the assertion is immediate.
For arbitrarily small $\xi>0$ and sufficiently large $L$, typical
$L$-words form a family $\mathcal W_L$ satisfying
\[
 \#\mathcal W_L\ge e^{L(h-\xi)},\qquad
 e^{-L(\chi+\xi)-D}\le |\psi_w'(x)|
 \le e^{-L(\chi-\xi)+D}\quad(x\in Y).
\]
Append the same interior guard $c$ to every word.  Since $M\ge c$,
these remain words over the finite alphabet.  Their compact images lie
strictly inside distinct $L$-cylinders and hence have positive mutual
separation.  In \eqref{eq:blockbounds} one may take
$
 a=L(\chi-\xi)-D>0$ and $
 b=L(\chi+\xi)+D+\gamma(c).$
Thus 
\[
\frac{\log m}{b}
\ge
\frac{L(h-\xi)}
{L(\chi+\xi)+D+\gamma(c)}.
\]
Choosing $\xi$ small and then $L$ large makes this at least
$s_M-\epsilon$.  The block family is fixed after these choices;
no assertion of exact preservation for a fixed finite block length is
needed.  Counting free letters instead of blocks in
\eqref{eq:prefixpenalty} is equivalent, since their common length is fixed.
\end{proof}

\subsection{Full-dimensional realization of every admissible spectrum}
The static classification and the guarded insertion lemma combine to
give a realization theorem for the entire function of the scale
parameter. Unlike the polynomial special case stated below, the
result allows the more general logarithmic-power contraction law.

\begin{theorem}\label{thm:allprofiles}
Assume \textup{(G1), (G2), (R)} and \textup{(L$_q$)} for some
$q\ge1$. For every scale-admissible $h:[0,1]\to[0,1]$,
\[
 \dim_H\{y\in K_\infty:p(y)=0,\quad
       \mathsf D^-_\theta(R(y,Y))=\mathsf D^+_\theta(R(y,Y))
                 =h(\theta)\ (0\le\theta\le1)\}=\dim_H K.
\]
\end{theorem}
\begin{proof}
Fix $M\ge c$ and $\epsilon>0$, and choose guarded blocks of common
length $L$ as in Lemma~\ref{lem:rigorous-forcing}.  Use the cyclic
epoch construction in Theorem~\ref{thm:attain}, with every packet
length replaced by $\ell_j=L\lceil j^{q+3}/L\rceil$.  Each packet
consists of freely chosen complete guarded blocks, so its physical
duration lies between $\ell_j$ and $M\ell_j$.  Prescribe the starting
times of successive packets by the same static spacing, increasing a
fixed spacing factor if necessary.  After a packet has been read,
a single return fills the positive remaining gap to the next starting
time.  The final return fills the gap to exactly $2^{j+1}$.
This is a causal rule.  Small fixed packing constants leave room for
these jumps even for the constant supporting type $a=1$.
Put a baseline of at least $2^{\sqrt j}$ free letters in complete
blocks at the beginning of every stage.  It occupies $2^{o(j)}$ time
and has zero spectrum.

The packet estimates in Lemma~\ref{lem:uniformdiagonal} hold with a
different fixed power of $j$: their proofs only use the local count
bound, the number of packets, and their containing span.  Bounded
free roof values and rounding by $L$ change these estimates by
polynomial factors.  Theorem~\ref{thm:attain} then gives matching
bounds for the sums over all distant shells, yielding exactly $h$,
including $\theta=0$.
All visit counts remain $o(2^j)$.

We check the prefix condition without assuming comparability between
successive active types.  Every inserted return in stage $j$ has
length at most $2^{O(j)}$ and contraction contribution at most $Cj^q$.
Except for a bounded number of stage bridges, each such return follows
at least $\ell_j$ new free letters.  After any $k'$ free letters of
the current stage, its accumulated inserted contribution is therefore
at most
$
 Cj^{-3}k'+Cj^q.
$
Before that stage, the baseline alone has supplied at least
$2^{\sqrt{j-1}}$ free letters.  The error $Cj^q$ divided by this
quantity tends to zero.  Given any $\eta>0$, discard finitely many
stages so that $Cj^{-3}<\eta$; summing the preceding estimate, and
using $\sum_{i\le j}i^q=O(j^{q+1})=o(2^{\sqrt{j-1}})$, bounds the
total contribution at every prefix by $\eta$ times its free-letter
count plus a negligible term.  This proves \eqref{eq:prefixpenalty}.

For clarity the limiting cases use the same explicit choices as in
the static construction.  If $h(1)=1$ and $h(0)<1$, use
$a_j=1-j^{-1/2}$, $b_j=1-2j^{-1/2}$ on even shells, and the constant
type $(1,h(0))$ on odd shells when $h(0)>0$; otherwise use the varying
type on every shell.  These parameters resolve all polynomial factors.
If $h(0)=1$, choose any fixed $p>1$ and put
$Q_j\asymp2^j/j^{p+q+3}$ packets of $\ell_j\asymp j^{q+3}$ letters,
with starting times spaced by $\asymp j^{p+q+3}$ across a fixed
fraction of the shell.  The count is $\asymp2^j/j^p$, and the same
local mass bound as in Theorem~\ref{thm:attain} gives constant
spectrum one.  The preceding prefix estimate applies in both cases.
For $h\equiv0$, the baseline and one terminal bridge suffice.

Lemma~\ref{lem:rigorous-forcing} now gives dimension at least
$s_M-\epsilon$.  Let $\epsilon\downarrow0$ and $M\to\infty$ and
use \textup{(G2)}.
\end{proof}

The optional polynomial lower bound (P) is a concrete sufficient
condition for the general theorem, not an additional geometric
assumption. The original two-parameter realization statement follows
without a separate packet construction.
\begin{corollary}\label{cor:polynomial-realization}
Assume \textup{(G1), (G2), (R)} and the optional
condition \textup{(P)}. For any
$0<b<a<1$, let
\[
\mathscr E(a,b)=\bigl\{y\in K_\infty:p(y)=0,\quad
 \lowD_\theta(R(y,Y))=\upD_\theta(R(y,Y))
   =P_{a,b}(\theta)\ \text{for all }\theta\in[0,1]\bigr\}.
\]
Then $\dim_H\mathscr E(a,b)=\dim_H K$. In particular, this
dimension equals one whenever $\dim_H K=1$.
\end{corollary}
\begin{proof}
By \eqref{eq:polynomialbranch},
$-\log\inf_Y|\psi_n'|\le p\log(n+1)-\log A$,
so \textup{(L$_1$)} holds after increasing its constant.
The function $P_{a,b}$ is scale-admissible; apply
Theorem~\ref{thm:allprofiles}.
\end{proof}

We now apply the general realization theorem to a classical family of intermittent maps whose first-return structure satisfies the required assumptions.
\subsection{Pomeau--Manneville maps}\label{sec:lsv}

For $\alpha>0$, consider
\[
 T_\alpha(x)=\begin{cases}
 x(1+2^\alpha x^\alpha),&0\le x<1/2,\\
 2x-1,&1/2\le x\le1,
 \end{cases}\qquad Y=[1/2,1].
\]
The first-return map $F=T_\alpha^\tau$ has one (modulo endpoints)
branch of every return length $n\ge1$; its inverse branches are uniformly
contracting and have uniform distortion bounds, with
\[
 C_\alpha^{-1}(n+1)^{-1-1/\alpha}
 \le |\psi_n'(y)|\le
 C_\alpha(n+1)^{-1-1/\alpha}\quad(y\in Y).
\]
These are the standard properties of the LSV inducing construction
\cite{LSV,BSTV}. Finite-subsystem approximation gives
$\lim_{M\to\infty}\dim_H K_M=1$ by the countable conformal IFS theory
\cite{MU}.

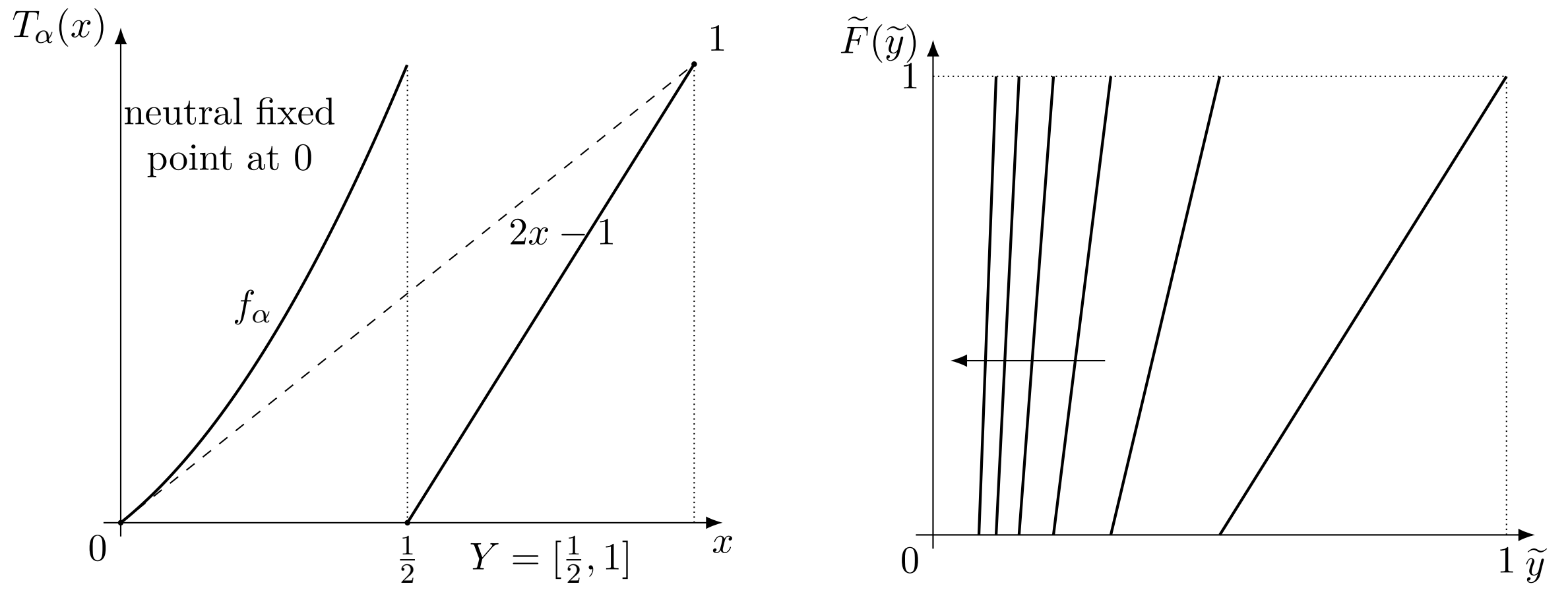
\begin{figure}[t]
\centering
\begin{minipage}[t]{0.48\textwidth}
\centering
\begin{tikzpicture}[x=5.75cm,y=4.6cm,>=Latex]
  \draw[->] (-0.03,0) -- (1.05,0) node[below] {$x$};
  \draw[->] (0,-0.03) -- (0,1.08) node[left] {$T_\alpha(x)$};
  \draw[dashed] (0,0) -- (1,1);
  \draw[thick,domain=0:0.5,samples=120,smooth] plot (\x,{\x*(1+2*(\x))});
  \draw[thick,domain=0.5:1,samples=2] plot (\x,{2*\x-1});
  \draw[densely dotted] (0.5,0) -- (0.5,1);
  \draw[densely dotted] (1,0) -- (1,1);
  \fill (0,0) circle (0.8pt) node[below left] {$0$};
  \fill (0.5,0) circle (0.8pt) node[below] {$\tfrac12$};
  \fill (1,1) circle (0.8pt) node[above right] {$1$};
  \node[below] at (0.75,0) {$Y=[\tfrac12,1]$};
  \node at (0.23,0.47) {$f_\alpha$};
  \node at (0.77,0.63) {$2x-1$};
  \node[align=center] at (0.19,0.84) {neutral fixed\\point at $0$};
\end{tikzpicture}

\smallskip
\small (a) The map $T_\alpha$, with neutral branch $f_\alpha(x)=x(1+2^\alpha x^\alpha)$, drawn for $\alpha=1$.
\end{minipage}
\hfill
\begin{minipage}[t]{0.48\textwidth}
\centering
\begin{tikzpicture}[x=5.75cm,y=4.6cm,>=Latex]
  \draw[->] (-0.03,0) -- (1.05,0) node[below] {$\widetilde y$};
  \draw[->] (0,-0.03) -- (0,1.08) node[left] {$\widetilde F(\widetilde y)$};
  \draw[densely dotted] (0,1) -- (1,1);
  \draw[densely dotted] (1,0) -- (1,1);
  \foreach \a/\b in {0.50/1.00,0.31/0.50,0.21/0.31,0.15/0.21,0.11/0.15,0.08/0.11}{
     \draw[thick] (\a,0) -- (\b,1);
  }
  \node[below left] at (0,0) {$0$};
  \node[below] at (1,0) {$1$};
  \node[left] at (0,1) {$1$};
  \draw[->] (0.30,0.38) -- (0.03,0.38);
\end{tikzpicture}

\smallskip
\small (b) Induced first-return map on $Y$ with one full branch for each return length.
\end{minipage}
\caption{The Pomeau--Manneville model and its induced first-return structure. The right panel is schematic, with $Y=[1/2,1]$ rescaled to $[0,1]$ on both axes: $\widetilde y=2y-1$ and $\widetilde F(\widetilde y)=2F((\widetilde y+1)/2)-1$. The branches accumulate at the left endpoint.}
\label{fig:pm-inducing}
\end{figure}
\begin{corollary}\label{cor:allalpha}
For every $\alpha>0$ and every scale-admissible $h:[0,1]\to[0,1]$,
\[
 \dim_H\{y\in Y:p(y)=0,\quad
 \mathsf D^-_\theta(R(y,Y))=\mathsf D^+_\theta(R(y,Y))=h(\theta)
 \text{ for all }\theta\}=1.
\]
In particular this holds for every two-parameter spectrum
$P_{a,b}$, $0<b<a<1$.
\end{corollary}
\begin{proof}
The branch estimates above verify \textup{(R)}, \textup{(G1)},
\textup{(G2)}, and \textup{(P)}, hence the hypotheses of
Theorem~\ref{thm:allprofiles} with $q=1$. The induced limit set
has dimension one.  No finite invariant measure for the original map
is required.  For $\alpha\ge1$ the assertion still concerns
exceptional itineraries and their Hausdorff dimension.
\end{proof}

\begin{corollary}
For every $\alpha>0$ there is a set of initial points of Hausdorff
dimension one such that their physical return-time set $R$ has zero
natural density and
\[
 \lowD_0(R)=\upD_0(R)=0,\quad
 \lowD_{1/2}(R)=\upD_{1/2}(R)=\tfrac14,\quad
 \lowD_1(R)=\upD_1(R)=\tfrac12.
\]
Indeed, take $a=5/8$ and $b=1/2$.
\end{corollary}

\section{Stretched-exponential branch geometry}\label{sec:stretched}
Long excursions become geometrically expensive if inverse branches
shrink at stretched-exponential rates. The obstruction below requires
only the stated contraction lower bounds; the subsequent full-dimensional
realizations use the geometric conditions (G1) and (G2) together
with (R), but do not assume the polynomial bound (P).
The argument proceeds from an energy obstruction for arbitrary
itineraries to dimension-preserving constructions with balanced
micro-packets, and then to exact counting laws and their critical
logarithmic correction.  The polynomial packets used for static
spectrum realization are not sufficient in this contraction regime:
when interpreted as first returns, their long separators can force
zero Hausdorff dimension.  The micro-packet construction controls
this accumulated contraction while retaining the required spectrum.

Let \(F\) be a countable full-branch map of a compact interval \(Y\), with one inverse branch \(\psi_n\) per allowable first-return length \(n\ge1\). Suppose for some \(\delta>0\), \(c,C>0\),
\begin{equation}\label{eq:stretchedbound}
\|\psi_n'\|_\infty\le C e^{-cn^\delta}.
\end{equation}
Recall that $K_\infty\subset K\subset Y$ consists of points whose
orbits return to $Y$ infinitely often. For $y\in K_\infty$, let
$\tau_i(y)$ denote the number of iterates of $T$ between its
$(i-1)$-st and $i$-th visits to $Y$, counting the visit at time zero
as the initial one. By the genuine first-return coding in
Section~\ref{sec:guards}, $\tau_i(y)$ is the $i$-th branch label. Set
\[
 t_0(y)=0,\qquad t_k(y)=\sum_{i=1}^k\tau_i(y)\quad(k\ge1).
\]
Then the positive return-time set defined in the Introduction is
$R(y,Y)=\{t_k(y):k\ge1\}$. Throughout this section all level sets
are restricted to $K_\infty$; points making only finitely many
returns are excluded.

For a uniformly discrete set $E$, Li and Xi
\cite[Proposition~2.4]{LX} identify the upper endpoint of the
macroscopic intermediate dimensions with the upper mass dimension.
In particular, since an unbounded subset $E\subset\mathbb N$ is
uniformly discrete, their result yields
\[
\mathsf D^+_1(E)=\limsup_{R\to\infty}
 \frac{\log\#(E\cap[1,R])}{\log R}.
\]
This standard endpoint identity is used below without repeating its
proof.

The appropriate universal invariant for dimension loss is the accumulated
geometric penalty, rather than only the logarithmic counting exponent.
The following criterion is insensitive to fluctuations between different
subsequences of return generations.
\begin{theorem}\label{thm:energy}
Under \eqref{eq:stretchedbound}, let
\[
 \mathscr Z_\delta=\left\{y\in K_\infty:
 \limsup_{k\to\infty}\frac1k\sum_{i=1}^k\tau_i(y)^\delta=\infty\right\}.
\]
Then $\dimH\mathscr Z_\delta=0$ for every $\delta>0$.
\end{theorem}
\begin{proof}
Fix $s>0$. Write $L(w)=\sum_{i=1}^k n_i^\delta$ for a word
$w=(n_1,\ldots,n_k)$ and let
$Z_s=C^s\sum_{n\ge1}e^{-csn^\delta/2}<\infty$.
For a fixed $H>0$, the length-$k$ cylinders having $L(w)\ge Hk$
satisfy
\[
 \sum_{\substack{|w|=k\\L(w)\ge Hk}}|Y_w|^s
 \le |Y|^s e^{-csHk/2}Z_s^k.
\]
Choose $H$ so large that $e^{-csH/2}Z_s<1$.
Consequently the sum of these $s$-contents over all $k$ is finite.
The set of points belonging to infinitely many such cylinder families
has zero Hausdorff $s$-measure by the covering definition.
Every point of $\mathscr Z_\delta$ belongs to infinitely many families
with $L(w)\ge Hk$. Since $s>0$ was arbitrary, its Hausdorff dimension
is zero. No assumption on the relative location of the branches is
needed for this upper bound.
\end{proof}

\begin{remark}\label{rem:generalpenalty}
The energy argument is more general than a stretched-exponential
assumption. If $\|\psi_n'\|_\infty\le C e^{-\Phi(n)}$ with
$\Phi(n)\ge0$ and
$\sum_n e^{-s\Phi(n)/2}<\infty$ for every $s>0$, then
\[
 \dim_H\left\{y:\limsup_{k\to\infty}
 \frac1k\sum_{i=1}^k\Phi(\tau_i(y))=\infty\right\}=0.
\]
Indeed, repeat the proof of Theorem~\ref{thm:energy} with
$L(w)=\sum_i\Phi(n_i)$. Thus the zero-dimensionality mechanism is an
abstract \emph{summable inverse-branch energy principle}, not a feature
of one particular rate.
\end{remark}

\begin{theorem}\label{thm:barrier}
Assume \eqref{eq:stretchedbound}.
\begin{enumerate}[label=\textup{(\roman*)},leftmargin=*]
\item If $0<\delta\le1$, then
\begin{equation}\label{eq:secondorderbarrier}
\dimH\left\{y\in K_\infty:
\limsup_{k\to\infty}\frac{t_k(y)^\delta}{k}=\infty\right\}=0.
\end{equation}
In particular, for every $0\le b<\delta$,
$ \dimH\bigl\{y\in K_\infty:\upD_1(R(y,Y))\le b\bigr\}=0. $
\item If $\delta\ge1$, the infinite-return points having zero
physical-time visit frequency satisfy
\[ \dimH\left\{y\in K_\infty:
 \lim_{R\to\infty}\frac{\#(R(y,Y)\cap[1,R])}{R}=0\right\}=0. \]
\end{enumerate}
\end{theorem}
\begin{proof}
For $0<\delta\le1$ concavity yields
$\sum_{i\le k}\tau_i^\delta\ge(\sum_{i\le k}\tau_i)^\delta=t_k^\delta$.
The first conclusion follows from Theorem~\ref{thm:energy}.
If $\upD_1(R)\le b<\delta$, Li and Xi~\cite[Proposition~2.4]{LX} show that,
for every $b'$ strictly between $b$ and $\delta$, one has
$k\le t_k^{b'}$ for all large $k$. Thus
$t_k^\delta/k\ge k^{\delta/b'-1}\to\infty$.
For $\delta\ge1$, Jensen's inequality gives
$k^{-1}\sum_{i\le k}\tau_i^\delta\ge(t_k/k)^\delta$.
Zero visit frequency is equivalent to $k/t_k\to0$, so the
energy-average diverges. Apply Theorem~\ref{thm:energy} again.
\end{proof}

\begin{remark}\label{rem:criticalcount}
For $0<\delta<1$, the condition
$\#(R(y,Y)\cap[1,R])=o(R^\delta)$ implies
$t_k^\delta/k\to\infty$ and hence forces dimension zero.
The unbounded-$\limsup$ condition in \eqref{eq:secondorderbarrier}
is stronger: even recurrent subsequences of exceptionally sparse
returns already yield a zero-dimensional exceptional set.
This is a \emph{second-order} obstruction at the counting
exponent $\delta$; the equality of exponents alone does not imply it.
\end{remark}

Recall the bounded-gap cluster spectrum $\varphi_a$ defined by \eqref{yingshe}. It describes clusters with $\asymp2^{aj}$
visits in intervals of length $\asymp2^{aj}$ near physical time $2^j$.

The critical line is settled by allowing a slowly diverging number of
free visits beyond the leading power. This correction is invisible to
all macroscopic intermediate dimensions but decisive for Euclidean
Hausdorff dimension.
\begin{theorem}\label{thm:clusterphase}
Assume the guarded full-branch conditions \textup{(G1), (G2), (R)},
with $K$ the induced limit set, and the two-sided geometric bounds
\begin{equation}\label{eq:twosidedstretched}
 c_{\mathrm g}n^\delta-C_{\mathrm g}\le -\log\|\psi_n'\|_\infty,
 \qquad -\log\inf_Y|\psi_n'|\le C_{\mathrm g}(1+n^\delta),
 \qquad n\ge1,
\end{equation}
where $0<\delta<1$ and $c_{\mathrm g},C_{\mathrm g}>0$ are fixed
geometric constants, independent of any prescribed counting coefficient.
For $0<a<1$ set
\[
\mathscr E(\varphi_a)=\left\{y\in K_\infty:p(y)=0,\quad
\lowD_\theta(R(y,Y))=\upD_\theta(R(y,Y))=\varphi_a(\theta)
\quad(0\le\theta\le1)\right\}.
\]
Then
\[ \dimH\mathscr E(\varphi_a)=
\begin{cases}0,&0<a<\delta,\\ \dim_H K,&\delta\le a<1.
\end{cases} \]
\end{theorem}
\begin{proof}
If $a<\delta$, the endpoint of the spectrum is
$\varphi_a(1)=a<\delta$; Theorem~\ref{thm:barrier} gives zero dimension.

To prove full dimension, fix $M\ge c$ and $\epsilon>0$ and use the guarded free-block alphabet
from Lemma~\ref{lem:rigorous-forcing}, writing $L$ for its common
block length.  The lemma needs only the uniform insertion-energy
condition.  Take positive factors $L_j\ge1$ such that
$\log L_j=o(j)$ and $2^{aj}L_j=o(2^j)$; for example
$L_j=j^3$ at equality $a=\delta$, and $L_j=1$ for $a>\delta$.
During physical stage $[2^j,2^{j+1})$ make
\begin{equation}\label{eq:criticalblock}
K_j=L\left\lfloor\frac{2^{aj}L_j}{L}\right\rfloor
\end{equation}
consecutive freely chosen bounded roof symbols in guarded blocks.
Then insert one prescribed bridge
$B_j=2^j-\sum_{\text{free stage }j}\tau_i$,
which is positive, uniformly comparable to $2^j$, and brings the
itinerary exactly to $2^{j+1}$. Make an initial prescribed jump to
the starting time $2^{j_0}$. The return times within each stage are
bounded-gap clusters of $K_j=2^{aj+o(j)}$ points in an interval of
length $2^{aj+o(j)}$.

Their common macroscopic dimension spectrum is $\varphi_a$: to see this directly,
replace every $2^{aj}$ by $2^{aj+o(j)}$ in the proof of
Lemma~\ref{lem:clusters}. The upper estimates acquire factors
$2^{o(j)}$, while the lower bound on the occupied shell
$j=N+O(1)$ becomes
$2^{N[a-1+\theta(1-s)]+o(N)}$.
At $\theta=0$, the single-cluster covering contributions are
$2^{-(1-a)sj+o(j)}$, summable for every $s>0$.
The visit frequency vanishes since $a<1$.

By \eqref{eq:twosidedstretched}, the geometric penalty of the new
bridge is at most $C2^{\delta j}$, while the number of free symbols
in a completed stage is $K_j\asymp 2^{aj}L_j$.
If $a>\delta$, choose $L_j\equiv1$; if $a=\delta$, take
$L_j=j^3$. In both cases,
$
\frac{2^{\delta j}}{K_{j-1}}
 \lesssim \frac{2^{(\delta-a)j}}{L_{j-1}}\rightarrow0.
$
The same inequality, together with the geometric growth of the numbers
of free letters in successive stages, bounds the total geometric
penalty divided by the number of free letters at \emph{every} intermediate prefix, including the
prefix immediately after a long bridge. Thus the hypothesis
\eqref{eq:prefixpenalty} holds for the inserted-bridge construction.
The guarded deletion argument yields a set of initial points of
Hausdorff dimension at least $\dim_HK_M-\epsilon$. Letting
$\epsilon\downarrow0$ and $M\to\infty$
gives $\dim_H K$, which is also the trivial upper bound.
\end{proof}

\begin{corollary}\label{cor:criticalsecondorder}
Under the assumptions of Theorem~\ref{thm:clusterphase}, at $a=\delta$
the full-dimensional realization may be chosen so that, uniformly
along its itineraries,
$\#(R(y,Y)\cap[1,R])\asymp R^\delta(\log R)^3$ as $R\to\infty.$
In contrast, under the upper contraction assumption alone the set
of itineraries satisfying
$\#(R(y,Y)\cap[1,R])=o(R^\delta)$
has Hausdorff dimension zero.
\end{corollary}
\begin{proof}
Choose $L_j=j^3$ in \eqref{eq:criticalblock}. The count of
visits completed in stages before $j$ is comparable to
$2^{\delta j}j^3$, uniformly in free choices.
Since $R\in[2^j,2^{j+1})$ has $R^\delta(\log R)^3$
within a fixed multiplicative factor of $2^{\delta j}j^3$,
the first assertion holds even during the long bridge. The last
assertion is Remark~\ref{rem:criticalcount}.
\end{proof}

\subsection{Polynomial packets and fractal micro-packets}
The preceding estimates also explain why the original polynomial packets cannot retain dimension in the stretched-exponential regime.
\begin{proposition}\label{prop:packetzero}
Assume \eqref{eq:stretchedbound} with $0<\delta\le1$, and fix
$0<b<a<1$. Consider a dynamical first-return realization of the
two-scale packet construction introduced before
Lemma~\ref{lem:packets}. At stage $j$, the packets have lengths
$\ell_j=j^{O(1)}$, there are $Q_j\asymp 2^{bj}/\ell_j$
packets filled with freely chosen short-return symbols, and
consecutive packets are separated by prescribed first-return
branches of length $G_j\asymp 2^{(a-b)j}\ell_j$.

If the free symbols range over a fixed finite alphabet, then the
resulting Cantor set of initial points has Hausdorff dimension zero.
\end{proposition}
\begin{proof}
At the end of stage $j$, the number of freely varied symbols is
$O(2^{bj})$, so the number of available prefix words is at most
$M^{O(2^{bj})}$ for some fixed integer $M$. The prescribed packet
separators alone give a logarithmic contraction of at least
$
cQ_jG_j^\delta
\asymp c\,2^{[b+\delta(a-b)]j}\ell_j^{\delta-1}.
$
For fixed $\delta>0$ and polynomial $\ell_j$, this dominates
$2^{bj}$ by an exponentially growing factor. Consequently, for
every $s>0$, the sum of the $s$-powers of the stage-cylinder
diameters is bounded by
\[
\exp\!\left(
C_s2^{bj}-c_s2^{[b+\delta(a-b)]j}\ell_j^{\delta-1}
\right)\longrightarrow0.
\]
These cylinders cover the constructed set, which therefore has
Hausdorff dimension zero.
\end{proof}

\begin{remark}
Theorem~\ref{thm:attain} provides a realization scheme for the \emph{integer-set} spectra defined by the macroscopic covering functionals, and for intermittent first-return systems with logarithmic-power derivative penalty. Theorem~\ref{thm:barrier} is a universal dynamical impossibility result in the stretched-exponential class. Proposition~\ref{prop:packetzero} rules out the original polynomial-packet construction in the stretched-exponential class, whereas Theorem~\ref{thm:twosharp} proves that a balanced Cantor micro-packet arrangement achieves full Hausdorff dimension exactly when $b\ge\delta$. Thus the geometry inside the packets, and not merely the macroscopic spectrum, decides whether a particular construction is efficient.
\end{remark}

The polynomial-size packets in Proposition~\ref{prop:packetzero}
are geometrically inefficient when a long branch has penalty $n^\delta$.
The following alternative replaces each dense packet by a balanced
Cantor arrangement of freely variable short-return blocks.  It retains
exactly the same macroscopic dimension spectrum while lowering the sum of the
$\delta$-powers of its prescribed gaps.  The elementary discretization
is recorded separately to make the critical case $b=\delta$ explicit.

\begin{lemma}\label{lem:balancedseeds}
For every $0<b<1$ and integer $m\ge1$ there is a set
$C_{b,m}=\{c_1,\ldots,c_{N_m}\}\subset\{0,\ldots,2^m-1\}$,
where $c_1<\cdots<c_{N_m}$, such that, with constants
independent of $m$, $c_b 2^{bm}\le \,\#C_{b,m}\le C_b2^{bm},$ and
\begin{align}
 &\#(C_{b,m}\cap I)\le C_b(1+|I|)^b
 \quad\text{for every interval }I,\label{eq:seedlocal}\\
 &\sum_{i=1}^{N_m-1}(c_{i+1}-c_i)^\delta
 \le C\,2^{bm}\sum_{k=0}^m 2^{(\delta-b)k},
 \qquad\delta>0,\label{eq:seedenergy}
\end{align}
Moreover the diameter of $C_{b,m}$ is comparable to $2^m$.
\end{lemma}
\begin{proof}
Let $A_b=\{i\ge0:\lceil b(i+1)\rceil-\lceil bi\rceil=1\}$
and take
\[
 C_{b,m}=\left\{\sum_{i\in A_b\cap[0,m)}\varepsilon_i2^i:
 \varepsilon_i\in\{0,1\}\right\}.
\]
Every string of $q$ successive binary positions contains $bq+O(1)$
active positions, with an absolute discrepancy of at most two.
Thus $\#C_{b,m}=2^{bm+O(1)}$.  A dyadic interval of length $2^k$
contains at most $2^{bk+O(1)}$ members, and an arbitrary interval of
length $r$ meets at most two or three dyadic intervals of length
$2^{\lceil\log_2(1+r)\rceil}$, proving \eqref{eq:seedlocal}.
At most $C2^{b(m-k)}$ dyadic intervals of length $2^k$ meet
$C_{b,m}$. Consequently at most $C2^{b(m-k)}$ consecutive gaps have
length at least $2^k$: each such gap separates occupied dyadic
intervals, apart from at most one boundary gap per occupied interval.
Summing the $\delta$-powers of gaps by dyadic size gives
\[
 \sum_{i=1}^{N_m-1}(c_{i+1}-c_i)^\delta
 \le C\sum_{k=0}^m2^{b(m-k)}2^{\delta k}
 \le C2^{bm}\sum_{k=0}^m2^{(\delta-b)k},
\]
which is \eqref{eq:seedenergy}.  In particular, the right side is
$O(2^{bm})$ for fixed $\delta<b$, $O(m2^{bm})$ when
$\delta=b$, and $O(2^{\delta m})$ for fixed $\delta>b$.
The unsimplified geometric sum is uniform even when $b$ varies
through $\delta$, since it is always at most
$(m+1)2^{\max\{b,\delta\}m}$.  Finally, the largest active binary
position lies within a bounded distance of $m-1$, so the largest
member of $C_{b,m}$ is bounded below by $c_b2^m$.
\end{proof}

The next estimate is uniform when the micro-packet parameters
vary from one physical shell to another.  This is the ingredient
needed for the general-spectrum diagonalization at the critical
counting endpoint.
\begin{lemma}\label{lem:uniform-micropenalty}
	Fix $0<\delta<1$, let $R=2^j$, $0<b<a\le1$, and put
	$e=(a-b)/(1-b)$. Suppose a prescribed return of length $n$
	contributes at most $C(1+n^\delta)$ to logarithmic contraction.
	Consider $Q\lesssim R^a/W$ consecutive cells of width
	$W=d2^m\asymp R^e$, where $1\le d\le W\le R$, with
	free-block starts given by translates of $dC_{b,m}$.
	Assume that intervening returns are positive and that there are
	$O(1)$ additional bridges of length $O(R)$; for $a=1$ require
	$QW\le\eta R$ with fixed $\eta<1$.
	Then the total inserted contraction contribution $A_j$ satisfies
	$
	A_j\lesssim (1+j)R^{\max\{b,\delta\}}.
	$
	The same bound holds for every stage prefix. All implicit
	constants are uniform in $a,b,j$ and the free choices.
\end{lemma}

\begin{proof}
	By Lemma~\ref{lem:balancedseeds}, the within-cell gaps have
	total contribution at most
	$CQd^\delta2^{bm}S_m(b,\delta)$, where
	\[
	 S_m(b,\delta)=\sum_{\ell=0}^m2^{(\delta-b)\ell}
	 \le (m+1)\max\{1,2^{(\delta-b)m}\}.
	\]
	The factor $2^{(\delta-b)m}$ must be retained when $b<\delta$.
	The cell-to-cell gaps contribute at most $CQW^\delta$,
	while the bounded number of stage bridges contribute $CR^\delta$.
	Hence
	$
	A_j\le C\{Qd^\delta2^{bm}S_m(b,\delta)+QW^\delta+R^\delta\}.
	$
	Since $W=d2^m$, this is bounded by
	$C(m+1)(Qd^\delta2^{bm}+QW^\delta)+CR^\delta$.
	Since $W=d2^m\asymp R^e$, $d\le R$, and $1+m\le C(1+j)$, we have 
	$
	Qd^\delta2^{bm}
	\le C R^{a-e}d^{\delta-b}R^{eb}
	=C R^b d^{\delta-b}
	\le C R^{\max\{b,\delta\}}.
	$
	Furthermore,
	\[
	QW^\delta\le C R^{a-e(1-\delta)},\qquad
	a-e(1-\delta)
	=\frac{1-a}{1-b}\,b
	+\frac{a-b}{1-b}\,\delta
	\le\max\{b,\delta\}.
	\]
	The same estimates apply to any prefix because every summand
	is nonnegative. The proof uses the unsimplified geometric sum
	in Lemma~\ref{lem:balancedseeds}; therefore no factor
	$|b-\delta|^{-1}$ is introduced as $b\to\delta$.
\end{proof}

The vanishing endpoint requires a stronger diagonal choice than the
usual condition that type-dependent constants are $2^{o(j)}$.
The following elementary estimate records the required uniformity.
\begin{lemma}\label{lem:zeroendpoint-diagonal}
Let $R_j=2^j$, and in shell $[R_j,2R_j)$ let $E_j$ consist of a
baseline of diameter at most $C R_j^b j^P$, with fixed $b<1$,
and at most one additional packet.  Assume that each nonconstant
packet has parameters $0<b_i<a_i<1$, diameter at most
$C_i j^{C_0}R_j^{a_i}$, and that at the shell where type $i$ is active
\begin{equation}\label{eq:zeroendpoint-epoch}
 j(1-a_i)\ge j^{1/2},\qquad
 \log_2(1+C_i)\le j^{1/4}.
\end{equation}
A constant type $a_i=1$ is permitted provided its number of
points is at most $C_i j^{C_0}R_j^{b_i}$, where $b_i\le c<1$,
and the same bound on $C_i$ holds.  Then, for every $s>c$ with
$0<s<1$, the covering sum over shells $n\ge N$ at $\theta=0$
of the baseline and all these packets tends to zero as the outer
index $N\to\infty$, provided $b<1$.  In particular, this covers the baseline and
constant support types of Theorem~\ref{thm:generalphase} with
$c=h(0)$; no inequality $b\le c$ is needed.
\end{lemma}
\begin{proof}
Under the uniform epoch conditions \eqref{eq:zeroendpoint-epoch},
at $\theta=0$ the permitted radii in shell $j+1$ range from $1$
to a quantity comparable with $R_j$, independently of the outer
index $N$.  One interval covers a nonconstant packet, so its
$s$-contribution is at most
$
 C_s(1+C_i)^s j^{sC_0}2^{-sj(1-a_i)}
 \le C_s j^{sC_0}2^{s j^{1/4}-s\sqrt j}.
$
This bound is summable in $j$, independently of the active type.
For a constant type use unit intervals instead; its contribution
is at most $C_s(1+C_i)j^{C_0}2^{-(s-b_i)j}$, which is summable
uniformly since $s>c\ge b_i$.  A baseline of diameter at most
$CR_j^b j^P$ contributes at most
$C_sj^{Ps}2^{-s(1-b)j}$ if $b<1$; a subpower baseline is even
smaller.  Summing the respective estimates over $j\ge N-O(1)$
gives a quantity tending to zero.  No comparison with a selected
subsequence of shells is used.
\end{proof}

\begin{theorem}\label{thm:twosharp}
Assume the guarded full-branch conditions \textup{(G1), (G2), (R)}
and suppose that
\eqref{eq:twosidedstretched} holds with $0<\delta<1$.
For $0<b<a<1$ define
\[
 \mathscr E(P_{a,b})=\bigl\{y\in K_\infty:
 p(y)=0,\ \lowD_\theta(R(y,Y))=\upD_\theta(R(y,Y))
 =P_{a,b}(\theta)\ \text{for all }0\le\theta\le1\bigr\}.
\]
Then the following \emph{sharp two-parameter transition} holds:
\begin{equation}\label{eq:twosharp}
 \dimH\mathscr E(P_{a,b})=
 \begin{cases}
 0,&0<b<\delta,\\ \dim_H K,&\delta\le b<a<1.
 \end{cases}
\end{equation}
In particular the transition is governed by the counting endpoint
$P_{a,b}(1)=b$, independently of the spread exponent $a$.
\end{theorem}
\begin{proof}
When $b<\delta$, Li and Xi~\cite[Proposition~2.4]{LX} and
Theorem~\ref{thm:barrier} give dimension zero. We prove full
dimension for $\delta\le b$ by constructing a nonstationary
family of genuine first-return itineraries.

Put
\[
 e=\frac{a-b}{1-b},\qquad 0<e<a<1,
 \qquad R_j=2^j.
\]
Fix $M\ge c$ and $\epsilon>0$, and choose a guarded-block alphabet
from Lemma~\ref{lem:rigorous-forcing} with dimension loss at most
$\epsilon$ and all roof values bounded by $M$. Choose $P>2/(1-\delta)$ and let $H_j$ be a multiple of the
fixed guarded-block length such that $H_j\asymp j^P$.
Choose a large fixed $A>4M$ and put
\[
 d_j=AH_j,\quad
 m_j=\left\lfloor ej-\log_2d_j\right\rfloor,\quad
 W_j=d_j2^{m_j}\asymp R_j^e,
 \quad Q_j=\left\lfloor \eta_0R_j^a/W_j\right\rfloor,
\]
where $\eta_0>0$ is a sufficiently small fixed packing constant.
In each of the $Q_j$
consecutive cells of length $W_j$, place the prospective starting
times of free blocks at
$t_{j,q,v}=R_j+qW_j+d_jv$ with $ 0\le q<Q_j$ and $ v\in C_{b,m_j}.$
At each such starting time run $H_j$ bounded-roof free symbols,
arranged as complete guarded blocks.  Between this free block and
the next prescribed starting time insert one return symbol of
exactly the required gap length.  This symbol is positive because
the entire free block occupies at most $MH_j<d_j/4$ units of
physical time, whereas distinct starting times are separated by at
least $d_j$.  Its value may depend on the previous free choices;
that dependence is causal, as required by Lemma~\ref{lem:rigorous-forcing}.
After the last free block of stage $j$, one additional prescribed
bridge brings the time to exactly $2^{j+1}$.  Since the entire
configuration has width $O(R_j^a)=o(R_j)$, the bridge length is
comparable to $R_j$.  This procedure realizes actual first-return
itineraries, not merely a subset of integer times.

Let $E_j$ be the resulting set of physical return times in
$[R_j,2R_j)$, and put $E=\bigcup_jE_j$.
The number of returns in stage $j$ is, uniformly in the free choices,
\begin{equation}\label{eq:microcount}
 K_j:=\#E_j\asymp Q_jH_j2^{bm_j}
 \asymp R_j^b H_j^{1-b}=R_j^{b+o(1)}.
\end{equation}
Their total span is $O(R_j^a)$, and their occupied cells extend
across an interval of length comparable to $R_j^a$.
For any interval $I$ of radius $r$ and the normalized counting
measure $\nu_j$ on $E_j$, the balanced-set estimate
\eqref{eq:seedlocal} and the bounded duration of each free block
imply
\begin{align}
 \nu_j(I)&\le R_j^{o(1)}(r/R_j)^b,
       &&1\le r\le W_j,\label{eq:microlocal1}\\
 \nu_j(I)&\le C r/R_j^a,
       &&W_j\le r\le R_j,\label{eq:microlocal2}
\end{align}
with an $o(1)$ uniform over every allowed choice of free letters.
Indeed, for $r\le W_j$ the interval meets only a bounded number
of cells and at most
$C H_j(1+(r/d_j)^b)$ return times; dividing by
\eqref{eq:microcount} gives \eqref{eq:microlocal1}.
For $r\ge W_j$ it meets $O(1+r/W_j)$ cells, each containing
$O(H_j2^{bm_j})$ returns, and
$Q_jW_j\asymp R_j^a$, giving \eqref{eq:microlocal2}.

We check the lower bound for the covering sum over shells
$n\ge N$ for every fixed $0<\theta\le1$ and
$0<s<P_{a,b}(\theta)$, so $s<b$ and
$1-a-\theta(1-s)<0$. Take the occupied shell $S_{j+1}$ with
$j=N$ and permitted radius $1\le r\le U=2^{j+1-N\theta}$.
For $r\le W_j$, \eqref{eq:microlocal1} yields
\[
 \frac{\nu_j(B(x,r))}{(r/R_j)^s}
 \le R_j^{o(1)}(W_j/R_j)^{b-s}
 =R_j^{-(1-e)(b-s)+o(1)}.
\]
For $W_j<r\le U$ (when this range is nonempty),
\eqref{eq:microlocal2} yields
\[
 \frac{\nu_j(B(x,r))}{(r/R_j)^s}
 \le C R_j^{s-a}U^{1-s}
 \le R_j^{1-a-\theta(1-s)+o(1)}.
\]
Both quantities decrease exponentially in $j$.
The mass-distribution argument from Lemma~\ref{lem:packets}
therefore gives $\cH_N^{s,\theta}(E)\ge2^{c_sN}$.
For the upper bound, unit balls and the containing interval give,
respectively, $K_jR_j^{-s}=R_j^{b-s+o(1)}$ and the same
maximal-radius bound as \eqref{eq:wide}--\eqref{eq:one}.
Summing over $j\ge N-O(1)$ proves
$
 \lowD_\theta(E)=\upD_\theta(E)=P_{a,b}(\theta)
 ~~(0<\theta\le1).
$
At $\theta=0$ the containing interval yields dimension zero;
\eqref{eq:microcount} and $b<1$ also imply zero natural density.

The loss of a prescribed branch of roof value $n$ is at most
$C(1+n^\delta)$ by \eqref{eq:twosidedstretched}.
Within one cell, successive seed gaps have physical length at most
$d_j$ times their integer difference.  Lemma~\ref{lem:balancedseeds}
shows that the sum of prescribed penalties within all cells is at most
$
 C Q_j d_j^\delta 2^{bm_j}(1+m_j).
$
The cell-to-cell gaps have length $O(W_j)$ and the final bridge
has length $O(R_j)$, so the total penalty $A_j$ incurred in stage $j$
obeys
\begin{equation}\label{eq:micropenalty}
 A_j\le C\bigl[
 Q_j d_j^\delta2^{bm_j}(1+m_j)
 +Q_jW_j^\delta+R_j^\delta\bigr].
\end{equation}
Comparing \eqref{eq:micropenalty} with \eqref{eq:microcount}, and using
$W_j=d_j2^{m_j}$, yields
\[ \frac{A_j}{K_j}
 \le C\left[
 d_j^{\delta-1}(1+m_j)
 +d_j^{\delta-1}2^{(\delta-b)m_j}
 +\frac{R_j^{\delta-b}}{H_j^{1-b}}
 \right]\longrightarrow0. \]
The last limit holds also at $b=\delta$: the middle factor is then
$d_j^{b-1}$, and the final factor is $H_j^{-(1-b)}$.
The chosen $P>2/(1-\delta)$ absorbs the factor $1+m_j=O(j)$.
Since $K_{j-1}\asymp 2^{b(j-1)}H_{j-1}^{1-b}$ is comparable to
$K_j$ up to constants and $A_j=o(K_j)$, all penalties accumulated
through any partial stage $j$ are $o$ of the number of free letters
already accumulated before that stage.  The same is true over the
previous stages by summing the geometric progression in
\eqref{eq:microcount}.  Thus the full prefix condition
\eqref{eq:prefixpenalty} holds uniformly in the free choices,
including immediately after the final bridge.

Lemma~\ref{lem:rigorous-forcing} now gives Hausdorff dimension at
least $\dim_HK_M-\epsilon$ for the initial points realizing the
constructed spectrum.  Letting $\epsilon\downarrow0$ and $M\to\infty$ gives $\dim_H K$ by \ref{hyp:G2}; the
opposite inequality is immediate since the level set is contained
in $K$.  This proves \eqref{eq:twosharp}.
\end{proof}

\subsection{A general-spectrum phase transition}

The sharp two-parameter result admits a useful strengthening: the
transition is determined by the counting endpoint even for a general
admissible macroscopic dimension spectrum.  The proof combines the supporting-curve
decomposition of Lemma~\ref{lem:support}, the balanced micro-packets,
and a dense free baseline.  The baseline is permitted because it has a
spectrum lying below every scale-admissible function with the same
value at $\theta=1$.

\begin{theorem}\label{thm:generalphase}
Under the assumptions of Theorem~\ref{thm:twosharp}, let
$h:[0,1]\to[0,1]$ be any scale-admissible function, and put
$b=h(1)$. Define $\mathscr E(h)$ as the set of points with zero
physical-time visit frequency and
$\lowD_\theta(R)=\upD_\theta(R)=h(\theta)$ for every
$\theta\in[0,1]$. Then
\begin{equation}\label{eq:generalphase}
 \dimH\mathscr E(h)=
 \begin{cases}
 0,&0\le b<\delta,\\
 \dim_H K,&\delta\le b\le1.
 \end{cases}
\end{equation}
This includes discontinuous spectra at $\theta=0$, and the boundary
cases $b=\delta$ and $b=1$.
\end{theorem}
\begin{proof}
\noindent\textit{Step 1. Energy obstruction and the free baseline.}
When $b<\delta$, Theorem~\ref{thm:barrier} applies because the
upper counting dimension is $b$ by Li and Xi~\cite[Proposition~2.4]{LX}.
We construct full-dimensional realizations for $b\ge\delta$.
Fix $M\ge c$ and $\epsilon>0$, and use guarded blocks with
dimension loss at most $\epsilon$ in Lemma~\ref{lem:rigorous-forcing}.

First take $\delta\le b<1$.  The scale comparison with
$\theta_2=1$ gives, for every $0<\theta\le1$,
\begin{equation}\label{eq:baselinebeneath}
 \varphi_b(\theta)\le h(\theta)
 \quad(0<\theta\le1),\qquad
 \varphi_b(0)=0\le h(0).
\end{equation}
In every physical shell $[R_j,2R_j)$ begin with a baseline of
$
 B_j=\lfloor R_j^b j^{P}\rfloor
$
freely variable bounded-roof letters in complete guarded blocks,
where $P>4+2/(1-\delta)$ is fixed.  This block fits in a time
interval of length $O(B_j)=o(R_j)$; its complete return-time spectrum
is $\varphi_b$ by the one-cluster argument of
Theorem~\ref{thm:clusterphase}.  Thus, by
\eqref{eq:baselinebeneath}, the baseline does not raise the
desired dimension at any parameter.

Enumerate the supporting curves $P_{a_i,b_i}$ from
Lemma~\ref{lem:support}; every $b_i\le b$.  If $h(0)>0$, also
include the constant supporting type $(a_i,b_i)=(1,h(0))$.
Allocate one active type per shell in the explicit cyclic epochs of
the proof of Theorem~\ref{thm:attain}, choosing $J_k\ge k^3$ and slowing the epochs so that, for every
active nonconstant type, with $e_i=(a_i-b_i)/(1-b_i)$,
\begin{equation}\label{eq:microepochs}
 j\min\{b_i,1-b_i,1-a_i,e_i,a_i-e_i\}\ge\sqrt j.
\end{equation}
For the count, mass, and diameter bounds, choose one constant $C_i\ge1$
for each fixed supporting type $i$, including the constant type.
Enlarge each epoch start $J_k$ further so that
\begin{equation}\label{eq:improved-micro-epochs}
 \log_2(1+C_i)\le j^{1/4}
 \qquad (i\le k,\quad j\ge J_k).
\end{equation}
The stronger $j^{1/4}$ bound, combined with
$j(1-a_i)\ge\sqrt j$ in \eqref{eq:microepochs}, is essential
at $\theta=0$; $2^{o(j)}$ alone would not be sufficient when
$a_i\uparrow1$.  Each $C_i$ is fixed by its type; \eqref{eq:improved-micro-epochs}
absorbs this dependence into the shell index.  Each epoch constrains
only finitely many types, so its start can be enlarged while
maintaining $J_{k+1}-J_k\ge k^2$ and cyclic recurrence.  Lemma~\ref{lem:epoch-recurrence} therefore still
provides every fixed type at $j=N+o(N)$ after all large $N$.  The quantities $e_i$ and
$a_i-e_i$ in \eqref{eq:microepochs} ensure that both the cell depth
and the integer number of cells tend to infinity; conditions only
on $a_i-b_i$ would not suffice.  Following the initial baseline, allocate an
interval of length at most $cR_j^{a_i}$ to a micro-packet of type
$(a_i,b_i)$, using the construction of
Theorem~\ref{thm:twosharp} with $H_j\asymp j^P$.
When $a_i=1$, use one macroscopic cell of length
$W_j\asymp cR_j$ by taking
$m_j=\lfloor\log_2(cR_j/d_j)\rfloor$; do not use the
small packing constant in a formula that would make the integer
number of cells vanish.  A prescribed jump from the baseline to the
first micro-packet block and a final bridge to $2R_j$ are positive
and at most $O(R_j)$, provided the fixed packing constant $c$ is
small.  The first bridging jump is allowed to depend on the preceding
free letters; no noncausal choice is needed.

\medskip
\noindent\textit{Step 2. Control of the accumulated contraction.}
Lemma~\ref{lem:uniform-micropenalty} applies to every active
supporting type, including the constant type $a_i=1$.  Its bound
comes from the unsimplified finite geometric sum and therefore does
not acquire a factor such as $|b_i-\delta|^{-1}$ when $b_i$
approaches $\delta$.  Since
$b_i\le b$ and $\delta\le b$, it gives the uniform stage bound
$
 A_{j,i}\le C(1+j)R_j^{\max\{b_i,\delta\}}
           \le C(1+j)R_j^b.
$
The baseline in the preceding shell contains
$B_{j-1}\asymp R_{j-1}^b(j-1)^P$ free letters.  Thus
$
 \frac{A_{j,i}}{B_{j-1}}\le Cj^{1-P}\rightarrow0.
$
This bound holds also immediately after any inserted bridge,
because the lemma bounds all prefixes.  In addition,
$\sum_{i\le j}A_i\le C(1+j)R_j^b$ by the geometric growth of
$R_i=2^i$, whereas the preceding baseline blocks already supply
$\sum_{i<j}B_i\gtrsim R_{j-1}^b(j-1)^P$ free letters.
The ratio is $O(j^{1-P})=o(1)$ uniformly in the active types
and the free choices.  Hence the global prefix condition
\eqref{eq:prefixpenalty} follows without imposing constants
that deteriorate as $b_i\to\delta$.

\medskip
\noindent\textit{Step 3. Verification of the prescribed spectrum.}
We verify the lower covering estimate, including the constant
support type, without appealing only to a fixed-parameter
subsequence.  Fix $0<\theta\le1$ and $0<s<h(\theta)$.  If the
baseline has $\varphi_b(\theta)>s$, its one-cluster measure gives the
required lower bound in every shell.  Otherwise, if
$s<h(0)$, the recurring constant type $a_i=1$, $b_i=h(0)$
may be used.  For this fixed type, put $c=h(0)$ and let $\nu_j$
be normalized counting measure on its single balanced cell.
The seed estimate and the bounded-roof free blocks yield
\[
 \nu_j(B(x,r))\le C\bigl[(H_j/R_j)^c+(r/R_j)^c\bigr],
 \qquad 1\le r\le R_j,
\]
where $H_j\asymp j^P$.  Choose an occurrence of this fixed type
at $j=N+o(N)$ using Lemma~\ref{lem:epoch-recurrence}, and retain
the permitted upper radius $U=2^{j+1-N\theta}$.  For all large $N$
we have $U\le R_j$, and division by $(r/R_j)^s$ gives
\[
 \sup_{\substack{x\in\mathbb R\\1\le r\le U}}
 \frac{\nu_j(B(x,r))}{(r/R_j)^s}
 \le C\left(H_j^cR_j^{s-c}
             +2^{-N\theta(c-s)}\right).
\]
Here the first term uses $r^{-s}\le1$, whereas the second uses
$r\le U$ and $c-s>0$.  Both terms decay exponentially in $N$
for each fixed $\theta>0$.  The case $\theta=0$ is treated below
by a positive shell-content bound and divergence of the full tail.
If neither of these alternatives applies, the supporting-curve
representation supplies a \emph{fixed} nonconstant type
$(a_i,b_i)$ with $P_{a_i,b_i}(\theta)>s$.
By Lemma~\ref{lem:epoch-recurrence}, it occurs for every sufficiently
large outer index $N$ in a shell $j=N+o(N)$.
For this fixed type the local estimates
\eqref{eq:microlocal1}--\eqref{eq:microlocal2}, with
$e_i=(a_i-b_i)/(1-b_i)$ and $U=2^{j+1-N\theta}$, give
\[
 \sup_{\substack{x\in\mathbb R\\1\le r\le U}}
 \frac{\nu_j(B(x,r))}{(r/R_j)^s}
 \le 2^{o(j)}\left(
 2^{-(1-e_i)(b_i-s)j}
 +2^{[1-a_i-\theta(1-s)]j+o(j)}\right).
\]
Both exponents have a fixed negative margin because
$s<P_{a_i,b_i}(\theta)$.  The shift $j-N=o(N)$ is
absorbed in the $o(j)$ term.  In each case the mass-distribution
principle gives a positive (indeed exponentially growing) contribution
to the $s$-covering sum over shells $n\ge N$ for every sufficiently
large $N$.
Conversely, fix $0<\theta\le1$ and
$h(\theta)<s_0<s<1$.  For any active nonconstant type with
$b_i\le s_0$, unit-radius covers give at most
$C_i j^{C_0}2^{-(s-s_0)j}$.  If $b_i>s_0$, then the supporting
inequality $P_{a_i,b_i}(\theta)\le h(\theta)<s_0$
forces
$
 1-a_i\ge\theta(1-h(\theta))>0.
$
When the largest admissible radius $U=2^{j+1-N\theta}$ does
not exceed the packet span, maximal-radius covers give
\[
 \mathcal H_{j+1,N}^{s,\theta}(E_{j,i})
 \le C_i j^{C_0}2^{-(1-a_i)j+N\theta(1-s)}
 \le C_i j^{C_0}2^{-\theta(s-h(\theta))j+O(1)},
\]
using $j\ge N-O(1)$.  If $U$ exceeds the span, a single ball
gives $C_i j^{C_0}2^{-s(1-a_i)j}$, again exponentially small
uniformly in $j$.  The constant type satisfies
$b_i=h(0)\le h(\theta)<s_0$ and so belongs to the first case.
The baseline obeys the same two alternatives, with both counting
and span exponent $b$, using $\varphi_b(\theta)\le h(\theta)$.
The factor $C_i j^{C_0}$ is bounded uniformly over all active
types by $2^{j^{1/4}+O(\log j)}$, by
\eqref{eq:improved-micro-epochs}.  For fixed $\theta>0$ and
$s>h(\theta)$, the preceding three covering cases have decay
$2^{-\gamma j}$ with $\gamma>0$ depending only on
$\theta,s,h(\theta)$, and not on the active type.  Therefore the
full bound is $2^{-\gamma j+j^{1/4}+O(\log j)}$.
Summing over \emph{all} $j\ge N-O(1)$, not just occupied shells,
shows that the sum of $s$-covering contributions over all shells
tends to zero exponentially in $N$.
Thus $\lowD_\theta(E)=\upD_\theta(E)=h(\theta)$ for
all $\theta>0$.  At $\theta=0$, Lemma~\ref{lem:zeroendpoint-diagonal} and
\eqref{eq:improved-micro-epochs} give, for every $s>h(0)$, a
summable uniform upper bound over \emph{all} shell indices, including
types with $a_i\uparrow1$.  In detail, a nonconstant type contributes
at most $C_sj^{C_s}2^{s j^{1/4}-s\sqrt j}$, a constant type with
$b_i=h(0)$ contributes at most
$C_sj^{C_s}2^{j^{1/4}-(s-h(0))j}$, and the baseline contributes
$C_sj^{Ps}2^{-s(1-b)j}$.  These estimates apply to the sum over all shells beyond $N$,
uniformly in $N$.  If $h(0)>0$, the
constant type $a_i=1$, $b_i=h(0)$ recurs at shells $j=N+o(N)$.
For $s<h(0)$, let $H_j\asymp j^P$ be the bounded-roof free
block length and $d_j\asymp H_j$ the seed mesh.  The one-cell
counting estimate from \eqref{eq:seedlocal}, with this \emph{fixed}
constant type, gives
\[
 \nu_j(B(x,r))\le C\bigl((d_j/R_j)^{h(0)}
                    +(r/R_j)^{h(0)}\bigr)
 \le C'(r/R_j)^{h(0)}\quad (r\ge H_j).
\]
For $1\le r\le H_j$, its normalized $s$-ball mass is bounded
by $C H_j^{h(0)}R_j^{-(h(0)-s)}=o(1)$.
Thus all admissible radii have $s$-ball ratio bounded by a fixed
constant, so each occupied shell has $s$-content bounded below by
a positive constant.  There are infinitely many such shells beyond
any outer index $N$, so the sum over all shells $n\ge N$
diverges.  If $h(0)=0$, the upper estimate alone establishes
the endpoint.
This proves the whole prescribed spectrum.  The total number of
visits in shell $j$ is $R_j^{b+o(1)}$, so its physical-time
frequency is zero.  Lemma~\ref{lem:rigorous-forcing} gives dimension at least
$s_M-\epsilon$.  Letting $\epsilon\downarrow0$ and $M\to\infty$
gives $\dim_H\mathscr E(h)=\dim_H K$.

\medskip
\noindent\textit{Step 4. The endpoint $b=1$.}
Scale admissibility forces
$h(\theta)=1$ for every $\theta>0$.  If $h(0)=c<1$, use in every
shell a free baseline consisting of $\lfloor R_j2^{-\sqrt j}\rfloor$
consecutive bounded-roof letters; it has dimension zero at
$\theta=0$ and dimension one at all positive parameters.
If $c>0$, insert on an infinite subsequence of shells the
constant-type micro-packets $(a_i,b_i)=(1,c)$, as above.  Their
forced penalty is at most
$O((1+j)R_j^{\max\{c,\delta\}})=o(R_j2^{-\sqrt j})$,
since $\max\{c,\delta\}<1$.
Hence the deletion lemma still applies, and the spectrum at zero is
exactly $c$.
Finally, if $h(0)=1$, put in every shell
$Q_j\asymp R_j/j^{p+q}$ packets of $\ell_j=j^q$ bounded-roof
free letters, separated by roof values $G_j\asymp j^{p+q}$,
where $p>1$ and $q$ is chosen so large that
$p-(p+q)(1-\delta)<0$.  The number of visits is
$R_j/j^p=o(R_j)$, the packet span is comparable to $R_j$,
and the forced $\delta$-energy is
$O(R_jj^{-(p+q)(1-\delta)}+R_j^\delta)
=o(R_j/j^p)$.
The normalized counting-measure estimate
$\nu_j(B(x,r))\le C(j^{p+q}/R_j+r/R_j)$ shows that the
$\theta=0$ content of every shell is bounded below for every
$s<1$, while unit covers give the upper bound one.
Monotonicity therefore yields the constant-one spectrum.
All these constructions retain Hausdorff dimension $\dim_H K$, proving
\eqref{eq:generalphase}.
\end{proof}

\begin{remark}\label{rem:microvpoly}
The fractal micro-packets explain why Proposition~\ref{prop:packetzero}
should not be interpreted as an impossibility theorem for the
specified two-parameter spectrum.  Although the two constructions
have the same complete macroscopic intermediate-dimension function,
the microscopic placement of their returns makes their forced
$\delta$-energies radically different.  For $b\ge\delta$, the new
construction has penalty $o(K_j)$ and hence retains full Euclidean
Hausdorff dimension.
\end{remark}

\subsection{An affine full-branch tower}
An affine model realizes the general transition and will also serve
as the setting for the energy-level formula in Section~\ref{sec:secondorder}.

Fix $0<\delta<1$ and $c>0$. Set
$
Z_\delta=\sum_{n=1}^\infty e^{-cn^\delta}$ and $
r_n=Z_\delta^{-1}e^{-cn^\delta}.$
Partition $Y=[0,1]$ (up to endpoints) into consecutive intervals
$I_n$ of lengths $r_n$, and let $F|_{I_n}$ be the increasing affine
bijection onto $(0,1)$. Its inverse branches have constant derivatives
$|\psi_n'|=r_n$ and zero distortion. Associate the integer roof
$\tau|_{I_n}=n$, and let $T$ be the first-return tower map on
$
\Delta=\{(y,\ell):y\in Y,\ 0\le\ell<\tau(y)\},
$
which increments $\ell$ until the roof and then applies $F$ and
returns to level zero. The base $Y\times\{0\}$ has the usual
Euclidean geometry of $Y$. We always compute the Hausdorff dimension
of \emph{base coordinates}, not a metric-dependent dimension of the
tower space.
\begin{corollary}\label{cor:lineartower}
For this tower and every scale-admissible $h$,
\[
 \dim_H\mathscr E(h)=
 \begin{cases}0,&h(1)<\delta,\\1,&h(1)\ge\delta.\end{cases}
\]
In particular, this includes the spectra $\varphi_a$ and $P_{a,b}$.
For $h=\varphi_\delta$, a full-dimensional subset can also satisfy
$N_y(t)\asymp t^\delta(\log t)^3$.
\end{corollary}
\begin{proof}
Every integer roof value occurs, the branch images have disjoint
interiors, and any interval $I_m$ with $m\ge2$ is compactly contained
in $(0,1)$. The contraction satisfies
$-\log r_n=cn^\delta+\log Z_\delta$. For each finite subsystem,
the similarity dimension $s_M$ solves
$\sum_{n\le M}r_n^{s_M}=1$; since
$\sum_{n\ge1}r_n=1$, we have $s_M\uparrow1$.
Theorem~\ref{thm:generalphase} and
Corollary~\ref{cor:criticalsecondorder} apply with $\dim_H K=1$.
\end{proof}

\subsection{Exact counting and the critical logarithmic transition}
The first result sharpens the critical energy obstruction under an
exact physical-time counting asymptotic.
\begin{theorem}\label{thm:exactcritical}
Assume the induced inverse branches satisfy the uniform bound
$\|\psi_n'\|_\infty\le Ce^{-cn^\delta}$, where
$0<\delta<1$.  For every $c_0>0$,
\begin{equation}\label{eq:exactcritical}
 \dimH\bigl\{y\in K_\infty:
 N_y(t)\sim c_0t^\delta\quad(t\to\infty)\bigr\}=0.
\end{equation}
More generally, the same conclusion holds for any family of infinite
return itineraries satisfying $t_k^\delta/k\asymp1$ and
$\max_{1\le i\le k}\tau_i=o(t_k)$.
\end{theorem}
\begin{proof}
For an orbit with $N_y(t)\sim c_0t^\delta$, inversion of the
monotone counting function gives
$t_k\sim(k/c_0)^{1/\delta}$.  In particular
$t_{k-1}/t_k\to1$, so
$\tau_k=t_k-t_{k-1}=o(t_k)$.
In fact $M_k:=\max_{1\le i\le k}\tau_i=o(t_k)$:
for any $\varepsilon>0$ choose $i_0$ so that
$\tau_i\le\varepsilon t_i$ for $i\ge i_0$, and note that
$\max_{i<i_0}\tau_i/t_k\to0$.
Since $0<\delta<1$,
$
 \sum_{i=1}^k\tau_i^\delta
 =\sum_{i=1}^k\frac{\tau_i}{\tau_i^{1-\delta}}
 \ge \frac{t_k}{M_k^{1-\delta}}.
$
Consequently,
\[ \frac1k\sum_{i=1}^k\tau_i^\delta
 \ge\frac{t_k^\delta}{k}
 \left(\frac{t_k}{M_k}\right)^{1-\delta}\longrightarrow\infty, \]
because $t_k^\delta/k\to c_0^{-1}$.  The general assertion follows
from the same inequality.  Theorem~\ref{thm:energy} now yields
\eqref{eq:exactcritical}.
\end{proof}

The same argument also excludes finite mean energy, independently
of the branch geometry.
\begin{proposition}\label{prop:exactcount-energy}
If $0<b\le\delta<1$ and $c_0>0$, every infinite return itinerary
with $N_y(t)\sim c_0t^b$ satisfies
\[
 \frac1k\sum_{i=1}^k\tau_i(y)^\delta\longrightarrow\infty.
\]
Thus the joint exact-counting and finite-energy level is empty.
\end{proposition}
\begin{proof}
For $b=\delta$, this is the energy estimate in the proof of
Theorem~\ref{thm:exactcritical}.  If $b<\delta$, then
$t_k\sim(k/c_0)^{1/b}$ and subadditivity gives
\[
 \frac1k\sum_{i=1}^k\tau_i^\delta\ge\frac{t_k^\delta}{k}
 \sim c_0^{-\delta/b}k^{\delta/b-1}\longrightarrow\infty. \qedhere 
\]
\end{proof}

\begin{remark}\label{rem:exact-vs-critical}
Theorem~\ref{thm:exactcritical} does not contradict the full-dimensional
critical realization in Theorem~\ref{thm:clusterphase}.  The latter
allows a slowly diverging correction, such as
$N_y(t)\asymp t^\delta(\log t)^3$, whereas
\eqref{eq:exactcritical} fixes the complete leading asymptotic with a
finite positive constant.  The exponent by itself cannot detect this
difference.
\end{remark}

The next theorem determines the dimension of every exact power-law
counting level.

\begin{theorem}\label{thm:exactpower}
Under the assumptions of Theorem~\ref{thm:twosharp}, for every $c_0>0$
and $0<b<1$ the level set $\{y\in K_\infty:N_y(t)\sim c_0t^b\}$ is
nonempty and
\begin{equation}\label{eq:exactpowerphase}
 \dim_H\{y\in K_\infty:N_y(t)\sim c_0t^b\}
 =\begin{cases}0,&0<b\le\delta,\\ \dim_H K,&\delta<b<1.\end{cases}
\end{equation}
\end{theorem}
\begin{proof}
For $b<\delta$, apply Theorem~\ref{thm:barrier}; for $b=\delta$,
apply Theorem~\ref{thm:exactcritical}.  Nonemptiness in these cases
follows by prescribing successive return times
$t_k=\lfloor(k/c_0)^{1/b}\rfloor$ after a finite initial adjustment:
their differences are eventually positive and every integer roof is
available.

Suppose $b>\delta$.  Choose $q>(1-b)/(b-\delta)$ and
$p=b(q+1)-1$, so that $0<q\delta<p<q$.  Fix $M\ge c$ and
$\epsilon>0$, and select guarded free blocks of common length $L$
from Lemma~\ref{lem:rigorous-forcing}.  At stage $j$, concatenate
$\ell_j=L\lfloor B j^p/L\rfloor$ free letters in complete blocks,
then insert the single prescribed return $n_j=\lfloor A j^q\rfloor$.
Start at a sufficiently large $j$.  Choose $A,B>0$ with
\[
 \frac{B}{p+1}\left(\frac{q+1}{A}\right)^b=c_0.
\]
Uniformly in all free choices, the free-symbol count $K_J$ and elapsed
physical time $T_J$ after stage $J$ satisfy
\[
 K_J\sim\frac{B}{p+1}J^{p+1},\qquad
 T_J=\sum_{j\le J}n_j+O(K_J)\sim\frac{A}{q+1}J^{q+1}.
\]
The $J$ inserted returns are $o(K_J)$.  Since
$\ell_J=o(K_J)$ and $n_J+M\ell_J=o(T_J)$, monotonicity interpolates
these asymptotics to every physical time and gives $N_y(t)\sim c_0t^b$.
The cumulative prescribed contraction is
$O(J^{q\delta+1}+J)=o(K_J)$.  Also $K_{J-1}/K_J\to1$, so this
estimate holds at every intermediate prefix, including immediately
after each bridge.  Lemma~\ref{lem:rigorous-forcing} gives dimension
at least $s_M-\epsilon$.  Let $\epsilon\downarrow0$ and then
$M\to\infty$; the reverse bound is inclusion in $K$.
This completes the proof of \eqref{eq:exactpowerphase}.
\end{proof}

The following theorem shows that, even with the complete macroscopic dimension spectrum fixed, a logarithmic correction to the critical return-time counting law produces a sharp transition between zero and full Hausdorff dimension.
\begin{theorem}\label{thm:logcritical}
	Under the same hypotheses, let $\mathscr C_\beta(c_0)$ be the critical
	level set defined in \eqref{eq:critical-log-level}, where
	$\beta\in\mathbb R$ and $c_0>0$.  Then $\mathscr C_\beta(c_0)$
	is nonempty, and
	\begin{equation}\label{eq:logcriticalphase}
		\dim_H\mathscr C_\beta(c_0)
		=\begin{cases}0,&\beta\le0,\\ \dim_H K,&\beta>0.\end{cases}
	\end{equation}
	The same formula holds if the macroscopic-spectrum constraint is omitted.
	Thus the phase transition persists while both the leading coefficient
	and the entire macroscopic dimension spectrum are fixed.
\end{theorem}

\begin{proof}
If $\beta<0$, the counting function is $o(t^\delta)$, so
Theorem~\ref{thm:barrier} gives dimension zero.  If $\beta=0$,
Theorem~\ref{thm:exactcritical} applies.  This is the genuine critical counting regime: when the positive
logarithmic correction is absent, the exact counting asymptotic forces
the mean $\delta$-energy of successive returns to diverge, as shown
in Theorem~\ref{thm:exactcritical}.  We give a construction for
all $\beta$ and verify its dimension when $\beta>0$.

Write $F(t)=c_0t^\delta(\log t)^\beta$ for large $t$.  Choose
$0<\rho<1$; for $\beta>0$ require
\begin{equation}\label{eq:logrho}
 \rho>\frac{1-\delta}{\beta+1-\delta}.
\end{equation}
Set $T_j=\lfloor e^{j^\rho}\rfloor$ and $\Delta_j=T_{j+1}-T_j$.
Fix a finite guarded block family of common length $L$, with roofs
bounded by $M$.  At physical time $T_j$ concatenate
$
 \ell_j=L\left\lfloor\frac{F(T_{j+1})-F(T_j)}L\right\rfloor
$
free letters in complete guarded blocks.  Then insert the return
\[
 B_j=\Delta_j-\sum_{\text{free letters in stage }j}\tau.
\]
The elementary difference asymptotics are
\begin{equation}\label{eq:logtimes}
\begin{aligned}
 \Delta_j\sim\rho e^{j^\rho}j^{\rho-1},\quad 
 \ell_j\sim c_0\delta\rho e^{\delta j^\rho}
                       j^{\rho\beta+\rho-1}.
\end{aligned}
\end{equation}
In particular $M\ell_j=o(\Delta_j)$, so $B_j>0$ and
$B_j\sim\Delta_j$ uniformly over free choices.  Every insertion
is causal and lies between complete blocks.  An initial prescribed
return reaches the first $T_j$.

Telescoping and the $O(J)$ rounding and bridge counts show that
\begin{equation}\label{eq:logcount}
 N_y(T_{J+1})\sim K_J:=\sum_{j\le J}\ell_j
 \sim F(T_{J+1})\asymp e^{\delta J^\rho}J^{\rho\beta}.
\end{equation}
Since $T_{j+1}/T_j\to1$ and $F(T_{j+1})/F(T_j)\to1$, monotonicity
squeezes the counting function between consecutive stage endpoints.
This proves the exact asymptotic at every physical time, uniformly in
the freely chosen words.

We next verify the full macroscopic spectrum.  A dyadic shell of size
$R=2^n$ meets only $(\log R)^{O(1)}$ stages, by
\eqref{eq:logtimes}.  Each free cluster it meets has
$R^\delta(\log R)^{O(1)}$ points and span of this same order, with
successive gaps between $1$ and $M$.  Here and below this notation
means upper and lower bounds by possibly different fixed powers of
$\log R$.  To see that every sufficiently large shell contains a complete
cluster, recall that $S_n=[R/2,R)\cap\N$. Since
$T_{j+1}/T_j\to1$, there is a stage start
$T_j\in[5R/8,3R/4]$ for every sufficiently large $R=2^n$.
By \eqref{eq:logtimes}, $\Delta_j=o(R)$ and
$M\ell_j=o(\Delta_j)$ for these stages, uniformly in the
free choices. Thus the entire free cluster beginning at $T_j$
lies in $S_n$. A shell can intersect partial clusters only at its
two endpoints, so boundary fragments affect the upper bound by at
most two clusters.
For $0<s<1$, cover each cluster by maximal permitted balls if
$U=2^{n-N\theta}$ is smaller than its span, and otherwise by one ball.
Together with unit-ball covers this gives
\[
 \mathcal H_{n,N}^{s,\theta}(R(y,Y))
 \le C n^C\left(2^{(\delta-1)n+N\theta(1-s)}
                     +2^{-(1-\delta)sn}\right).
\]
The constant $C$ may depend on $s,\beta,\rho,M$ but not on $n,N$
or the free choices.  The isolated bridge endpoints obey the same
upper estimate.  Summing over all $n\ge N$ gives exponential decay
if $s>\varphi_\delta(\theta)$; at $\theta=0$ this holds for every $s>0$.
For the converse, a complete cluster in shell $n=N$ contains at least
$2^{\delta N}N^{-C}$ points.  Any admissible cover has
$\sum r_i\ge c2^{\delta N}N^{-C}$, since a ball of radius $r_i\ge1$
covers at most $3r_i$ integers.  Hence
\[
 \mathcal H_N^{s,\theta}(R(y,Y))
 \ge cN^{-C}2^{N[\delta-1+\theta(1-s)]}.
\]
This grows exponentially for $s<\varphi_\delta(\theta)$ whenever that
value is positive.  This proves the asserted spectrum, including
$\theta=0$.

Finally, by the upper geometric bound and \eqref{eq:logtimes}, the
prescribed contraction through stage $J$ is at most
\begin{equation}\label{eq:logenergy}
 A_J\le C\sum_{j\le J}(1+B_j^\delta)
 \le C e^{\delta J^\rho}J^{(1-\delta)(1-\rho)}.
\end{equation}
To justify the last estimate, compare the sum of
$e^{\delta j^\rho}j^{-\delta(1-\rho)}$ with its integral, substitute
$u=x^\rho$, and integrate the exponentially increasing integrand;
the resulting power is
$1-\rho-\delta(1-\rho)=(1-\delta)(1-\rho)$.
Equations \eqref{eq:logcount}--\eqref{eq:logenergy} yield
$
 \frac{A_J}{K_J}\le
 C J^{(1-\delta)(1-\rho)-\rho\beta}\longrightarrow0
$
by \eqref{eq:logrho}.  Since $K_{J-1}/K_J\to1$, the same bound
controls every partial stage and the position immediately after its
bridge.  For $\beta>0$ use Lemma~\ref{lem:rigorous-forcing} with
arbitrary $M\ge c$ and $\epsilon>0$, then let
$\epsilon\downarrow0$ and $M\to\infty$.  This proves the full-dimensional case of
\eqref{eq:logcriticalphase}.  For $\beta\le0$, repeating a single interior
guarded block in the same construction proves nonemptiness and the
claimed spectrum.  Dropping the spectrum constraint cannot decrease
dimension, and the zero-dimensional upper bounds already apply to the
unconstrained levels.  The final assertion follows.
\end{proof}

\begin{remark}\label{rem:log-endpoint-universal}
The distinction at $\beta=0$ is not an artefact of the chosen
return-time construction.  The zero-dimensional conclusion at this
endpoint follows for every itinerary with the specified exact
counting asymptotic, since its mean $\delta$-energy diverges by
Theorem~\ref{thm:exactcritical}.  By contrast, any fixed positive
logarithmic correction permits an insertion scheme whose prescribed
contraction is negligible relative to the free entropy.  The
critical assertion is thus a consequence of the counting and
contraction estimates, not of the particular sequence of stages.
\end{remark}

\subsection{A two-branch interval map with stretched-exponential returns}
\begin{figure}[tbp]
\centering
\begin{minipage}[t]{0.48\textwidth}
\centering
\begin{tikzpicture}[x=5.55cm,y=4.55cm,>=Latex]
  \fill[black!5] (0.56600,0) rectangle (1,1);
  \draw[->] (-0.025,0) -- (1.06,0) node[right] {$x$};
  \draw[->] (0,-0.025) -- (0,1.06) node[left] {$T_\delta(x)$};
  \draw[densely dashed,thin] (0,0) -- (1,1);
  \draw[thick] plot[smooth] coordinates {
  (0.00000,0.00000) (0.00674,0.00866) (0.01348,0.01762) (0.02021,0.02676) (0.02695,0.03603) (0.03369,0.04541) (0.04043,0.05489) (0.04717,0.06447)
  (0.05390,0.07412) (0.06064,0.08386) (0.06738,0.09367) (0.07412,0.10356) (0.08086,0.11352) (0.08759,0.12354) (0.09433,0.13363) (0.10107,0.14379)
  (0.10781,0.15401) (0.11455,0.16430) (0.12128,0.17464) (0.12802,0.18505) (0.13476,0.19552) (0.14150,0.20604) (0.14824,0.21663) (0.15498,0.22727)
  (0.16171,0.23798) (0.16845,0.24874) (0.17519,0.25956) (0.18193,0.27044) (0.18867,0.28138) (0.19540,0.29238) (0.20214,0.30343) (0.20888,0.31455)
  (0.21562,0.32572) (0.22236,0.33695) (0.22909,0.34824) (0.23583,0.35958) (0.24257,0.37099) (0.24931,0.38246) (0.25605,0.39398) (0.26278,0.40557)
  (0.26952,0.41721) (0.27626,0.42892) (0.28300,0.44069) (0.28974,0.45252) (0.29647,0.46441) (0.30321,0.47636) (0.30995,0.48838) (0.31669,0.50046)
  (0.32343,0.51261) (0.33016,0.52482) (0.33690,0.53709) (0.34364,0.54943) (0.35038,0.56184) (0.35712,0.57432) (0.36385,0.58686) (0.37059,0.59947)
  (0.37733,0.61215) (0.38407,0.62491) (0.39081,0.63773) (0.39754,0.65063) (0.40428,0.66360) (0.41102,0.67664) (0.41776,0.68976) (0.42450,0.70296)
  (0.43124,0.71623) (0.43797,0.72958) (0.44471,0.74302) (0.45145,0.75653) (0.45819,0.77012) (0.46493,0.78380) (0.47166,0.79756) (0.47840,0.81141)
  (0.48514,0.82535) (0.49188,0.83937) (0.49862,0.85349) (0.50535,0.86770) (0.51209,0.88200) (0.51883,0.89639) (0.52557,0.91089) (0.53231,0.92548)
  (0.53904,0.94018) (0.54578,0.95497) (0.55252,0.96987) (0.55926,0.98488) (0.56600,1.00000)
  };
  \draw[thick] (0.56600,0) -- (1,1);
  \draw[dotted] (0.56600,0) -- (0.56600,1);
  \fill (0,0) circle (1pt);
  \draw[fill=white] (0.56600,1) circle (1.2pt);
  \fill (0.56600,0) circle (1.2pt);
  \node[below] at (0.56600,0) {$a$};
  \node[below] at (0.985,-0.015) {$1$};
  \node[above left] at (0.26,0.40) {$f_\delta$};
  \node[right] at (0.72,0.30) {$g$};
  \node[align=left,anchor=west] at (0.04,0.83) {neutral fixed\\point at $0$};
  \node at (0.80,0.90) {$Y=[a,1]$};
\end{tikzpicture}

\smallskip
\small (a) Two branches of $T_\delta$: the neutral left branch $f_\delta$ and the linear right branch $g$.
\end{minipage}
\hfill
\begin{minipage}[t]{0.48\textwidth}
\centering
\begin{tikzpicture}[x=5.55cm,y=4.55cm,>=Latex]
  \draw[->] (-0.025,0) -- (1.04,0);
  \node[below] at (0.54,-0.11) {$\frac{x-a}{1-a}$};
  \draw[->] (0,-0.025) -- (0,1.06) node[left] {$\frac{F(x)-a}{1-a}$};
  \draw[thick] plot[smooth] coordinates {
    (0.07867,0.00000) (0.07991,0.03226) (0.08114,0.06452) (0.08235,0.09677) (0.08354,0.12903) (0.08472,0.16129) (0.08588,0.19355) (0.08702,0.22581)
    (0.08815,0.25806) (0.08926,0.29032) (0.09035,0.32258) (0.09143,0.35484) (0.09250,0.38710) (0.09355,0.41935) (0.09459,0.45161) (0.09561,0.48387)
    (0.09662,0.51613) (0.09762,0.54839) (0.09860,0.58065) (0.09957,0.61290) (0.10053,0.64516) (0.10147,0.67742) (0.10241,0.70968) (0.10333,0.74194)
    (0.10423,0.77419) (0.10513,0.80645) (0.10602,0.83871) (0.10689,0.87097) (0.10775,0.90323) (0.10860,0.93548) (0.10944,0.96774) (0.11027,1.00000)
  };
  \draw[thick] plot[smooth] coordinates {
    (0.11027,0.00000) (0.11212,0.03226) (0.11394,0.06452) (0.11573,0.09677) (0.11751,0.12903) (0.11926,0.16129) (0.12098,0.19355) (0.12269,0.22581)
    (0.12437,0.25806) (0.12603,0.29032) (0.12767,0.32258) (0.12928,0.35484) (0.13088,0.38710) (0.13246,0.41935) (0.13402,0.45161) (0.13556,0.48387)
    (0.13708,0.51613) (0.13858,0.54839) (0.14006,0.58065) (0.14153,0.61290) (0.14297,0.64516) (0.14440,0.67742) (0.14581,0.70968) (0.14721,0.74194)
    (0.14858,0.77419) (0.14994,0.80645) (0.15129,0.83871) (0.15261,0.87097) (0.15393,0.90323) (0.15522,0.93548) (0.15650,0.96774) (0.15777,1.00000)
  };
  \draw[thick] plot[smooth] coordinates {
    (0.15777,0.00000) (0.16058,0.03226) (0.16337,0.06452) (0.16612,0.09677) (0.16883,0.12903) (0.17152,0.16129) (0.17417,0.19355) (0.17680,0.22581)
    (0.17939,0.25806) (0.18196,0.29032) (0.18450,0.32258) (0.18700,0.35484) (0.18948,0.38710) (0.19194,0.41935) (0.19436,0.45161) (0.19676,0.48387)
    (0.19913,0.51613) (0.20147,0.54839) (0.20379,0.58065) (0.20608,0.61290) (0.20835,0.64516) (0.21059,0.67742) (0.21281,0.70968) (0.21500,0.74194)
    (0.21717,0.77419) (0.21932,0.80645) (0.22144,0.83871) (0.22354,0.87097) (0.22561,0.90323) (0.22767,0.93548) (0.22970,0.96774) (0.23171,1.00000)
  };
  \draw[thick] plot[smooth] coordinates {
    (0.23171,0.00000) (0.23618,0.03226) (0.24061,0.06452) (0.24500,0.09677) (0.24935,0.12903) (0.25366,0.16129) (0.25793,0.19355) (0.26216,0.22581)
    (0.26635,0.25806) (0.27050,0.29032) (0.27461,0.32258) (0.27868,0.35484) (0.28271,0.38710) (0.28671,0.41935) (0.29067,0.45161) (0.29460,0.48387)
    (0.29848,0.51613) (0.30233,0.54839) (0.30615,0.58065) (0.30993,0.61290) (0.31367,0.64516) (0.31738,0.67742) (0.32106,0.70968) (0.32470,0.74194)
    (0.32831,0.77419) (0.33188,0.80645) (0.33542,0.83871) (0.33893,0.87097) (0.34240,0.90323) (0.34584,0.93548) (0.34925,0.96774) (0.35263,1.00000)
  };
  \draw[thick] plot[smooth] coordinates {
    (0.35263,0.00000) (0.36017,0.03226) (0.36767,0.06452) (0.37513,0.09677) (0.38253,0.12903) (0.38990,0.16129) (0.39722,0.19355) (0.40449,0.22581)
    (0.41172,0.25806) (0.41890,0.29032) (0.42604,0.32258) (0.43314,0.35484) (0.44019,0.38710) (0.44720,0.41935) (0.45416,0.45161) (0.46109,0.48387)
    (0.46796,0.51613) (0.47480,0.54839) (0.48159,0.58065) (0.48834,0.61290) (0.49505,0.64516) (0.50171,0.67742) (0.50833,0.70968) (0.51491,0.74194)
    (0.52144,0.77419) (0.52793,0.80645) (0.53438,0.83871) (0.54079,0.87097) (0.54715,0.90323) (0.55348,0.93548) (0.55976,0.96774) (0.56600,1.00000)
  };
  \draw[thick] plot[smooth] coordinates {
    (0.56600,0.00000) (0.58000,0.03226) (0.59400,0.06452) (0.60800,0.09677) (0.62200,0.12903) (0.63600,0.16129) (0.65000,0.19355) (0.66400,0.22581)
    (0.67800,0.25806) (0.69200,0.29032) (0.70600,0.32258) (0.72000,0.35484) (0.73400,0.38710) (0.74800,0.41935) (0.76200,0.45161) (0.77600,0.48387)
    (0.79000,0.51613) (0.80400,0.54839) (0.81800,0.58065) (0.83200,0.61290) (0.84600,0.64516) (0.86000,0.67742) (0.87400,0.70968) (0.88800,0.74194)
    (0.90200,0.77419) (0.91600,0.80645) (0.93000,0.83871) (0.94400,0.87097) (0.95800,0.90323) (0.97200,0.93548) (0.98600,0.96774) (1.00000,1.00000)
  };
  \node[below] at (0,0) {$0$};
  \node[below] at (1,0) {$1$};
  \node[fill=white,inner sep=1pt] at (0.88,0.44) {$n=1$};
  \node[fill=white,inner sep=1pt] at (0.53,0.38) {$2$};
  \node[fill=white,inner sep=1pt] at (0.38,0.42) {$3$};
  \node[anchor=west] at (0.04,0.89) {$n\to\infty$};
  \draw[->,thin] (0.145,0.86) -- (0.025,0.42);
\end{tikzpicture}

\smallskip
\small (b) Induced first-return map $F=T_\delta^\tau|_Y$. The full branches $Y_n$ accumulate at $a$.
\end{minipage}
\caption{The two-branch intermittent map of Proposition~\ref{prop:intervalmodel} and its genuine first-return map. The curves are drawn for $\delta=0.65$; the branch labels $n$ in (b) are the actual return times. The right-hand axes rescale the target interval $Y=[a,1]$ to $[0,1]$.}
\label{fig:stretched-interval-map}
\end{figure}
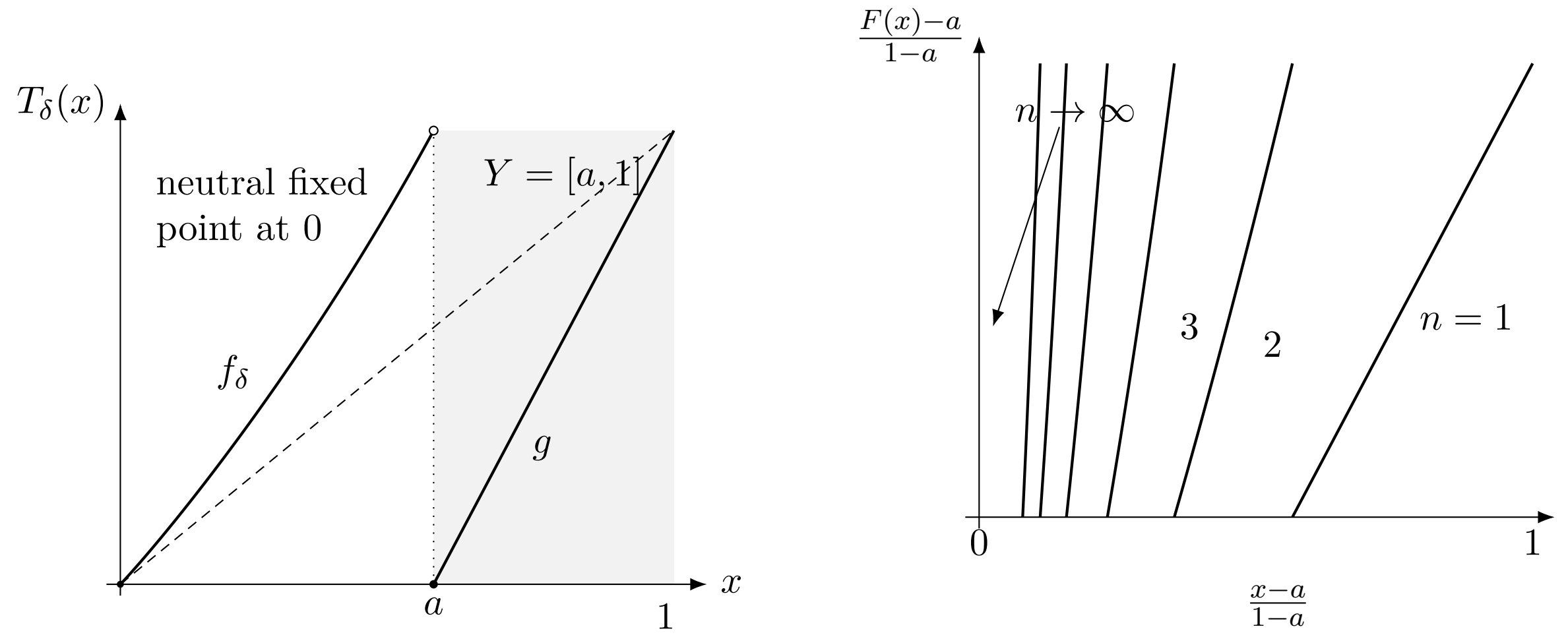

We now realize the required roof and branch geometry as actual first
returns of a map of one interval.  In particular the preceding phase
transitions are not confined to a tower endowed with a chosen roof.
Figure~\ref{fig:stretched-interval-map} illustrates the two branches
of the interval map and its induced first-return branches.

\begin{proposition}\label{prop:intervalmodel}
Fix $0<\delta<1$, put $a=\exp(1-2^\delta)$, and define
\begin{equation}\label{eq:intervalmap}
 T_\delta(x)=
 \begin{cases}
 0,&x=0,\\
 \displaystyle\exp\!\left(1-
       \left[(1-\log x)^{1/\delta}-1\right]^\delta\right),
       &0<x<a,\\
 (x-a)/(1-a),&a\le x\le1.
 \end{cases}
\end{equation}
This map has two increasing branches, a neutral fixed point at zero,
and a genuine first-return map to $Y=[a,1]$ with one branch of every
return length $n\ge1$.  Its inverse first-return branches are
\begin{equation}\label{eq:intervalinverse}
 \psi_n(y)=a+(1-a)\exp\left(1-[u(y)+n-1]^\delta\right),
 \qquad u(y)=(1-\log y)^{1/\delta}.
\end{equation}
They satisfy \textup{(G1), (G2)} and, uniformly for $y\in Y$,
\begin{equation}\label{eq:intervalderivative}
 |\psi_n'(y)|\asymp n^{\delta-1}e^{-n^\delta},\qquad
 -\log|\psi_n'(y)|=n^\delta+(1-\delta)\log n+O(1).
\end{equation}
The induced limit set has Hausdorff dimension one.
\end{proposition}
\begin{proof}
Let $f$ denote the left branch and set
$x_m=\exp(1-(m+1)^\delta)$ for $m\ge0$.  Then $x_0=1$,
$x_1=a$, $x_m\downarrow0$, and $f(x_m)=x_{m-1}$ for $m\ge2$,
with the same identity at $m=1$ using the left limit.  The inverse
iterates of the left branch satisfy
$
 f^{-m}(y)=\exp(1-[u(y)+m]^\delta), ~~ y\in[a,1].
$
For a point of $Y$, the first application of $T_\delta$ is the right
branch.  If its image belongs to $(x_n,x_{n-1})$, it then makes
exactly $n-1$ left-branch iterates before returning to $Y$.
Inverting this itinerary gives \eqref{eq:intervalinverse}, including
$n=1$.  The intervals $\psi_n(Y)$ have disjoint interiors and cover
$Y$ up to endpoints and the accumulation point $a$.

Since $1\le u(y)\le2$, differentiation gives
\[
 \psi_n'(y)=(1-a)\frac{e^{1-[u(y)+n-1]^\delta}}{y}
                  \left(\frac{u(y)}{u(y)+n-1}\right)^{1-\delta}.
\]
This is comparable to $n^{\delta-1}e^{-n^\delta}$ uniformly on $Y$,
because $(n+u-1)^\delta-n^\delta=O(1)$; taking logarithms
also gives the second estimate in \eqref{eq:intervalderivative}.
Both factors after $1-a$
are at most one, so $\sup_n\|\psi_n'\|_\infty\le1-a<1$.
The derivatives of $\log\psi_n'$ on $Y$ are uniformly bounded:
$u,u'$ are bounded there, $u\ge1$, and the powers of $u+n-1$
that arise have nonpositive exponents.  Summing this Lipschitz bound
along contracting compositions gives uniform bounded distortion.
Moreover $\psi_2(Y)$ is compactly contained in the interior of $Y$.

The depth-$k$ branch intervals partition $Y$ apart from countably many
endpoints and have diameters at most $|Y|(1-a)^k$.  If $s<1$, then
\[
 \sum_{|w|=k}|\psi_w(Y)|^s
 \ge [|Y|(1-a)^k]^{s-1}|Y|\longrightarrow\infty.
\]
Bounded distortion therefore permits a finite collection of depth-$k$
words with $\sum_w(\inf_Y|\psi_w'|)^s>1$.  The finite conformal IFS
formed by these words has dimension greater than $s$, by the
finite-system pressure criterion, and is contained in $K_M$ for some
$M$.  Thus $\sup_M\dim_HK_M=1$.  Equivalently, every point except
the countable backward orbit of partition endpoints has an infinite
first-return itinerary, so $\dim_HK_\infty=1$.

Finally, direct differentiation on the left branch shows
\[
 f'(x)=\frac{f(x)}x
 \left(\frac{(1-\log x)^{1/\delta}}
 {(1-\log x)^{1/\delta}-1}\right)^{1-\delta}>1,
 \qquad f'(x)\longrightarrow1\quad(x\downarrow0).
\]
Also $f(x)/x\to1$, proving differentiability at zero and $f'(0)=1$.
The map is smooth in the interior of each branch and $C^1$ at zero
from the right.  Its precise regularity limitation is
$f'(x)-1\asymp(\log(1/x))^{-(1-\delta)/\delta}$; hence the left
branch is not $C^{1+\eta}$ at zero for any $\eta>0$.
The induced branches, on which all our hypotheses are imposed, have
the uniform regularity just proved.
\end{proof}

\begin{corollary}\label{cor:intervalphases}
For $T_\delta$ and $Y=[a,1]$, every scale-admissible return-time spectrum
$h$ has a zero-frequency level of dimension zero if $h(1)<\delta$ and
one if $h(1)\ge\delta$.  Exact counting levels $N_y(t)\sim c_0t^b$
have dimension zero for $b\le\delta$ and one for $b>\delta$.
At the counting threshold, the levels
$N_y(t)\sim c_0t^\delta(\log t)^\beta$ with complete spectrum
$\varphi_\delta$ have dimension zero for $\beta\le0$ and one for $\beta>0$.
\end{corollary}
\begin{proof}
Apply Theorems~\ref{thm:generalphase}, \ref{thm:exactpower} and
\ref{thm:logcritical} using Proposition~\ref{prop:intervalmodel}.
\end{proof}

\section{Exact counting and intrinsic energy levels}\label{sec:secondorder}
Above the counting threshold, prescribing an intrinsic mean return
energy distinguishes complete refined levels with different Hausdorff
dimensions.  We first compute the spectra of the polynomial schedules
used in their construction, then prove the dimension formula for the
entire joint level.  A final construction illustrates the different
role of a prescribed insertion schedule.

We now work with the affine full-branch return tower from
Corollary~\ref{cor:lineartower}, with
\begin{equation}\label{eq:affineratio2}
 r_n=Z^{-1}e^{-\kappa n^\delta},\qquad
 Z=\sum_{m\ge1}e^{-\kappa m^\delta},\qquad
 0<\delta<1,\quad\kappa>0.
\end{equation}

\subsection{Polynomial schedules and macroscopic spectra}
The following estimates include the outer covering scale and are
uniform in the bounded free symbols.  They will impose the same
complete spectrum on every point of the energy-level construction.

\begin{lemma}\label{lem:powerprofile}
Fix $0<p<q$ and $A,B>0$.  Suppose an integer return-time
sequence is obtained by concatenating, at stage $j$,
$\ell_j=\lfloor Bj^p\rfloor+O(1)$ free return symbols from a
fixed finite alphabet of positive integers, followed by one
prescribed return symbol $n_j=\lfloor Aj^q\rfloor\vee n_0$,
where $n_0\ge1$ is fixed.  Let $E$ be the resulting set of
physical return times.  Put
\[
 b_{p,q}=\frac{p+1}{q+1},\qquad
 d_{p,q}=\frac1{q+1-p}.
\]
Then the lower and upper macroscopic intermediate dimensions agree
for every $\theta\in[0,1]$, and their common value is
\begin{equation}\label{eq:general-powerprofile}
 \mathfrak p_{p,q}(0)=d_{p,q},\qquad
 \mathfrak p_{p,q}(\theta)=
 \max\left\{d_{p,q},
       \left(1-\frac{1-b_{p,q}}{\theta}\right)_+\right\}
 \quad(0<\theta\le1).
\end{equation}
The transition between the constant and
increasing parts occurs at $\theta_*=1-p/(q+1)$.
All estimates are uniform over the free symbol choices; the spectrum
is independent of $A,B$ and the particular fixed finite alphabet.
\end{lemma}
\begin{proof}
We give the covering estimates, including the outer scale $N$.
Put $u=1/(q+1)$, $v=p/(q+1)$ and $w=q/(q+1)$, so that
$u+w=1$, $b_{p,q}=u+v$ and $d_{p,q}=u/(1-v)$.
Write $R=2^n$.  In a fixed interior portion of every sufficiently
large dyadic shell, the prescribed construction consists of
$\asymp R^u$ occupied blocks, each containing $\asymp R^v$
integer return times with gaps bounded above and below by
positive constants determined by the fixed free alphabet.  Consecutive
blocks are separated by $\asymp R^w$, where $u+w=1$ and
$0<v<w$.  The same estimates, with uniform constants, hold over
each entire dyadic shell after dropping at most two boundary blocks.

Put equal mass on the return times in the chosen interior portion of
one shell and normalize it to a probability $\nu_n$.  An interval
of radius $r$, $1\le r\le R$, satisfies
\begin{equation}\label{eq:mass3regimes}
\nu_n(B(x,r))\lesssim
\begin{cases}
 r/R^{b_{p,q}},&1\le r\le R^v,\\
 R^{-u},&R^v\le r\le R^w,\\
 r/R,&R^w\le r\le R.
\end{cases}
\end{equation}
In the first case the ball intersects at most one occupied block,
up to harmless uniform constants; in the last case it intersects
$O(r/R^w)$ blocks.  The constants are independent of the free
symbol sequence.  The mass-distribution inequality
$\sum_i(r_i/R)^s\gtrsim
 \inf_{1\le r\le R^{1-\theta}}(r/R)^s/
 \sup_x\nu_n(B(x,r))$
then gives, for $0<s<1$ and an outer scale $N=n$,
\[ \mathcal H_{n,n}^{s,\theta}(E)\gtrsim
\begin{cases}
 \min\{R^{u-s(1-v)},R^{\theta(1-s)}\},
    &0<\theta\le1-v,\\
 R^{b_{p,q}-1+\theta(1-s)},
    &1-v<\theta\le1.
\end{cases} \]
For the second line, the maximal allowed radius is already
smaller than $R^v$, so only the first regime of
\eqref{eq:mass3regimes} is relevant.  Whenever
$s<\mathfrak p_{p,q}(\theta)$ and $\theta>0$, the appropriate lower
bound grows exponentially in $n$ (for $1-v<\theta$ use
$s<1-(1-b_{p,q})/\theta$).  This proves the lower bound for both
macroscopic intermediate dimensions at every positive $\theta$.

At $\theta=0$, the estimates \eqref{eq:mass3regimes} give a
positive shell content, uniformly in $n$, whenever $s<d_{p,q}$:
indeed the infimum of $(r/R)^s/\sup_x\nu_n(B(x,r))$ over
$1\le r\le R$ is bounded below by a constant.  Since the
covering sum includes \emph{all} shells $n\ge N$, it diverges.
Thus $\mathsf D_0^\pm(E)\ge d_{p,q}$.

For the upper bounds, cover every occupied block of shell $n$
by one ball of radius $\asymp2^{vn}$.  This produces a shell
contribution $\lesssim2^{[u-s(1-v)]n}$ whenever
$2^{vn}\lesssim2^{n-N\theta}$; the sum over admissible shells
decays when $s>d_{p,q}$.  If the block radius exceeds the
window, cover the return times by balls of maximal permitted
radius $2^{n-N\theta}$, resulting in
$
 \mathcal H_{n,N}^{s,\theta}(E)
 \lesssim 2^{(b_{p,q}-1)n+N\theta(1-s)}.
$
The small-radius regime is confined to
$N\le n\lesssim N\theta/(1-v)$, and its sum decays when
$s>1-(1-b_{p,q})/\theta$.  Hence $\mathcal H_N^{s,\theta}(E)$
decays exponentially for every
$s>\mathfrak p_{p,q}(\theta)$, including $\theta=0$.  The two bounds
prove \eqref{eq:general-powerprofile}.
\end{proof}

\subsection{Intrinsic energy and complete level sets}
We prescribe the mean of the \emph{total} $\delta$-energy of the
return symbols, an observable intrinsic to the orbit.  The resulting
dimension formula holds for the complete joint level.  Its main
content is that the counting and spectrum constraints do not lower
the dimension of the underlying energy level.

Write
\[
 Z(\lambda)=\sum_{n\ge1}e^{-\lambda n^\delta},\qquad
 \mathcal E(\lambda)=\frac{\sum_{n\ge1}n^\delta e^{-\lambda n^\delta}}
                          {Z(\lambda)},\qquad \lambda>0.
\]
The function $\mathcal E$ is strictly decreasing from infinity to $1$,
so each $E>1$ corresponds to a unique $\lambda(E)>0$ with
$\mathcal E(\lambda(E))=E$.  Define
\begin{equation}\label{eq:intrinsic-dim}
 \mathfrak d_\kappa(E)=
 \frac{\lambda(E)E+\log Z(\lambda(E))}
      {\kappa E+\log Z(\kappa)}.
\end{equation}
The numerator is the largest Shannon entropy of any probability law
$p=(p_n)_{n\ge1}$ satisfying $\sum_n p_n n^\delta=E$;
indeed
\begin{equation}\label{eq:entropy-variational}
 H(p)\le\lambda E+\log Z(\lambda),
 \quad\text{with equality at }p_n=Z(\lambda)^{-1}e^{-\lambda n^\delta}.
\end{equation}

Fix a counting exponent $\delta<b<1$.  Choose
\begin{equation}\label{eq:thinbridgeq}
 q>\frac{1-b}{b-\delta},\qquad
 p=b(q+1)-1.
\end{equation}
Then $0<q\delta<p<q$.  Put
\[
 u=\frac1{q+1},\quad v=\frac p{q+1}=b-u,\quad
 d_{q,b}=\frac{u}{1-v}=\frac1{1+(1-b)(q+1)},
\]
and let
\begin{equation}\label{eq:intrinsic-profile}
 h_{q,b}(0)=d_{q,b},\qquad
 h_{q,b}(\theta)=\max\left\{d_{q,b},
 \left(1-\frac{1-b}{\theta}\right)_+\right\}
 \quad(0<\theta\le1).
\end{equation}
The lower bound will use the following exact-energy
approximation, with no limiting error in the prescribed mean.
\begin{lemma}\label{lem:finite-exact-energy}
Let $0<\delta<1$ and $E>1$.  With
$H(E)=\lambda(E)E+\log Z(\lambda(E))$, there are
finitely supported probability vectors $P^{(M)}$ on $\N$ such that
$
 \sum_nP_n^{(M)}n^\delta=E$ and $ H(P^{(M)})\rightarrow H(E).$
Fix any guard symbol $g\in\N$.  For all sufficiently large
block lengths $L$, there are probability laws on finitely many
words $B=(n_1,\ldots,n_{L-1},g)$ such that
\[
 \frac1L\sum_Bp_B\sum_{i=1}^L n_i^\delta=E,
 \qquad \lim_{L\to\infty}\frac{H(p)}L=H(E).
\]
For the affine ratios $r_n=Z(\kappa)^{-1}e^{-\kappa n^\delta}$,
the expected geometric contraction per letter for these words is
exactly $\kappa E+\log Z(\kappa)$.
\end{lemma}
\begin{proof}
Let $p_n=Z(\lambda(E))^{-1}e^{-\lambda(E)n^\delta}$.
Condition this law on $n\le M$ to obtain $p^{[M]}$ with
mean $E_M<E$.  Then $E_M\uparrow E$ and
$H(p^{[M]})\to H(p)=H(E)$ by convergence of the weighted
sums defining the entropy.  Choose a fixed integer $m_+>E^{1/\delta}$
and put
\[
 \varepsilon_M=\frac{E-E_M}{m_+^\delta-E_M},\qquad
 P^{(M)}=(1-\varepsilon_M)p^{[M]}+\varepsilon_M\delta_{m_+}.
\]
The new mean equals $E$ exactly, $\varepsilon_M\to0$, and
$H(P^{(M)})\to H(E)$ by the entropy bound for a binary mixture.
For the guarded assertion put
$E_L=(LE-g^\delta)/(L-1)>1$ for large $L$.  Apply the same
finite-support construction at energy $E_L$, selecting the
finite law $P_L$ so that its entropy differs from $H(E_L)$ by
at most $L^{-1}$.  Choose the first $L-1$ symbols independently
with law $P_L$ and append $g$.  Their joint entropy is
$(L-1)H(P_L)$ and their expected total energy is $LE$.
Continuity of $H(E)$, a consequence of the strictly decreasing
Gibbs energy parameter, gives the required entropy limit.
The contraction identity follows by summing
$-\log r_n=\kappa n^\delta+\log Z(\kappa)$ over the word.
\end{proof}

The next theorem calculates the Hausdorff dimension of the
\emph{entire} level set with all three quantities specified, rather
than merely constructing a subset of a given dimension.
\begin{theorem}\label{thm:intrinsicpressure}
For the affine full-branch tower \eqref{eq:affineratio2}, let
$0<\delta<1$, $\delta<b<1$, $c_0>0$, and choose $q$ as in
\eqref{eq:thinbridgeq}.  For each $E>1$ define the joint level set
\begin{equation}\label{eq:secondorder-levelset}
\mathscr L(c_0,b,E;q)=\left\{y\in K_\infty:\;
\begin{aligned}
&N_y(t)\sim c_0t^b,\qquad
\lim_{k\to\infty}\frac1k\sum_{i=1}^k\tau_i(y)^\delta=E,\\
&\mathsf D^-_\theta(R(y,Y))=\mathsf D^+_\theta(R(y,Y))
=h_{q,b}(\theta)\quad(\forall\,\theta\in[0,1])
\end{aligned}\right\}.
\end{equation}
Then
\begin{equation}\label{eq:intrinsicpressure-theorem}
 \dim_H\mathscr L(c_0,b,E;q)=\mathfrak d_\kappa(E).
\end{equation}
Moreover $\mathfrak d_\kappa(E)$ is continuous,
$0<\mathfrak d_\kappa(E)\le1$,
and it equals $1$ precisely at
$E=\mathcal E(\kappa)$.  It tends to zero as
$E\downarrow1$ or $E\to\infty$.  Thus the exact second-order
energy levels, even with the \emph{same} counting coefficient and
complete macroscopic spectrum, realize every Hausdorff dimension
in $(0,1]$.
\end{theorem}
\begin{proof}
\emph{Upper bound.}
Fix $E>1$ and $\lambda=\lambda(E)$.  At symbolic depth $k$, every
word $w=(n_1,\ldots,n_k)$ in the energy window
$|k^{-1}\sum n_i^\delta-E|\le\varepsilon$ has
\[
 |I_w|=Z(\kappa)^{-k}
      \exp\left(-\kappa\sum_{i=1}^k n_i^\delta\right)
 \le\exp[-k(\kappa(E-\varepsilon)+\log Z(\kappa))].
\]
The number of such words is at most
$
 \exp\bigl(k[\lambda(E+\varepsilon)+\log Z(\lambda)]\bigr),
$
because $\sum_{|w|=k}e^{-\lambda\sum_i n_i^\delta}
=Z(\lambda)^k$ and each word in the window carries weight
at least $e^{-\lambda(E+\varepsilon)k}$.
For every $s>\mathfrak d_\kappa(E)$ choose $\varepsilon>0$
small enough that the sum of $s$-powers of these cylinder
diameters decays exponentially in $k$.  All points having average
energy $E$ belong to these window families for all sufficiently
large $k$.  Hence their Hausdorff dimension is at most
$\mathfrak d_\kappa(E)$, without requiring the counting or macroscopic dimension-spectrum
constraints.

\emph{Lower bound.}
We construct a family of points obeying all three constraints and
having dimension arbitrarily close to the right side.
For positive $A,B$, at stage $j$ arrange $\ell_j=\lfloor B j^p\rfloor$
free inducing digits, followed by one prescribed digit
$n_j=\lfloor Aj^q\rfloor\vee5$.
Choose $B$ from
\[ \frac{B}{p+1}\left(\frac{q+1}{A}\right)^b=c_0. \]
The number of free symbols through stage $J$ is
$K_J\sim BJ^{p+1}/(p+1)$, and physical time is
$T_J\sim AJ^{q+1}/(q+1)$.  Bounded free roof values give
$N_y(t)\sim c_0t^b$ uniformly over all free choices, since
$n_J=o(T_J)$ and $\ell_J=o(K_J)$.
The prescribed long digits contribute
\begin{equation}\label{eq:vanishingenergy}
 \sum_{j\le J}n_j^\delta=O(J^{q\delta+1})
   =o(J^{p+1})=o(K_J).
\end{equation}
In fact their accumulated negative logarithmic contraction is
$O(J^{q\delta+1}+J)$, and a single prescribed symbol has
contraction and energy $O(1+J^{q\delta})$.  Both estimates are
uniform over free choices.  In particular, once a finite alphabet
and a guarded block length $L$ have been fixed, and $k$ counts
completed free blocks, the insertion penalty in
Lemma~\ref{lem:rigorous-forcing} can be taken to be
\begin{equation}\label{eq:power-prefix-penalty}
 A(k)=C_{M,L}(1+k)^\rho,\qquad
 \rho=\frac{q\delta+1}{p+1}<1.
\end{equation}
Indeed, the stage containing the $(k+1)$st free block satisfies
$J\lesssim_{M,L}(1+k)^{1/(p+1)}$; the bound includes a prescribed
symbol immediately before that block.  Thus the explicit estimate
\eqref{eq:power-prefix-penalty} controls \emph{every} free-block
prefix, not only the stage endpoints.

Apply Lemma~\ref{lem:finite-exact-energy} with an interior
guard $g$.  It supplies a finitely supported law on complete
$L$-letter guarded blocks whose mean energy is \emph{exactly}
$E$ per letter and whose entropy per letter approaches $H(E)$.
Distinct guarded cylinder intervals have positive separation.
Take $\ell_j$ divisible by $L$; this changes none of the
counting or energy asymptotics.  The resulting independent blocks
have finite variance in both information and energy, so the
strong law applies at all free-letter prefixes after observing that
a partial block contains at most $L$ symbols and $L$ is fixed.
More precisely, almost surely the energy of the first $m$ free
symbols is $Em+o(m)$.  Up to physical stage $J$ there are only
$J=o(K_J)$ prescribed return symbols, with total energy
$o(K_J)$ by \eqref{eq:vanishingenergy}.  Within a stage, at most
$O(J^p)$ further free symbols and one prescribed symbol of energy
$O(J^{q\delta})$ remain, both $o(K_J)$.  Hence the full return
symbol sequence, at \emph{every} prefix length $k$, satisfies
$\lim_{k\to\infty}\frac1k\sum_{i=1}^k\tau_i^\delta=E$ almost surely.

For the resulting Bernoulli measure, the almost-sure information
per free symbol equals the chosen block entropy divided by its length;
these values approach $H(E)$ as the alphabet and block length grow.
The almost-sure negative logarithmic contraction per free symbol
equals $\kappa E+\log Z(\kappa)$.  The prescribed digits contribute
$o(K_J)$ by \eqref{eq:vanishingenergy}; the maximal one-step
contraction is likewise $o(K_J)$.  The weighted insertion estimate \eqref{eq:weighted-insertion},
restricted to this full-measure set, therefore gives the Euclidean
lower local dimension at least the free entropy divided by
$\kappa E+\log Z(\kappa)$.  Letting the finite guarded block laws
approach the entropy limit in Lemma~\ref{lem:finite-exact-energy}
gives dimension arbitrarily close to
$H(E)/(\kappa E+\log Z(\kappa))$.

It remains to verify the complete macroscopic spectrum.  In every
shell of physical size $R$, the occupied free-return blocks have
$\asymp R^u$ components of length $\asymp R^v$, separated by gaps
$\asymp R^{q/(q+1)}$.  Since $0<v<q/(q+1)$ and free roof
values are bounded, the three-regime counting-measure estimate
\eqref{eq:mass3regimes} applies with $v=p/(q+1)$.
The general form of Lemma~\ref{lem:powerprofile}, applied with
$p=b(q+1)-1$, gives exactly \eqref{eq:intrinsic-profile}
for every itinerary, including its value at $\theta=0$.  This establishes the lower bound with all
constraints imposed.

Finally, the entropy variational inequality
\eqref{eq:entropy-variational} and the relative-entropy identity
\[
 \kappa E+\log Z(\kappa)-H(E)
 =\sum_{n\ge1}P_\lambda(n)
     \log\frac{P_\lambda(n)}{r_n}\ge0,
 \qquad P_\lambda(n)=Z(\lambda)^{-1}e^{-\lambda n^\delta},
\]
show that $\mathfrak d_\kappa(E)\le1$ with equality only when
$\lambda(E)=\kappa$.  The endpoint limits follow from
$H(E)\to0$ as $E\downarrow1$ and
$H(E)=O(\log(2+E))$ as $E\to\infty$.
Continuity yields the remaining assertions.
\end{proof}

\begin{remark}\label{rem:intrinsicnumerics}
For an explicit illustration take $\delta=1/2$, $\kappa=1$,
$b=3/4$, and $q=2$, so $p=5/4$.  Then, for every fixed
$c_0>0$, the complete physical-time spectrum is
\[
 h(0)=\frac47,\qquad
 h(\theta)=\max\left\{\frac47,
 \left(1-\frac1{4\theta}\right)_+\right\}
 \quad(0<\theta\le1).
\]
It is independent of the intrinsic energy $E$.  The formula
\eqref{eq:intrinsicpressure-theorem} gives, approximately,
$\mathfrak d_1(1.5)=0.860$,
$\mathfrak d_1(2)=0.988$, and
$\mathfrak d_1(4)=0.892$, while its unique maximal value $1$
occurs at $E=\mathcal E(1)\approx2.313$.  Hence the same
exact power-law counting coefficient and the same macroscopic
spectrum admit both nontrivial intermediate and full-dimensional
intrinsic-energy levels.
\end{remark}

\begin{remark}\label{rem:energyintrinsic}
The threshold $b>\delta$ in Theorem~\ref{thm:intrinsicpressure}
is necessary: Proposition~\ref{prop:exactcount-energy} shows that,
at and below the critical exponent, an exact power count is
incompatible with finite mean energy.  The unrestricted critical
counting level has dimension zero by Theorem~\ref{thm:exactcritical}.
\end{remark}

\subsection{Prescribed schedules with nonzero insertion energy}
We now fix a particular insertion schedule whose prescribed
branches contribute a positive logarithmic contraction per free
symbol.  The pressure equation below computes the Hausdorff
dimension of the \emph{resulting Cantor subset}; it is not another
formula for the \emph{entire} joint energy level, whose dimension
was determined in Theorem~\ref{thm:intrinsicpressure}.

We fix two nonadjacent free symbols, for example $\mathcal A=\{2,4\}$.
The corresponding fundamental intervals have a positive relative
gap because $I_3$ separates $I_2$ from $I_4$.

Fix $\delta<b<1$ and write
\[ q=\frac{1-b}{b-\delta}>0,\quad
 u=\frac1{q+1},\quad v=\frac{q\delta}{q+1},\quad
 w=\frac q{q+1},\quad b=u+v. \]
For positive constants $A,B$ choose integers
\begin{equation}\label{eq:bridge-design}
 n_j=\lfloor A j^q\rfloor\vee5,
 \qquad \ell_j=\lfloor B j^{q\delta}\rfloor\vee1.
\end{equation}
At stage $j$, place $\ell_j$ freely chosen symbols from
$\mathcal A$, followed by the single prescribed symbol $n_j$.
The construction retains every sequence of free choices.  Its image
under the affine coding map is denoted $K(A,B)$.

The physical return-counting asymptotic and the cumulative contraction
of the prescribed branches are uniform over all free choices.
\begin{lemma}\label{lem:powerbridge}
For every point of $K(A,B)$, the physical return-counting function
satisfies
\begin{equation}\label{eq:exactbcount}
 N_y(t)\sim C(A,B)t^b,\qquad
 C(A,B)=\frac{B}{q\delta+1}
             \left(\frac{q+1}{A}\right)^b.
\end{equation}
If $k$ counts free symbols, then the cumulative negative logarithmic
contraction $G(k)$ of the prescribed $n_j$-branches prior to the $k$th
free symbol satisfies
\begin{equation}\label{eq:gamma-rate}
 \frac{G(k)}k\longrightarrow
 \gamma(A,B):=\frac{\kappa A^\delta}{B},
 \qquad
 \frac{\max\{-\log r_{n_j}:n_j\text{ inserted before }k\}}k
 \longrightarrow0.
\end{equation}
The same $o(k)$ bound holds for the contraction of any one inserted
branch near the $k$th free symbol.
\end{lemma}
\begin{proof}
Let $K_J=\sum_{j\le J}\ell_j$ and let $T_J$ be the
physical time elapsed after the $J$th prescribed excursion.  The
bounded free roof values $2$ and $4$ give, uniformly in the free
choices,
\[
 K_J\sim\frac{B}{q\delta+1}J^{q\delta+1},\qquad
 T_J=\sum_{j\le J}n_j+O(K_J)
 \sim\frac A{q+1}J^{q+1}.
\]
Here $q\delta+1<q+1$, so the free roof contribution is lower
order.  Moreover $n_J+\ell_J=o(T_J)$ and
$\ell_J=o(K_J)$.  These estimates interpolate the two displayed
asymptotics to every intermediate physical time, proving
\eqref{eq:exactbcount}.
Since $-\log r_n=\kappa n^\delta+\log Z$,
\[
 G(K_J)=\sum_{j\le J}(\kappa n_j^\delta+\log Z)
 \sim\frac{\kappa A^\delta}{q\delta+1}J^{q\delta+1}.
\]
The largest individual insertion has size $O(J^{q\delta}+1)$,
which is $o(K_J)$.  The estimates continue to hold at every
intermediate free-symbol prefix, since the number of free symbols
in one stage is $o(K_J)$.
\end{proof}

Lemma~\ref{lem:powerprofile}, with $p=q\delta$, gives the
common spectrum of $K(A,B)$ as
$h_{q,\delta}=\mathfrak p_{q\delta,q}$, namely
\begin{equation}\label{eq:powerprofile}
 h_{q,\delta}(0)=d_0,\qquad
 h_{q,\delta}(\theta)=
 \max\left\{d_0,\left(1-\frac{1-b}{\theta}\right)_+\right\}
 \quad(0<\theta\le1),
\end{equation}
where
\[
 b=\frac{1+q\delta}{1+q},\qquad
 d_0=\frac1{1+q(1-\delta)}.
\]

\begin{proposition}\label{prop:prescribed-schedule}
Consider the affine return tower \eqref{eq:affineratio2}, and fix
$\delta<b<1$, $c_0>0$ and $\gamma>0$.  There exists a Cantor set
$K_{b,c_0,\gamma}$ of initial base points such that every
$y\in K_{b,c_0,\gamma}$ satisfies
\[ N_y(t)\sim c_0t^b,\qquad
 \mathsf D^-_\theta(R(y,Y))=
 \mathsf D^+_\theta(R(y,Y))=h_{q,\delta}(\theta)
 \quad\text{for every }\theta\in[0,1], \]
where $q=(1-b)/(b-\delta)$ and $h_{q,\delta}$ is given by
\eqref{eq:powerprofile}.  The prescribed-branch geometric loss per
free symbol converges to $\gamma$, and
\begin{equation}\label{eq:pressure-intermediate}
 \dim_H K_{b,c_0,\gamma}=s_\gamma,
 \qquad r_2^{s_\gamma}+r_4^{s_\gamma}
              =e^{\gamma s_\gamma}.
\end{equation}
The solution $s_\gamma$ is unique and belongs to $(0,1)$.
The function $\gamma\mapsto s_\gamma$ is continuous and strictly
decreasing, with
$\lim_{\gamma\to\infty}s_\gamma=0$ and
$\lim_{\gamma\downarrow0}s_\gamma=s_0^{\rm free}$,
where $r_2^{s_0^{\rm free}}+r_4^{s_0^{\rm free}}=1$.
These dimensions describe the specified Cantor subsets
$K_{b,c_0,\gamma}$, with $b,c_0$ and the complete spectrum fixed.
\end{proposition}
\begin{proof}
Use the bridge construction \eqref{eq:bridge-design}.  The two
conditions $C(A,B)=c_0$ and $\gamma(A,B)=\gamma$ can be solved
explicitly: take
\[ A=\left(
 \frac{\kappa(q+1)^b}{\gamma(q\delta+1)c_0}
 \right)^{1/(b-\delta)},\qquad
 B=\frac{\kappa A^\delta}{\gamma}. \]
Then Lemma~\ref{lem:powerbridge} yields the precise counting
asymptotic and insertion energy, and Lemma~\ref{lem:powerprofile}
gives the spectrum for all $\theta$ simultaneously.

For the Hausdorff dimension let $k$ denote the number of free
letters and write $G(k)$ for the prescribed contraction prior to
that free position.  At a free-digit cylinder of length $k$ the
Euclidean diameter, up to a fixed multiplicative constant, equals
$
 e^{-G(k)}\prod_{i=1}^k r_{a_i}$ with $ a_i\in\{2,4\}.$
The $s$-dimensional contribution of the complete family of such
cylinders is therefore
\[
 e^{-sG(k)}(r_2^s+r_4^s)^k
 =\exp\bigl(k[\log(r_2^s+r_4^s)-s\gamma+o(1)]\bigr).
\]
For $s>s_\gamma$ this tends to zero, proving the upper bound.

To prove the reverse inequality, define independent free digits
with probabilities
\[
 p_a=e^{-s_\gamma\gamma}r_a^{s_\gamma}
 \quad(a\in\{2,4\});\qquad p_2+p_4=1.
\]
The law of large numbers gives, for almost every coded point,
\[
 -\log\mu([a_1\ldots a_k])=kh+o(k),\quad
 -\log|I([a_1\ldots a_k])|=k(\chi+\gamma)+o(k),
\]
where
$h=-\sum_a p_a\log p_a$ and
$\chi=-\sum_a p_a\log r_a$.  The normalization of $p_a$ yields
$h=s_\gamma(\chi+\gamma)$.

We explain why the preceding cylinder calculation also controls
\emph{Euclidean balls}.  The distinct free first-level intervals
$I_2$ and $I_4$ are separated by the positive-length interval
$I_3$.  At the first differing free digit, the Euclidean separation
of two coded descendants is therefore at least a fixed multiple
of the length of their common similarity prefix.  The additional
prescribed contraction between successive free digits has
logarithm $o(k)$ by \eqref{eq:gamma-rate}.  Consequently, for
$\mu$-almost every point, the separation radii of its successive
free-digit cylinders have negative logarithm
$k(\chi+\gamma)+o(k)$, with no unaccounted exponential gaps
between successive radii.  At each such radius a ball can intersect
only the cylinder specified by the corresponding free prefix.
The mass distribution principle (or the lower local-dimension
criterion) therefore gives
$\underline{\dim}_{\rm loc}\mu\ge
h/(\chi+\gamma)=s_\gamma$ almost everywhere.
This proves the lower bound.

Finally, $s\mapsto\log(r_2^s+r_4^s)-\gamma s$ is strictly
decreasing, positive at $s=0$ and negative at $s=1$.
The implicit equation proves uniqueness, continuity and strict
monotonicity in $\gamma$; the endpoint limits follow directly.
\end{proof}

\begin{remark}\label{ex:intermediate-half}
Take $\delta=1/2$, $q=1$ and $\kappa=1$.  Then $b=3/4$,
$\theta_*=3/4$, and the complete macroscopic dimension spectrum is
\[
 h(0)=\frac23,\qquad
 h(\theta)=\max\left\{\frac23,
            \left(1-\frac1{4\theta}\right)_+\right\}
 \quad(0<\theta\le1).
\]
For \emph{every} fixed $c_0>0$ we obtain sets of base points with
$N_y(t)\sim c_0t^{3/4}$ and exactly this function $h$, while their
Hausdorff dimensions vary continuously with $\gamma$.  For example,
with the normalized ratios \eqref{eq:affineratio2},
$\gamma=0.2$ gives $s_\gamma\approx0.288$, whereas
$\gamma=1$ gives $s_\gamma\approx0.216$.
The numbers here are dimensions of the explicitly constructed
Cantor subsets, not an asserted formula for the full unconstrained
level set.
\end{remark}

\begin{remark}\label{rem:generalfree-pressure}
The free alphabet $\{2,4\}$ is chosen only to make Euclidean
separation transparent.  For a finite family $\mathcal W$ of
well-separated free words, write $r_w$ for the affine contraction
of word $w$.  The identical proof gives the unique pressure equation
$\sum_{w\in\mathcal W}r_w^s=e^{s\gamma}$, where $\gamma$ is measured
per free word.  This finite-alphabet pressure computation is consistent with
the theory of non-autonomous Moran constructions
\cite{HollandZhang,RempeUrbanski}.
\end{remark}

\section{Conditional pressure and nonlinear induced geometry}\label{sec:nonlinearpressure}

The exact entropy-over-contraction formula in Theorem~\ref{thm:intrinsicpressure}
uses the similarity character of the affine branches.  We now show which
parts persist for genuinely nonlinear induced systems and formulate an
exact conditional pressure principle under explicit thermodynamic
hypotheses.  In particular, uniform bounded distortion alone does
\emph{not} justify replacing the geometric Lyapunov exponent by a
branchwise return-energy average.

\subsection{Geometry and constrained pressure}
Let $Y$ be a compact interval and let $(\psi_n)_{n\ge1}$ be $C^{1+\beta}$
orientation-preserving contractions taking $Y$ into intervals with
pairwise disjoint interiors.  Assume the full-branch coding is one-to-one
outside a countable boundary set and that all finite compositions satisfy
uniform bounded distortion.  Let $\Sigma=\N^\N$, $\sigma$ denote the
left shift and $\pi:\Sigma\to Y$ the coding map.  As before, the symbol
$n$ is the \emph{actual first return length}.  Write
\[
 \Phi(n)=n^\delta,\qquad
 \chi(\omega)=-\log\left|\psi_{\omega_1}'
       (\pi(\sigma\omega))\right|,\quad 0<\delta<1.
\]
The chain rule gives the identity
$ -\log|\psi_{\omega_1\cdots\omega_k}'
        (\pi(\sigma^k\omega))|=S_k\chi(\omega). $
Our main geometric comparison assumption is that, for $\kappa>0$ and
$\eta\ge0$, the derivative of every branch satisfies
\begin{equation}\label{eq:nonlinear-perturb}
 e^{-\eta}r_n\le |\psi_n'(y)|\le e^{\eta}r_n,
 \qquad r_n=Z(\kappa)^{-1}e^{-\kappa n^\delta},\quad
 Z(t)=\sum_{n\ge1}e^{-t n^\delta}.
\end{equation}
We allow $\eta$ to be large in the pressure statements, but require
$\sup_n\|\psi_n'\|_\infty<1$.  Notice that the common geometric
potential $\chi$ depends on the full future itinerary, even though
$\Phi$ depends only on the current first-return symbol.

For a word $w=(n_1,\ldots,n_k)$ put
$M_w=\sup_{y\in Y}|\psi_w'(y)|$ and $\Phi(w)=\sum_{i=1}^k n_i^\delta$.
For $s>0$ and $t\in\R$ define the constrained upper pressure
\[ \mathcal P(s,t;E)=\limsup_{k\to\infty}\frac1k
 \log\sum_{|w|=k}M_w^s
       \exp\bigl(t[\Phi(w)-kE]\bigr). \]
Whenever $t<s\kappa$, the series is finite under
\eqref{eq:nonlinear-perturb} and bounded distortion.  One can
alternatively use the Gurevi\v c pressure of the potential
$-s\chi+t(\Phi-E)$; the two pressures agree because distortion
changes the partition sums by a bounded multiplicative factor.

The energy-constrained geometric pressure provides an upper bound without an equilibrium-state assumption.
\begin{proposition}\label{prop:nonlinear-upper}
For any $E>1$, let $\mathscr A(E)$ be the set of coded points with
$\lim_{k\to\infty}k^{-1}S_k\Phi=E$.  Under the assumptions above,
\[ \dim_H\mathscr A(E)\le
 \inf\bigl\{s>0:\ \mathcal P(s,t;E)<0
                    \text{ for some }t<s\kappa\bigr\}. \]
Consequently the same upper bound holds with any additional
physical counting and macroscopic return-dimension constraints.
\end{proposition}
\begin{proof}
Fix $s>0$ and $t$ with $\mathcal P(s,t;E)<0$.  On the cylinder
family $|S_k\Phi-kE|\le\varepsilon k$,
\[
 \sum |\psi_w(Y)|^s\le |Y|^s e^{|t|\varepsilon k}
   \sum_{|w|=k}M_w^s e^{t(\Phi(w)-kE)}.
\]
For sufficiently small $\varepsilon>0$ the right side decays
exponentially.  Every point of $\mathscr A(E)$ is in these cylinder
unions at all sufficiently large depths.  Uniform contraction makes
the diameters tend to zero, so its Hausdorff $s$-measure vanishes.
\end{proof}

\subsection{A quantitative stability theorem}\label{subsec:nonlinear-stability}
Fix $\delta<b<1$, $c_0>0$, and choose
\begin{equation}\label{eq:nonlinear-qchoice}
 q>\frac{1-b}{b-\delta},\qquad
 p=b(q+1)-1,\qquad
 d_{q,b}=\frac{1}{1+(1-b)(q+1)}.
\end{equation}
Then $q\delta<p<q$.  Recall the concrete admissible dimension spectrum
\begin{equation}\label{eq:nonlinear-h}
 h_{q,b}(0)=d_{q,b},\qquad
 h_{q,b}(\theta)=\max\left\{d_{q,b},
             \left(1-\frac{1-b}{\theta}\right)_+\right\}
             \quad(0<\theta\le1).
\end{equation}
Let $\mathscr L_{\psi}(c_0,b,E;q)$ denote the \emph{entire}
level set determined by $N_y(t)\sim c_0t^b$, the exact
$\Phi$-average $E$, and the complete spectrum \eqref{eq:nonlinear-h}.
Put
\[
 H(E)=\lambda(E)E+\log Z(\lambda(E)),\qquad
 L(E)=\kappa E+\log Z(\kappa),
\]
where $\lambda(E)>0$ is the unique solution of
$-Z'(\lambda)/Z(\lambda)=E$.

The following theorem gives both a genuinely nonlinear quantitative
statement and a robust test for the affine exact formula.
\begin{theorem}\label{thm:nonlinear-stability}
Assume \eqref{eq:nonlinear-perturb}, the uniform bounded distortion
hypothesis, and the existence of an interior guard branch.  Let
$E>1$, $\delta<b<1$, $c_0>0$ and let $q$ be as in
\eqref{eq:nonlinear-qchoice}.  Then
\begin{equation}\label{eq:nonlinear-stability-bound}
 \frac{H(E)}{L(E)+\eta}\le
 \dim_H\mathscr L_{\psi}(c_0,b,E;q)
 \le\min\left\{1,\frac{H(E)}{L(E)-\eta}\right\}
 \qquad\text{if }L(E)>\eta.
\end{equation}
The lower bound remains valid without the restriction $L(E)>\eta$.
In particular, for a family satisfying the assumptions uniformly
with $\eta\downarrow0$, the dimensions of these complete level sets
converge to $H(E)/L(E)$.
\end{theorem}
\begin{proof}
We prove the upper and lower bounds in
\eqref{eq:nonlinear-stability-bound} separately.
For the upper bound let $\varepsilon>0$.  On every energy-window
word of length $k$ we have by \eqref{eq:nonlinear-perturb}
\[
 |\psi_w(Y)|\le |Y|\prod_{i=1}^k e^\eta r_{n_i}
 \le |Y|\exp\{-k[L(E)-\eta-\kappa\varepsilon]\}.
\]
The entropy inequality for the countable alphabet gives
\[
 \#\{w:|\Phi(w)/k-E|\le\varepsilon\}
 \le\exp\{k[H(E)+\lambda(E)\varepsilon]\}.
\]
The $s$-dimensional sum therefore decays whenever
$s>[H(E)+\lambda(E)\varepsilon]/
    [L(E)-\eta-\kappa\varepsilon]$.
Taking $\varepsilon\downarrow0$ proves the upper bound.

For the lower bound, choose a finite law on $L$-letter guarded
blocks from Lemma~\ref{lem:finite-exact-energy}.  Its mean energy
per letter is exactly $E$, and its entropy per letter approaches
$H(E)$.  For each such law, \eqref{eq:nonlinear-perturb} gives
\[
 \frac1L\sum_Bp_B[-\log\inf_Y|\psi_B'|]
 \le \kappa E+\log Z(\kappa)+\eta=L(E)+\eta.
\]
Use the polynomial schedule from the proof of
Theorem~\ref{thm:intrinsicpressure}: place
$\ell_j=L\lfloor Bj^p/L\rfloor$ free letters in whole blocks
before $n_j=\lfloor Aj^q\rfloor\vee2$, with
\[
 \frac{B}{p+1}\left(\frac{q+1}{A}\right)^b=c_0.
\]
The exact count follows uniformly over free choices, and
Lemma~\ref{lem:powerprofile} gives \eqref{eq:nonlinear-h}.
The comparison \eqref{eq:nonlinear-perturb} bounds every prescribed
prefix by the same estimate \eqref{eq:power-prefix-penalty}, with
$\rho=(q\delta+1)/(p+1)<1$.  Both the number and total energy
of prescribed symbols are $o(k)$ at free-symbol depth $k$.
The strong law for each fixed block law therefore gives mean energy
$E$ along the full itinerary, including prefixes inside stages.

The weighted insertion estimate \eqref{eq:weighted-insertion}
applies to this full-measure set, using the positive separation
of the fixed finite guarded family.  It yields
\[
 \dim_H\mathscr L_{\psi}(c_0,b,E;q)
 \ge\frac{H(p)/L}{L(E)+\eta}.
\]
Passing through the finite block laws of
Lemma~\ref{lem:finite-exact-energy} proves the lower bound.
\end{proof}

\subsection{Exact pressure under thermodynamic regularity}
The bounds above require no equilibrium-state theory.  To obtain an
\emph{equality} in a nonlinear model, the pressure must retain the
full geometric potential $\chi$ rather than only the function
$\Phi(n)$.  The relevant background is the countable-state
thermodynamic formalism of Sarig \cite{Sarig} and the multifractal
variational methods of Iommi and Jordan \cite{IJ,IJQuotients}.

We formulate the additional, independently checkable assumption
that allows one to invoke this machinery.  Write
\[
 Q(s,t)=P_{\rm G}(-s\chi+t\Phi),\qquad s>0,\quad t<s\kappa.
\]
\emph{Thermodynamic regularity (TR).} We require the following.
\begin{enumerate}[label=\textup{(TR\arabic*)},leftmargin=*]
\item Throughout $s>0$, $t<s\kappa$, the pressure is finite,
differentiable, and represented by a mixing Gibbs equilibrium state
$\mu_{s,t}$, with
$\partial_tQ=\int\Phi\,d\mu_{s,t}$ and
$\partial_sQ=-\int\chi\,d\mu_{s,t}$.
\item For each $s>0$, the value of $\partial_tQ(s,t)$ ranges
continuously and strictly increasingly over $(1,\infty)$ as
$t$ traverses $(-\infty,s\kappa)$.
\item The constrained equilibria admit approximation in entropy,
geometric exponent, and energy by ergodic finite-alphabet measures.
The mean energy $E$ can be obtained exactly by arbitrarily small
finite-block adjustments.
\end{enumerate}

These hypotheses are stronger than uniform bounded distortion; in
particular, we do \emph{not} infer (TR) from bounded distortion alone.
They hold in the affine model and in smooth conjugates of that
model, and are natural to verify for sufficiently regular summable
H\"older perturbations using countable-state transfer operators.

The pressure zero below gives the dimension of the intrinsic
energy level.  The polynomial insertion construction shows that this
dimension persists after both physical-time constraints are imposed.

\begin{theorem}\label{thm:nonlinear-conditional}
Assume \eqref{eq:nonlinear-perturb}, uniform bounded distortion,
interior guarded separation and (TR).  Given $E>1$, define
\begin{equation}\label{eq:nonlinear-legendre}
 \mathcal F(s,E)=\inf_{t<s\kappa}
                 \bigl[Q(s,t)-tE\bigr].
\end{equation}
There is a unique $s_E\in(0,1]$ with
\begin{equation}\label{eq:nonlinear-root}
 \mathcal F(s_E,E)=0.
\end{equation}
For every $\delta<b<1$, $c_0>0$, and $q$ satisfying
\eqref{eq:nonlinear-qchoice}, the entire constrained level set has
\[ \dim_H\mathscr L_{\psi}(c_0,b,E;q)=s_E. \]
If $t_E$ minimizes \eqref{eq:nonlinear-legendre} at $s=s_E$,
then
\begin{equation}\label{eq:nonlinear-entropy-ratio}
 \int\Phi\,d\mu_{s_E,t_E}=E,\qquad
 s_E=\frac{h_{\mu_{s_E,t_E}}(\sigma)}
             {\int\chi\,d\mu_{s_E,t_E}}.
\end{equation}
In the affine case, \eqref{eq:nonlinear-root} reduces to the
explicit formula \eqref{eq:intrinsic-dim}.
\end{theorem}
\begin{proof}
By (TR), for each $s>0$ the infimum in
\eqref{eq:nonlinear-legendre} is attained at the unique $t=t(s,E)$
with $\int\Phi\,d\mu_{s,t}=E$.  The pressure variational identity
gives
\[
 \mathcal F(s,E)=
 h_{\mu_{s,t}}(\sigma)-s\int\chi\,d\mu_{s,t}.
\]
Here are details at the pressure endpoints.  Put
$L(E)=\kappa E+\log Z(\kappa)$ and let $\eta$ be as in
\eqref{eq:nonlinear-perturb}.  The two-sided derivative comparison,
applied to every finite word, yields
\[
 \left|Q(s,t)-\bigl[\log Z(s\kappa-t)
                 -s\log Z(\kappa)\bigr]\right|\le s\eta.
\]
Taking the infimum over $t<s\kappa$, or equivalently over
$\alpha=s\kappa-t>0$, and using the Gibbs entropy variational
identity \eqref{eq:entropy-variational}, we obtain
\begin{equation}\label{eq:pressure-sandwich}
 H(E)-s[L(E)+\eta]\le\mathcal F(s,E)
                    \le H(E)-s[L(E)-\eta].
\end{equation}
The comparison \eqref{eq:pressure-sandwich} in particular gives
$\lim_{s\downarrow0}\mathcal F(s,E)=H(E)>0$.
If $s_2>s_1$, the uniform contraction estimate
$\chi\ge\chi_{\min}>0$ gives
$Q(s_2,t)\le Q(s_1,t)-(s_2-s_1)\chi_{\min}$;
As the domain $t<s_1\kappa$ is contained in
$t<s_2\kappa$, taking the infimum shows that
$\mathcal F(s,E)$ is strictly decreasing with rate at least
$\chi_{\min}$.
Its continuity follows from (TR) and the unique, continuously
varying minimizer $t(s,E)$ (equivalently, from the envelope
identity $\partial_s\mathcal F=-\int\chi\,d\mu_{s,t(s,E)}$).
The uniqueness of $s_E$ therefore follows from strict monotonicity,
not merely from the existence of a pressure zero.
Thus there is a unique positive zero.  Finally, the depth-$k$
branch intervals are disjoint, and bounded distortion gives
$\sum_{|w|=k}M_w\le e^D$ after normalization of $|Y|$.
Consequently $Q(1,0)\le0$, hence
$\mathcal F(1,E)\le Q(1,0)\le0$ and the zero lies in $(0,1]$.
At the zero the pressure identity gives
\eqref{eq:nonlinear-entropy-ratio}.

For the upper bound, if $s>s_E$ then
$\mathcal F(s,E)<0$.  Selecting a minimizing $t$ yields
$\mathcal P(s,t;E)<0$ by bounded distortion, so
Proposition~\ref{prop:nonlinear-upper} gives
$\dim_H\mathscr L_{\psi}\le s_E$.

For the lower bound, (TR) supplies successively richer finite
laws $(p_B)$ on guarded blocks of common length $L$, with mean
energy exactly $E$, such that
\[
 \frac{H(p)}L\longrightarrow h_{\mu_{s_E,t_E}}(\sigma),
 \qquad
 \frac1L\sum_Bp_B[-\log\inf_Y|\psi_B'|]
 \longrightarrow\int\chi\,d\mu_{s_E,t_E}.
\]
Bounded distortion contributes $O(D/L)$; the guard and the
exact-energy adjustment contribute $o(1)$ per letter.  Each
finite block law is fixed before applying the dimension estimate.
In particular, its positive separation constant may depend on
the alphabet and $L$.

Insert the same polynomial schedule as in
Theorem~\ref{thm:nonlinear-stability}, with its constants $A,B$
chosen to give coefficient $c_0$.  If $k$ counts completed free
blocks, the stage before the next block is
$O_{M,L}((1+k)^{1/(p+1)})$.  The derivative bound in
\eqref{eq:nonlinear-perturb} therefore gives, at every prefix,
\begin{equation}\label{eq:nonlinear-prefix-penalty}
 A(k)\le C_{M,L}(1+k)^\rho,
 \qquad \rho=\frac{q\delta+1}{p+1}<1.
\end{equation}
The number and energy of inserted symbols are likewise negligible.
The strong law for independent blocks of the fixed law gives mean
energy $E$ along the full return itinerary.  The counting law is
uniform in the free choices, and Lemma~\ref{lem:powerprofile}
supplies the entire spectrum $h_{q,b}$.

All hypotheses of the weighted insertion estimate
\eqref{eq:weighted-insertion} now hold: finite guarded separation,
causal insertions, and a uniform prefix penalty $o(k)$.
For the full-measure set satisfying the energy law, it gives
\[
 \dim_H\mathscr L_{\psi}(c_0,b,E;q)
 \ge \frac{H(p)}{\sum_Bp_B[-\log\inf_Y|\psi_B'|]}.
\]
Passing through the finite approximations gives
$h_{\mu_{s_E,t_E}}(\sigma)/\int\chi\,d\mu_{s_E,t_E}=s_E$.
Together with the upper bound, this proves the formula for the
entire constrained level.
\end{proof}

\begin{corollary}\label{cor:nonlinear-conjugacy}
Let $(\phi_n)$ be the affine branches with ratios $r_n$ and let
$g:Y\to Y$ be an orientation-preserving $C^{1+\beta}$
diffeomorphism with $0<\inf|g'|\le\sup|g'|<\infty$.
Set $\psi_n=g\circ\phi_n\circ g^{-1}$, keeping the
first-return roof $\tau(n)=n$.  The $\psi_n$ are generally nonlinear,
with uniform bounded distortion, and for all $E>1$, $\delta<b<1$,
$c_0>0$, and $q$ as above,
\[
 \dim_H\mathscr L_{\psi}(c_0,b,E;q)
              =\frac{H(E)}{L(E)}.
\]
\end{corollary}
\begin{proof}
The geometric coding of the nonlinear model is obtained from the
affine one by $g$, which is bi-Lipschitz on $Y$.  It changes neither
the symbolic return times nor their physical counting and macroscopic
spectrum, and preserves Hausdorff dimension.  Equivalently, the
chain rule gives the bounded telescoping identity
\[
 S_k\chi_{\psi}(\omega)
 =k\log Z(\kappa)+\kappa S_k\Phi(\omega)
    +\log|g'(\pi_{\rm aff}(\sigma^k\omega))|
    -\log|g'(\pi_{\rm aff}(\omega))|.
\]
Indeed, writing $z_i=\pi_{\rm aff}(\sigma^i\omega)$, we have
$\psi_{\omega_1\cdots\omega_k}'(g(z_k))
 =g'(z_0)(\prod_{i=1}^k r_{\omega_i})/g'(z_k)$,
which confirms the signs in the identity.
Thus the geometric correction is a bounded
coboundary and the affine pressure formula remains exact.
\end{proof}

\subsection{A genuinely nonlinear model with an exact pressure formula}
To interpret the roof labels as \emph{genuine first-return times},
we first give the underlying dynamical realization.  This avoids
identifying a labelled iterated function system with a first-return
map without exhibiting the ambient dynamics.
\begin{lemma}\label{lem:nonlinear-return-tower}
Let $Y$ be a compact interval partitioned, up to endpoints, into
pairwise interior-disjoint intervals $I_n=\psi_n(Y)$, $n\ge1$,
where each $\psi_n:Y\to I_n$ is a homeomorphism.  Put
$F|_{\operatorname{int}I_n}=\psi_n^{-1}$ and
$\tau|_{\operatorname{int}I_n}=n$.  Then there is a countable-level
return-tower transformation $T:X\to X$ whose first-return map to
its base $Y$ is $F$ and whose first-return time at $y\in I_n$ is
exactly $n$.  The geometry and Hausdorff dimension on the base are
the original Euclidean ones.
\end{lemma}
\begin{proof}
Let $\mathcal B$ be the union of the partition boundaries and
all their iterated inverse images under $F$.  This is countable,
and $Y^*=Y\setminus\mathcal B$ is forward invariant under $F$.
Write $I_n^*=\operatorname{int}I_n\cap Y^*$ and define
\[
 X=\bigsqcup_{n\ge1}\bigsqcup_{\ell=0}^{n-1}
       (I_n^*\times\{\ell\}),
 \qquad Y_0=\bigsqcup_{n\ge1}(I_n^*\times\{0\}).
\]
Identify $Y_0$ with $Y^*$ via $(y,0)\mapsto y$, and put
\[
 T(y,\ell)=
 \begin{cases}
 (y,\ell+1),&0\le\ell<n-1,\quad y\in I_n,\\
 (F(y),0),&\ell=n-1,\quad y\in I_n.
 \end{cases}
\]
The intermediate levels are disjoint from the base, and hence
$T^k(y,0)\notin Y_0$ for $1\le k<n$ while
$T^n(y,0)=(F(y),0)\in Y_0$.  Iterating gives exactly the
successive roof symbols of the $F$-coding.  The base is measured
with its inherited interval metric.  This is a genuine first-return
\emph{tower}; the construction does not claim that $T$ is a smooth
single-interval map.
\end{proof}

The thermodynamic hypothesis (TR) can be checked in a natural
non-conjugate class, rather than imposed abstractly.  Let
$a_n=\sum_{m<n}r_m$, so that the intervals
$I_n=[a_n,a_n+r_n]$ partition $[0,1]$ up to endpoints.  Choose
$\varepsilon>0$ satisfying $0<\varepsilon<\min\{1,r_1^{-1}-1\}$ and put
\begin{equation}\label{eq:branchwise-nonlinear}
 \psi_n(x)=a_n+r_n\bigl[x+\varepsilon(-1)^n x(1-x)\bigr],
 \qquad 0\le x\le1.
\end{equation}
These are smooth full inverse branches with
$\eta=\max\{-\log(1-\varepsilon),\log(1+\varepsilon)\}$
in \eqref{eq:nonlinear-perturb}.  Their branchwise geometry depends
nontrivially on the future coordinate $x$.

\begin{proposition}\label{prop:nonlinear-example}
For the family \eqref{eq:branchwise-nonlinear}, with roof
$\tau(n)=n$ and $\delta\in(0,1)$, the associated return tower
of Lemma~\ref{lem:nonlinear-return-tower} has genuine first returns;
its countable full shift has uniform bounded distortion and
satisfies (TR).  Consequently,
for every $E>1$, $\delta<b<1$, $c_0>0$ and $q$ as in
\eqref{eq:nonlinear-qchoice}, the complete constrained level set
has the \emph{exact} dimension $s_E$ specified by
\eqref{eq:nonlinear-root}.  If $\varepsilon\ne0$, the system is not
$C^1$-conjugate, by a single global branch-label-preserving
conjugacy, to the affine
branch family $(x\mapsto a_n+r_nx)_n$.
\end{proposition}
\begin{proof}
Put $\lambda=r_1(1+\varepsilon)<1$ and
$D_0=2\varepsilon/(1-\varepsilon)$.  Since
\[
 \psi_n'(x)=r_n[1+\varepsilon(-1)^n(1-2x)],\qquad
 \operatorname{Lip}(\log\psi_n')\le D_0,
\]
the variation of $\log|\psi_w'|$ is at most
$D_0/(1-\lambda)$ for every word $w$.
For $s>0$ and $t<s\kappa$, consider on Lipschitz functions the operator
\[
 (\mathcal L_{s,t}f)(x)=\sum_{n\ge1}
       e^{tn^\delta}|\psi_n'(x)|^s f(\psi_n(x)),\qquad
 \alpha=s\kappa-t>0.
\]
The series and all parameter derivatives converge locally in Lipschitz
operator norm, since their terms are bounded by a fixed polynomial in
$n^\delta$ times $e^{-\alpha n^\delta}$.  We record the spectral and
uniformity arguments needed for (TR).

\emph{Positive eigenfunction and spectral gap.}
The cone of positive functions normalized at $x_0$ and satisfying
$\operatorname{Lip}(\log h)\le B_s:=sD_0/(1-\lambda)$ is preserved by
$f\mapsto\mathcal L_{s,t}f/\mathcal L_{s,t}f(x_0)$.
Indeed the logarithmic variation of each summand is at most
$sD_0+\lambda B_s=B_s$.  The normalized cone is convex, compact in
the uniform norm, and bounded away from zero.  Schauder's theorem
therefore supplies a positive Lipschitz eigenfunction $h$ and
$\varrho>0$ with $\mathcal L_{s,t}h=\varrho h$.  Its ratio bounds
depend on $s$ but are independent of $t<s\kappa$.
Normalize the operator by
\[
 Pf=\frac{\mathcal L_{s,t}(hf)}{\varrho h}
   =\sum_n p_n(x)f(\psi_n(x)),\qquad \sum_np_n(x)=1.
\]
The logarithmic Lipschitz constants of $p_n$ are uniformly bounded
in $n$; more precisely,
$\operatorname{Lip}(p_n)\le C_{s,t}\sup_xp_n(x)$ and
$\sum_n\operatorname{Lip}(p_n)<\infty$.
Indeed $\log p_n=s\log|\psi_n'|+\log h\circ\psi_n
-\log h-\log\varrho+t n^\delta$, whose spatial Lipschitz
constant is bounded independently of $n$.  The kernel therefore
preserves Lipschitz functions.  Iteration gives the
Lasota--Yorke estimate
\[
 \operatorname{Lip}(P^k f)\le\lambda^k\operatorname{Lip}(f)
                   +C\|f\|_\infty,
 \qquad \|P^kf\|_\infty\le\|f\|_\infty.
\]
To justify the spectral statement, equip
$\operatorname{Lip}([0,1])$ with
$\|f\|_{\operatorname{Lip}}=\|f\|_\infty+\operatorname{Lip}(f)$
and use $C([0,1])$ with the uniform norm as the weak space.
The unit ball of the former is relatively compact in the latter
by Arzel\`a--Ascoli.  The two displayed inequalities are precisely
the Doeblin--Fortet estimates for these strong and weak norms.
The essential-spectral-radius theorem for Lipschitz kernels
\cite{Hennion} therefore gives
$r_{\mathrm{ess}}(P|_{\operatorname{Lip}})\le\lambda<1$.
Since $P\mathbf 1=\mathbf 1$ and $\|P^k\|_{C\to C}\le1$, the
spectral radius is one.  In particular, the spectrum outside the
closed disk of radius $\lambda$ consists of isolated eigenvalues
of finite multiplicity.
There are no other peripheral eigenvalues.  If $Pf=zf$ with $|z|=1$,
choose a point maximizing $|f|$.  Equality in the triangle inequality
forces $f(\psi_w x)=z^k f(x)$ for every word of length $k$.
Here take $x_0$ to be a point where $|f|$ attains its
maximum.  Since the depth-$k$ branch intervals cover $[0,1]$
and have diameter at most $\lambda^k$, for each $x\in[0,1]$
there is a word $w$ of length $k$ such that
$|x-\psi_w(x_0)|\le\lambda^k$.
Equality in the convex combination defining $P^k f(x_0)$
gives $f(\psi_w(x_0))=z^k f(x_0)$ for every such word.
Uniform continuity then implies
$f(x)=\lim_{k\to\infty}z^k f(x_0)$; hence $f$ is constant, and
$z=1$.  This also shows the eigenspace at one is one-dimensional.
The power bound $\|P^k\|_{C\to C}\le1$ excludes a Jordan block at
one.  Together with quasi-compactness this proves a spectral gap.
By the fixed-point theorem for probabilities on a compact space,
$P^*$ has a stationary probability $\nu$.  Every nonempty open
interval contains a sufficiently deep branch interval, and the
positive transition weights give
$\nu(U)=\int P^k\one_U\,d\nu>0$;
thus $\nu$ has full support.  The spectral gap yields convergence
of $P^kf$ to $\int f\,d\nu$ for Lipschitz $f$.  Because $F\circ\psi_n=\operatorname{id}$, stationarity
$P^*\nu=\nu$ implies $F_*\nu=\nu$.
The positivity of the transition probabilities for symbols $1$
and $2$ gives $\sup_{x,n}p_n(x)\le1-\varepsilon_{s,t}$ for some
$\varepsilon_{s,t}>0$.  Hence $\nu$ has no atoms, so the
countable partition boundaries cause no ambiguity in the coding.
Let $\mu_{s,t}$ be the coding pullback of $\nu$ to the full shift.
It is shift invariant.  We spell out why the spectral gap
implies mixing for the symbolic coding, since cylinder indicators
are not themselves Lipschitz functions on the interval.  For a
cylinder $[w]$ of length $k$, let $A_w=\psi_w(Y)$, with endpoints
understood modulo the null partition boundaries.  The function
$g_w=P^k\one_{A_w}$ is Lipschitz: in the full-branch expansion of
$P^k$, only the inverse branch corresponding to $w$ contributes,
and its normalized weight is a positive Lipschitz function.  For
any second cylinder $[v]$, writing $A_v=\psi_v(Y)$ and using the
adjoint identity for $P$ and $F$, we obtain for $m\ge0$
\[
 \begin{aligned}
 \mu_{s,t}([w]\cap\sigma^{-(k+m)}[v])=\int\one_{A_w}(\one_{A_v}\circ F^{k+m})\,d\nu
 =\int (P^m g_w)\one_{A_v}\,d\nu.
 \end{aligned}
\]
The spectral gap gives $P^m g_w\to\int g_w\,d\nu
=\nu(A_w)$ uniformly, and hence this correlation tends to
$\mu_{s,t}([w])\mu_{s,t}([v])$.  Cylinder approximation then
establishes mixing on the full symbolic sigma-algebra.
The identity
\[
 p_n(x)=\frac{e^{t n^\delta}|\psi_n'(x)|^s
                     h(\psi_n x)}{\varrho h(x)}
\]
telescopes over words.  More explicitly, stationarity and the
full-branch partition give
\[
 \begin{aligned}
 \mu_{s,t}([w])=\nu(\psi_w(Y))
 &=\int P^k\mathbf1_{\psi_w(Y)}(x)\,d\nu(x)\\
 &=\varrho^{-k}e^{t\Phi(w)}
   \int |\psi_w'(x)|^s\frac{h(\psi_w x)}{h(x)}\,d\nu(x).
 \end{aligned}
\]
The order of inverse branches in the expansion of $P^k$ is
reversed, and hence the last formula has the \emph{forward}
cylinder word $w$ in the usual coding order.
The bounds on $h$ and uniform distortion therefore yield
\[
 C_{s,t}^{-1}\le
 \frac{\mu_{s,t}([w])}{\varrho^{-k}e^{t\Phi(w)}M_w^s}
 \le C_{s,t},
\]
with a constant independent of $k$ and $w$.
In particular, the one-symbol tails decay like
$e^{-(s\kappa-t)n^\delta}$, so the Shannon entropy, geometric
exponent, and parameter derivatives are finite.
The identification of $\log\varrho$ with pressure can also
be checked directly.  Writing
$Z_k(s,t)=\sum_{|w|=k}e^{t\Phi(w)}M_w^s$, bounded distortion
and the positive upper and lower bounds on $h$ give
$
 Z_k(s,t)\asymp\mathcal L_{s,t}^{k}\one(x)
 \asymp\varrho^k,
$
with constants independent of $k$ and $x$.  Thus the logarithmic
growth rate of the full-shift partition sums is $\log\varrho$.
Together with the Gibbs estimate, integrable energy and geometric
potential, and the countable-shift variational principle
\cite{Sarig,MU}, this identifies $Q(s,t)=\log\varrho$ as
Gurevi\v c pressure and gives the equilibrium identity
\[
 Q(s,t)=h_{\mu_{s,t}}(\sigma)
          -s\int\chi\,d\mu_{s,t}+t\int\Phi\,d\mu_{s,t}.
\]
Analytic perturbation of the isolated simple eigenvalue gives
analytic pressure and
$Q_t=\int\Phi\,d\mu_{s,t}$,
$Q_s=-\int\chi\,d\mu_{s,t}$.

\emph{Energy range and strict monotonicity.}
The bounds on $h$, the derivative comparison with $r_n$, and the
eigenfunction equation give, with $C_s$ independent of $t$,
\begin{equation}\label{eq:normalized-gibbs}
 C_s^{-1}\frac{e^{-\alpha n^\delta}}{Z(\alpha)}
 \le p_n(x)\le C_s\frac{e^{-\alpha n^\delta}}{Z(\alpha)}.
\end{equation}
They also hold for one-symbol Gibbs probabilities after integrating
against the stationary probability.  As $\alpha\downarrow0$,
$Z(\alpha)\to\infty$, so every fixed finite set of symbols has
probability tending to zero; consequently $Q_t\to\infty$.
As $\alpha\to\infty$, the upper bound on
$\sum_n(n^\delta-1)p_n(x)$ tends to zero, giving $Q_t\to1$.

For strict monotonicity, analytic differentiation of the
simple leading eigenvalue gives the Green--Kubo identity
\[
 Q_{tt}(s,t)=\lim_{k\to\infty}\frac1k
 \int\bigl(S_k\Phi-kQ_t(s,t)\bigr)^2\,d\mu_{s,t}.
\]
The identity follows by differentiating the normalized transfer
operator twice and summing correlations; the preceding spectral
gap makes these correlations summable.  All required moments are
finite by \eqref{eq:normalized-gibbs}.
For completeness, strict positivity can be seen by a martingale
argument.  The centered conditional reward
$A(x)=\sum_np_n(x)(n^\delta-Q_t)$ is Lipschitz and has stationary
mean zero.  The spectral gap solves $(I-P)u=A$ with bounded Lipschitz
$u$.  Consider the stationary Markov chain that, at state $x$,
chooses $n$ with probability $p_n(x)$ and moves to $\psi_n x$.
The corrected increment
$
 \Delta(x,n)=n^\delta-Q_t+u(\psi_n x)-u(x)
$
has conditional mean $A(x)+Pu(x)-u(x)=0$.
The centered reward sum differs from the corresponding martingale
sum of $\Delta$'s by the bounded coboundary
$u(x_0)-u(x_k)$.  Orthogonality of martingale differences therefore
gives
\[
 Q_{tt}(s,t)=\int\sum_{n\ge1}p_n(x)\Delta(x,n)^2\,d\nu(x).
\]
If this variance vanished, positivity of every $p_n$ and full
support of $\nu$, together with the continuity of $\Delta(\cdot,n)$,
would imply $\Delta(x,n)=0$ for every $x$ and $n$.
Evaluating at the fixed point of $\psi_n$ would give
$n^\delta=Q_t$ for every $n$, an impossibility.
Thus $Q_{tt}>0$ and the energy ranges strictly increasingly over
$(1,\infty)$.

\emph{Finite-alphabet approximation and exact energy.}
For a compact parameter rectangle $K$ inside $s>0$,
$t<s\kappa$, put $\alpha_K=\inf_K(s\kappa-t)>0$.
For $a+b\le2$, the derivatives of the tail satisfy
\[
 \sup_{(s,t)\in K}
 \bigl\|\partial_s^a\partial_t^b
   (\mathcal L_{s,t}-\mathcal L_{s,t}^{(M)})
 \bigr\|_{\operatorname{Lip}\to\operatorname{Lip}}
 \le C_K\sum_{n>M}(1+n^{2\delta})e^{-\alpha_K n^\delta}
 \longrightarrow0.
\]
The extra logarithmic derivative of $r_n$ is
$O(1+n^\delta)$ and is absorbed in the polynomial factor.
Here is the perturbation step in detail.  Fix $(s_0,t_0)$ in
the parameter region.  By the spectral-gap argument above, the
leading eigenvalue of $\mathcal L_{s_0,t_0}$ is simple and isolated.
Choose a circle around it that excludes the rest of the spectrum.
Operator-norm convergence of the truncated operators, together with
the resolvent identity, implies convergence of their Riesz
projections on this circle.  On a sufficiently small parameter
neighborhood the circle remains in the resolvent set, so the
corresponding simple eigenvalues $\varrho_M(s,t)$ and their
spectral projections depend analytically on $(s,t)$.  The tail
bounds for derivatives through order two, or equivalently Cauchy's
formula on a slightly larger complex parameter neighborhood, show
that $\varrho_M$ and its derivatives through order two converge
uniformly there.  Since $\varrho>0$, taking logarithms yields
locally uniform convergence of $Q_M=\log\varrho_M$ and its first
two derivatives.  A finite cover gives the stated uniformity on
any prescribed compact parameter set.  All constants here are
local in $(s,t)$; no spectral-gap bound uniform up to
$t=s\kappa$ is asserted.
For fixed $(s,t)$ with $Q_t=E>1$, take $t_-<t<t_+$ in the
parameter domain with $Q_t(s,t_-)<E<Q_t(s,t_+)$.
For all large $M$, the same strict inequalities hold for the
finite-alphabet pressures.  Their energy maps are continuous and strictly increasing when
$M\ge2$: the same martingale argument used above applies to the
finite alphabet, and distinct symbols have distinct $n^\delta$.
They therefore admit a root $t_M\in(t_-,t_+)$ with
$Q_{M,t}(s,t_M)=E$.  Every limit point solves $Q_t(s,t)=E$;
strict convexity implies $t_M\to t$.
For each such $M$, the energy identity is an equality, rather than
merely convergence to $E$; this is the input needed for the later
exact-energy block construction.  The finite-alphabet Gibbs measures
$\mu_M$ then have
\[
 \int\Phi\,d\mu_M=E,\quad
 \int\chi\,d\mu_M\longrightarrow\int\chi\,d\mu_{s,t},\quad
 h_{\mu_M}(\sigma)\longrightarrow h_{\mu_{s,t}}(\sigma).
\]
For the last convergence, use the equilibrium identity
$h_{\mu_M}=Q_M+s\int\chi\,d\mu_M
-t_M\int\Phi\,d\mu_M$ and convergence of $Q_M,t_M$.

We make the exact finite-block correction explicit.  Fix such an
$M$ and an interior guard symbol $g$.  Let $\pi_{M,L}$ be the
$(L-1)$-coordinate marginal of $\mu_M$, extended to words of
length $L$ by appending $g$.  Its entropy satisfies
$
 \frac{H(\pi_{M,L})}{L}\longrightarrow h_{\mu_M}(\sigma)$ as $L\to\infty$,
by the stationary entropy-rate theorem.  The expected mean energy
of these guarded words is
$E+(g^\delta-E)/L$.  Choose fixed symbols
$n_-<E^{1/\delta}<n_+$ and enlarge the alphabet to include them
and $g$.  The constant words consisting of $L-1$ copies of
$n_-$ or $n_+$ followed by $g$ have mean energies on opposite
sides of $E$ when $L$ is sufficiently large.  Write $E_L=E+(g^\delta-E)/L$ for the original mean energy
per symbol, and let $e_{\pm,L}$ denote the mean energies of the
two constant guarded words.  Choose the sign $\pm$ so that
$E$ lies between $E_L$ and $e_{\pm,L}$, and set
\[
 \omega_L=\frac{|E-E_L|}{|e_{\pm,L}-E_L|}=O(L^{-1}).
\]
The mixture $(1-\omega_L)\pi_{M,L}+\omega_L\delta_{B_\pm}$
has mean \emph{block} energy exactly $LE$.
Its entropy differs from $H(\pi_{M,L})$ by at most
$H_2(\omega_L)+\omega_L H(\pi_{M,L})=o(L)$.
Its expected contraction changes by $O(1)$: for fixed $M$,
all involved symbols belong to a fixed finite alphabet, so the
logarithmic contraction of either alternative word is $O(L)$,
and $\omega_L=O(L^{-1})$.

The mixing weight is $O(L^{-1})$ with $M$ fixed, so this
adjustment changes the entropy and contraction by $o(L)$; the
order of limits is first $L\to\infty$ at fixed $M$, then
$M\to\infty$.  For geometric convergence, write
$\widetilde\pi_{M,L}=(1-\omega_L)\pi_{M,L}
+\omega_L\delta_{B_\pm}$ for the corrected law. Bounded
distortion of length-$L$ compositions gives
\[
 \left|\mathbb E_{\widetilde\pi_{M,L}}
     \bigl[-\log\inf_Y|\psi_B'|\bigr]
       -\mathbb E_{\mu_M}S_{L-1}\chi\right|
 \le C_g+D+C_{M,E,g},
\]
where $C_g$ is the contraction contribution of the appended guard,
$D$ is the uniform distortion constant, and $C_{M,E,g}$ is
independent of $L$. Indeed, for the fixed augmented finite alphabet,
all length-$L$ words have logarithmic contraction $O_{M,E,g}(L)$;
the correction therefore contributes only
$O_{M,E,g}(\omega_LL)=O_{M,E,g}(1)$.
Dividing by
$L$ and then letting $L\to\infty$ proves that the entropy per
letter tends to $h_{\mu_M}$, the geometric exponent per letter
tends to $\int\chi\,d\mu_M$, and the mean energy is exactly $E$.
Finally let $M\to\infty$.  This proves the finite-block
approximation clause of (TR) without replacing exact energy by an
energy window.

\medskip
\noindent\emph{Non-conjugacy.}
The fixed point $x_n$ of $\psi_n$ lies in $I_n$
and is larger than $1/2$ for all sufficiently large $n$.  Hence
$
 \psi_n'(x_n)=r_n[1+\varepsilon(-1)^n(1-2x_n)]\ne r_n.
$
A common $C^1$ diffeomorphic conjugacy to the affine branches would
preserve each fixed-point multiplier, which is impossible.
\end{proof}

\begin{corollary}\label{cor:nonlinear-analytic}
In the branchwise nonlinear model \eqref{eq:branchwise-nonlinear},
assume the inverse branches cover $Y$ up to endpoints.  The dimension
function $E\mapsto s_E$ is real analytic on $(1,\infty)$ and
\[
 0<s_E\le1,\qquad s_E=1
  \ \Longleftrightarrow\ E=E_*:=\int\Phi\,d\mu_{1,0}.
\]
For the minimizing parameter $t_E$ in
\eqref{eq:nonlinear-legendre}, one has
\begin{equation}\label{eq:nonlinear-derivative}
 \frac{ds_E}{dE}
   =-\frac{t_E}{\displaystyle\int\chi\,d\mu_{s_E,t_E}}.
\end{equation}
At the unique full-dimensional level $E_*$,
$ds_E/dE=0$ and $d^2s_E/dE^2<0$.  For every allowed
$\varepsilon$ the function satisfies $s_E\to0$ as $E\downarrow1$
and $E\to\infty$; it is strictly increasing on $(1,E_*)$ and
strictly decreasing on $(E_*,\infty)$.
Thus this genuinely nonlinear model has a nonconstant analytic
family of complete return-time levels realizing every dimension
strictly between zero and one.
\end{corollary}
\begin{proof}
Analyticity follows from analytic pressure, positive variance
$\partial_{tt}Q>0$, and $\partial_s\mathcal F=-\int\chi\,d\mu<0$.
Differentiation of $\mathcal F(s_E,E)=0$ gives
\eqref{eq:nonlinear-derivative}, because the minimizing equation
$\partial_tQ=E$ cancels the derivative of $t_E$.
As the full branches cover the interval,
$Q(1,0)=0$, by the change-of-variables identity for the
Perron--Frobenius operator.  Therefore $\mathcal F(1,E)\le0$,
with equality exactly at
$E=\partial_tQ(1,0)$ by strict convexity of pressure.
At this point $s_E=1$ and $t_E=0$.  Differentiating the equation
$\partial_tQ(s_E,t_E)=E$ gives
$t_E'=1/\partial_{tt}Q(1,0)>0$ at $E=E_*$,
since $s_E'=0$.  The derivative identity then implies
$s_E''=-1/[\partial_{tt}Q(1,0)\int\chi\,d\mu_{1,0}]<0$.
Finally, uniform contraction gives
$\chi_{\min}:=-\log(r_1(1+\varepsilon))>0$.  The entropy
bound at fixed energy implies $s_E\le H(E)/\chi_{\min}\to0$
as $E\downarrow1$.  For large $E$, the stability upper bound gives
$s_E\le H(E)/[\kappa E+\log Z(\kappa)-\eta]\to0$, because
$H(E)=O(\log(2+E))$.  No additional smallness of $\varepsilon$ is
needed.  By \eqref{eq:nonlinear-derivative}, a stationary point
requires $t_E=0$.  Then $Q(s_E,0)=0$; strict decrease of $Q(s,0)$
and $Q(1,0)=0$ force $s_E=1$ and $E=E_*$.  The negative second
derivative at $E_*$ and absence of any other stationary point prove
the stated strict monotonicity.
\end{proof}

\begin{remark}\label{rem:nonlinear-scope}
Theorem~\ref{thm:nonlinear-conditional} is a conditional
variational principle. Uniform bounded distortion yields the
pressure upper bound, while the thermodynamic regularity
assumption (TR) provides the energy-constrained equilibrium
states and their finite-alphabet approximations needed for
the matching Hausdorff dimension lower bound.  For a
non-affine system not smoothly conjugate to a similarity IFS, the
relevant dimension is generally not $H(E)/L(E)$; its dependence on
geometry is encoded by the two-variable pressure $Q(s,t)$.  Proposition~\ref{prop:nonlinear-example} verifies (TR) for one explicit
non-conjugate full-branch model, whose roof symbols become genuine
first-return lengths in Lemma~\ref{lem:nonlinear-return-tower}.
That tower realization should not be confused with a smooth
single-interval map.  Further questions concern natural
first-return maps, including \eqref{eq:intervalmap}, and levels with
different lower and upper energy averages.
\end{remark}

\addtocontents{toc}{\protect\setcounter{tocdepth}{-2}}
\section*{Statements and Declarations}

\subsection*{Data Availability}
This article is theoretical. No research datasets were generated or analysed.

\subsection*{Competing Interests}
The author declares that he has no financial or non-financial interests
relevant to the content of this article.

\subsection*{Use of generative AI}
During the preparation of this manuscript, the author used ChatGPT
(OpenAI) to assist with language
editing, structural organization, and the checking and refinement of
mathematical exposition. The author takes full responsibility for
the mathematical arguments, calculations, references, and final wording
of the manuscript.

\end{document}